\documentclass[a4paper,9pt]{amsart}
\usepackage[english]{babel}
\usepackage[T1]{fontenc}
\usepackage{graphicx}
\usepackage{amsmath, amsfonts,amsthm, mathrsfs, amssymb}
\usepackage{verbatim}

\usepackage[OMLmathsfit]{isomath}

\usepackage{enumerate}
\usepackage[latin1]{inputenc}
\usepackage{geometry}

\usepackage{soul,xcolor,empheq}
\usepackage[many]{tcolorbox}
\tcbset{
  myhlight/.style={
    colback=yellow,
    arc=0pt,
    outer arc=0pt,
    boxrule=0pt,
    top=0pt,
    bottom=0pt,
    left=0pt,
    right=0pt,
  },
  highlight math style={myhlight}
}
\newtcolorbox{hllong}{breakable,myhlight}

\numberwithin{equation}{section}

\newcommand{\bfj}{{\bf j}}

\newcommand{\bfh}{{\bf h}}
\newcommand{\bfr}{{\bf r}}
\newcommand{\bfs}{{\bf s}}

\newcommand{\C}{{\mathbb C}}

\newcommand{\N}{{\mathbb N}}

\newcommand{\R}{{\mathbb R}}
\newcommand{\T}{{\mathbb T}}
\newcommand{\Z}{{\mathbb Z}}

\newcommand{\cB}{{\mathcal B}}

\newcommand{\cE}{{\mathcal E}}
\newcommand{\cF}{{\mathcal F}}
\newcommand{\cG}{{\mathcal G}}

\newcommand{\cI}{{\mathcal I}}

\newcommand{\cL}{{\mathcal L}}
\newcommand{\cM}{{\mathcal M}}
\newcommand{\cN}{{\mathcal N}}

\newcommand{\cP}{{\mathcal P}}

\newcommand{\cR}{{\mathcal R}}

\newcommand{\cT}{{\mathcal T}}
\newcommand{\cU}{{\mathcal U}}

\newcommand{\cW}{{\mathcal W}}

\newcommand{\cZ}{{\mathcal Z}}

\newcommand{\Ham}{\mathsf}

\newcommand{\la}{\left\langle}
\newcommand{\ra}{\right\rangle}

\newcommand{\ii}{\mathrm{i}}
\newcommand{\sleq}{\lesssim}

\newcommand{\di}{\mathrm{d}}

\newtheorem{theorem}{Theorem}[section]
\newtheorem{lemma}[theorem]{Lemma}
\newtheorem{corollary}[theorem]{Corollary}
\newtheorem{proposition}[theorem]{Proposition}
\newtheorem{definition}[theorem]{Definition}
\newtheorem{remark}[theorem]{Remark}

\usepackage[colorlinks=true]{hyperref}

\title[]{Energy cascade for the Klein-Gordon lattice}

\author{S. Pasquali$^{(\star)}$}
\address[$\star$]{Universit\'e Paris-Saclay, CNRS, Laboratoire de math\'ematiques d'Orsay, 91405, Orsay, France}
\email[$\star$]{stefano.pasquali@universite-paris-saclay.fr}

\begin{document}


\begin{abstract}
We study analytically the dynamics of a $d$-dimensional Klein-Gordon lattice with periodic boundary conditions, for $d \leq 3$. We consider initial data supported on one low-frequency Fourier mode. We show that, in the continuous approximation, the resonant normal form of the system is given by a small-dispersion nonlinear Schr\"odinger (NLS) equation. By exploiting a result about the growth of Sobolev norms for solutions of small-dispersion NLS, we are able to describe an energy cascade phenomenon for the Klein-Gordon lattice, where part of the energy is transferred to modes associated to higher frequencies. Such a phenomenon holds within the time-scale for which we can ensure the validity of the continuous approximation. \\

\end{abstract}
\maketitle
\tableofcontents

\section{Introduction} \label{intro}

In this paper we present an analytical study of the dynamics of a $d$-dimensional ($d \leq 3$) Klein-Gordon lattice with periodic boundary conditions, for initial data where only one low-frequency Fourier mode is initially excited. We give a rigorous result concerning the transfer of energy to modes associated to higher frequencies, resulting in an energy cascade phenomenon, according to the terminology used in \cite{ponno2005energy}. \\

The study of the dynamics of Hamiltonian lattices started with the numerical result by Fermi, Pasta and Ulam (FPU) \cite{fermi1995studies}, who investigated the dynamics of a one-dimensional chain of particles with nearest neighbour interaction. In the original simulations all the energy was initially given to a single low-frequency Fourier mode, with the aim of measuring the time of relaxation of the system to the `thermal equilibrium' by looking at the evolution of the Fourier spectrum. According to classical statistical mechanics (more precisely, by the theorem of equipartition of energy), the energy spectrum corresponding to the thermal equilibrium is a plateau. Fermi, Pasta and Ulam expected that the system would have reached thermal equilibrium within a short time-scale, However, the outcome of their numerical experiment was surprising: indeed, they obtained a Fourier spectrum which was far from being flat, exhibiting a lack of thermalization and a recurrent behaviour of the dynamics.

Both from a physical and a mathematical point of view, the studies on FPU-like systems have a long history: a survey of the vast literature on the subject (until 2006) is discussed in the monograph \cite{gallavotti2007fermi}. We also mention the numerical study in \cite{benettin2011time} about the time-scale for equipartition for two classes of FPU models: in this paper the authors showed that in the regime of small specific energy $\epsilon \ll 1$, the lattice first reaches a state very far from energy equipartition (the one observed in the FPU paper), and then it evolves towards equipartition of energy  among normal modes on a much longer time-scale. 

Unlike the dynamics of one-dimensional Hamiltonian lattices, the dynamics of higher-dimensional Hamiltonian lattices is far less clear; it is expected that different choices of the geometry of the lattice and of the specific energy regime could lead to different results. Moreover, when considering nearly integrable Hamiltonian lattices, the dynamics of ``a priori stable'' models (namely, when the unperturbed system is integrable, but does not present hyperbolicity) and ``a priori unstable'' models (when the unperturbed system is integrable and presents hyperbolicity) may be drastically different. 

Benettin and collaborators \cite{benettin2005time,benettin2008study} studied numerically a two-dimensional FPU lattice with triangular cells and different boundary conditions in order to estimate the equipartition time-scale. They found out that the equipartition is reached faster than in the one-dimensional case. 

Regarding the dynamics of ``a priori unstable'' models, there are some analytical results concerning Arnold diffusion for lattices of weakly-coupled penduli. In \cite{kaloshin2014arnold,huang2017energy} the authors consider a $1$-dimensional lattice of weakly-coupled penduli with periodic boundary conditions, and they prove by variational methods the existence of orbits which display an energy transfer among modes. More recently, in \cite{giuliani2023arnold} the authors construct orbits for $d$-dimensional lattices ($d \geq 1$) of weakly-coupled penduli which display Arnold diffusion. We also mention  \cite{bernard2016arnold,gidea2020general} about Arnold diffusion mechanisms for other classes of Hamiltonian systems. \\

The relation between the dynamics of Hamiltonian lattices and the dynamics of dispersive PDEs has been first highlighted in \cite{zabusky1965interaction}, by comparing the FPU model and the KdV equation (see also \cite{schneider2000counter,bambusi2006metastability,hong2021korteweg,gallone2021korteweg,gallone2022hamiltonian} and references therein, where metastability phenomena for the FPU model are described by approximating its dynamics with the dynamics of the KdV equation). We also mention \cite{bambusi2002nonlinear}, where the authors obtain the NLS equation as a resonant normal form for the FPU model. More recently, metastability phenomena have been proved for two-dimensional rectangular lattices \cite{gallone2021metastability}, by approximating their dynamics in some particular specific energy regimes with the dynamics of integrable dispersive PDEs. Regarding instability phenomena, we mention \cite{gallone2022burgers}, which describes the approach to thermalisation in the FPU model through a system of generalised Burgers equations.\\

In this paper we study a $d$-dimensional Klein-Gordon (KG) lattice with $(2N+1)^d$ sites, with polynomial nonlinearity of degree $2\ell+2$ ($\ell \geq 1$), in the case of periodic boundary conditions, for $d \leq 3$. We show the existence of an energy cascade phenomenon, in which part of the energy initially given to only one low-frequency mode is transferred to higher-frequency modes. More precisely, if we denote by $h \ll 1$ the wave-number of the Fourier mode initially excited and we consider small-amplitude solutions for which the specific energy $\epsilon \sim h^{2\alpha}$, with $\alpha$ belonging to a suitable interval $(\alpha_0,\alpha_1) \subset (0,1/\ell)$, we obtain that the lattice is approximated by a small-dispersion nonlinear Schr\"odinger (NLS) equation. We then exploit a result concerning the growth of Sobolev norms for solutions of the small-dispersion NLS equation in order to deduce a rigorous energy cascade result for the KG lattice. 

Up to the author's knowledge, this is the first rigorous result concerning instability phenomena for multidimensional ``a priori stable'' Hamiltonian lattices by continuous approximation arguments. \\

To prove our result we adopt the following strategy. The first step consists in the approximation of the dynamics of the lattice with the dynamics of a continuous system; this step gives also a natural perturbative order. Next, we perform a normal form canonical transformation and we notice that for $\frac{1}{2\ell} < \alpha < \frac{1}{\ell} $ the effective dynamics is given by a small-dispersion NLS equation. Next, we exploit a theorem by Kuksin about the growth of Sobolev norms for solutions of the small-dispersion NLS  (see Sec. 3 of \cite{kuksin1997oscillations}) in order to construct approximate solutions of the original discrete lattices, and we estimate the error with repect to a true solution with the corresponding initial datum (see Section 6 of \cite{gallone2021metastability} for a similar argument in the analytic setting). This allows us to estimate the specific energies for the modes of approximate solutions of the original lattice and to prove the energy cascade phenomenon.

The novelties of this work are: a mathematically rigorous proof of an energy cascade phenomenon for a Hamiltonian lattice through a growth of Sobolev norms the dynamics of a dispersive PDE, see Theorem \ref{dDNLSrThm} and its proof in Sec. \ref{ThmProof}; the normal form Theorem  \ref{gavthm}, which combines the techniques used in \cite{faou2011hamiltonian,faou2012geometric} with those used in \cite{bambusi2005galerkin,pasquali2018dynamics} (we mention that, unlike the results proved in \cite{bambusi2006metastability,gallone2021metastability}, the normal form theorem proved in the present paper holds for sequences in $\ell^1_s$-spaces, see also comment iv. later); we also mention the estimates in Sec. \ref{ApprSec} for bounding the error between the approximate solution and the true solution of the lattice, which need a more careful study than the ones appearing in \cite{schneider2000counter,bambusi2006metastability} for the one-dimensional case (see also Sec. 6 of \cite{gallone2021metastability}).  

Some comments are in order: 
\begin{enumerate}
\item[i.] as mentioned before, the specific energy $\epsilon$ of the system and the wave-number $h \sim \frac{1}{N}$ of the Fourier mode initially excited are linked through the relation $\epsilon \sim h^{2\alpha}$, with $\alpha$ belonging to a suitable interval $(\alpha_0,\alpha_1) \subset (0,1/\ell)$. In particular, this implies that the main result of the paper does not hold in the thermodynamic limit regime (namely, for large $N$ and for fixed small $\epsilon$ not depending on $N$), which is the most relevant regime for statistical mechanics;
\item[ii.] the averaging Theorem proved in Sec. \ref{galavsec} of the paper applies also to the endpoint case $\alpha=1/\ell$, as highlighted by \eqref{alphaBound0} and by Corollary \ref{NLScor}. It would be interesting in this context to use more recent results proved in \cite{colliander2010transfer,guardia2016growth,giuliani2022sobolev} about transfer of energy among Fourier modes for the NLS equation on the torus. However, those results hold for solutions with initial data supported on a large number of Fourier modes, and this is a first difference with respect to the initial data considered in the present paper. Moreover, some technical estimates, like the ones proved in Proposition \ref{ApprPropdNLS}, would require a more careful proof in the endpoint case;
\item[iii.] in the one-dimensional case the dynamics of the KG lattice with different boundary conditions has been studied in \cite{bambusi2009boundary}; since the regime investigated in \cite{bambusi2009boundary} corresponds to $\ell=1$ and $\epsilon \sim h^2$ in the notation of the present paper, it does not overlap with the results presented here;
\item[iv.] the mechanism of the energy cascade exhibited in this paper is obtained in four steps. First, we prove a normal form Theorem \ref{gavthm} in $\ell^1_s$-spaces, by which we can approximate the continuous approximation \eqref{KGseqc} of the Klein-Gordon lattice \eqref{Ham2KGs} with a small-dispersion NLS equation \eqref{NLSeqFT}. Next, we use a result of growth of Sobolev norms for solutions of small-dispersion NLS, see Theorem \ref{KuksinThm}. Then, we estimate the difference between the Fourier modes of the small-dispersion NLS and the normal modes of the lattice, see Proposition \ref{ApprPropdNLS}. Finally, we prove an energy transfer phenomenon in $\ell^2_m$-spaces, with $m>s$. We mention that the approach of using a normal form result in $\ell^1$-spaces in order to show growth in $\ell^2_m$-spaces has been used in recent results about growth of Sobolev norms and instability phenomena for Hamiltonian PDEs where dispersion and nonlinearity are of the same order of magnitude, see \cite{carles2013energy,hani2014long,guardia2016growth};
\item[v.] the result presented in this paper is a first example in which an instability phenomenon for a multidimensional Hamiltonian lattice is described by exploiting an instability result for a nonlinear Hamiltonian PDE. Apart from the results by Kuksin on nonlinear PDEs with small dispersion \cite{kuksin1995squeezing,kuksin1996growth,kuksin1997oscillations}, it would be interesting to apply to multidimensional Hamiltonian lattices more recent results regarding instability phenomena for nonlinear Hamiltonian PDEs, like the ones proved in \cite{gerard2010cubic,gerard2012effective,hani2015modified,gerard2018two,giuliani2021chaotic,baldi2024effective}. 
\end{enumerate}
\smallskip

The paper is organized as follows: in Section \ref{result} we introduce the mathematical setting of the model and we give an informal statement of the main result. In Section \ref{FuncSubsec} we describe the functional setting for stating in Section \ref{BNFsubsec} an abstract Averaging Theorem; we prove this theorem in Section \ref{BNFprsubsec}. In Section \ref{ApplSubsec} we apply the averaging Theorem to the KG lattice, deriving the NLS with small dispersion. In Section \ref{NLSdyn} we review a result about the dynamics of the NLS with small dispersion. In Section \ref{ApprSec} we use the normal form equations in order to construct approximate solutions, and we estimate the difference with respect to the true solutions with corresponding initial data. In Section \ref{ThmProof} we state and we prove the main result of the paper. In Appendix \ref{BNFest} we prove Lemma \ref{NFest}, while in Appendix \ref{ApprEstSec11} we prove Proposition \ref{NLSphiProp}. \\

\emph{Relevant notation.} We denote by $N$ the half length of the side of the lattice. We denote by $k=(k_1,\ldots,k_d) \in \mathbb{Z}^d$ the index of the wave-vector.

\section{Setting and informal statement of the Main Result} \label{result}

We consider a periodic $d$-dimensional lattice ($d=1,2,3$), called KG lattice, which combines the nearest-neighbour potential with an on-site one. We write
\begin{align} \label{Zd}
\Z^d_{N} &\coloneqq \{ (j_1,\ldots,j_d): j_1,\ldots,j_d \in \Z, |j_i| \leq N  \; \; \forall i=1,\ldots,d \};
\end{align}
we also write $( e_j^{(d)} )_{j=1,\ldots,d}$ for the vectors of the canonical basis of $\Z^d$. We denote by 
\begin{align} \label{Nd}
N_d &\coloneqq (2N+1)^d
\end{align}
the total numbers of sites in the lattice $\mathbb{Z}^d_N$.\\

The Hamiltonian describing the scalar KG lattice is given by
\begin{align}
H_{KG}(Q,P) &= \sum_{j \in \Z^d_{N}} \frac{P_j^2}{2} + \frac{1}{2} \sum_{j \in \Z^d_{N}} Q_j \; (-\Delta_1 Q)_j + \sum_{j \in \Z^d_{N}} U(Q_j), \label{Ham2KGs} \\
(\Delta_1 f)_j &\coloneqq \sum_{k=1}^d ( f_{j+e_k^{(d)}} - 2 f_j + f_{j-e_k^{(d)}} ), \label{Delta1Op} \\
U(x) &=  \frac{x^2}{2} + \beta \frac{x^{2\ell+2}}{2\ell+2},  \qquad \beta>0, \; \ell \geq 1, \label{potKG2}
\end{align}
The associated equations of motion are
\begin{align}
\ddot{Q_j} &= (\Delta_1 Q)_j -  Q_j - \beta Q_j^{2\ell+1}, \qquad j \in \Z^d_{N}. \label{dDKGseq} 
\end{align}
If we take $\ell=1$, we obtain a generalization of the one-dimensional $\phi^4$ model. We point out that \eqref{dDKGseq} corresponds to the discrete nonlinear Klein-Gordon equation with defocusing nonlinearity. \\

We also introduce the Fourier coefficients of $Q$ via the following standard relation,
\begin{equation} \label{fourierQ}
Q_j(t) \; \coloneqq\; \frac{1}{ \sqrt{N_d} } \sum_{k \in \Z^d_{N}} \hat{Q}_k(t) \, e^{2 \pi \ii \frac{j \cdot k}{ N_d }  } \, , \; \qquad j \in \Z^d_{N},
\end{equation}
and similarly for $P_j$. We denote by
\begin{align}
E_k &\coloneqq \frac{|\hat{P}_k|^2+\omega_k^2 |\hat{Q}_k|^2}{2}, \label{EnNormModeKG} \\
\omega_k^2 &\coloneqq  1+4 \sum_{j=1}^d \sin^2\left(\frac{k_j \, \pi}{2N+1} \right) , \label{FreqNormModeKG}
\end{align}
the energy and the square of the frequency of the mode at site $k=(k_1,\ldots,k_d) \in \Z^d_{N}$.  \\
We will restrict to solutions of \eqref{Ham2KGs} which have a discrete odd-symmetry, namely to solutions $Q=(Q_j)_{j \in \mathbb{Z}^d_N}$ such that
\begin{align*}
Q_{j_1,\ldots,j_d} &= -Q_{j_1,\ldots,-j_i,\ldots,j_d}, \; \; \forall i=1,\ldots,d, \; \; \forall j=(j_1,\ldots,j_d) \in \mathbb{Z}^N_d,
\end{align*}
and similarly for $P$ (we refer the interested reader to Sec. 7 of \cite{rink2003symmetric} for a study of invariant manifolds based on symmetries for Hamiltonian lattices and their continuous approximation). \\

As it is customary in lattices with a large number of degrees of freedom, especially in relation with statistical mechanics, we introduce the \emph{specific wave vector} $\kappa$ as
\begin{equation}\label{kappa}
\kappa \;\coloneqq\; \kappa(k) = \left( \frac{k_1}{ N+\frac{1}{2} }, \ldots, \frac{k_d}{ N+\frac{1}{2} } \right) \, , 
\end{equation} 
the \emph{specific energy of the specific normal mode $\kappa=\kappa(k)$} as
\begin{align}\label{enkappa}
\cE_\kappa &\coloneqq \frac{E_k}{ \left(N+\frac{1}{2}\right)^d } \, ,
\end{align}
and the \emph{specific energy} $\epsilon$ of the system, given by
\begin{align*}
\epsilon &\coloneqq E/N_d ,
\end{align*}
where $E$ is the total energy of the system.

The important role of the energy of normal modes can be understood by a well-known argument in classical equilibrium statistical mechanics: considering \eqref{Ham2KGs} at equilibrium and in the harmonic case $\beta=0$, by the Theorem of Equipartition of Energy the expected values of the harmonic energies $E_k$ at a given total energy $E$ do not depend on $k$ and are all equal to the so-called \emph{specific energy} $\epsilon$ of the system. The result does not change qualitatively for a slightly anharmonic system (namely, for $|\beta| \ll 1$ in \eqref{Ham2KGs}) or for small specific energy $\epsilon \ll 1$ (see pag. 155 in Chap. 4 of \cite{gallavotti2007fermi} for a more detailed explanation). \\

In the present paper we study the behaviour of small amplitude solutions of \eqref{Ham2KGs}, with initial data in which only one low-frequency Fourier mode is excited. We introduce the quantity
\begin{align} 
h &\coloneqq \frac{2}{2N+1}, \label{small} 
\end{align}
which plays the role of the small parameter in our construction: we will use it in the asymptotic expansion of the dispersion relation of the continuous approximation of the lattice in order to derive the approximating PDE in the regime we are considering. \\

We write $\kappa_0 \coloneqq  \frac{ v_0 }{N+\frac{1}{2}} = h \, v_0$, with $v_0 \coloneqq \sum_{j=1}^d e_j^{(d)}$. 

We now present an informal statement of the main result. 

\begin{theorem}[informal statement] \label{dDNLSrInf}
Consider \eqref{Ham2KGs} with $\beta > 0$ and $\ell \geq 1$.

Then for all sufficiently large $m \gg 1$ there exist an open interval $\mathcal{I} \subset (0,1/\ell)$ and $\delta_0>0$ such that the following holds true: for any $\delta \in ( \delta_0 , 1 )$ and for any $\alpha \in \mathcal{I}$ we have that for all $0 < h \ll 1$, if the initial datum satisfies
\begin{align*} 
\mathcal{E}_{\kappa_0}(0) = C_0 \; h^{2\alpha} , &\; \; \mathcal{E}_{\kappa}(0) = 0 , \; \; \forall \kappa \neq \kappa_0 ,
\end{align*}
then there exist a sufficiently large $T_m >0$ and  $c>0$ such that
\begin{align} \label{InfCascade}
\sum_{\kappa \in h \, \mathbb{Z}^d_N : |\kappa|  \leq h^{\delta} } |\kappa|^{2m} \, \cE_\kappa(T_m)  &\geq \frac{c}{8} \;  h^{2\alpha - 2(1-\ell\alpha) \mu +2m } ,
\end{align}
where $\mu=\mu(d,\ell,m)>0$.

\end{theorem}

We defer a rigorous statement of our main result, together with its proof, to Section \ref{ThmProof}. In particular, we give a rigorous stament of the main result in Theorem \ref{dDNLSrThm}. We also discuss the dependence of the quantities appearing in Theorem \ref{dDNLSrThm} with respect to parameters in Remark \ref{QuantDepRem}, and in Remark \ref{ConstantsRem} we point out where the quantities mentioned in the statement of Theorem \ref{dDNLSrThm} are defined in the paper.

\begin{remark} \label{InfRem}

First, we observe that the initial datum in Theorem \ref{dDNLSrInf} satisfies
\begin{align*}
|\kappa_0|^{2m} \, \mathcal{E}_{\kappa_0}(0) = C_0 \, |v_0|^{2m} \, h^{2\alpha+2m} , &\; \; |\kappa|^{2m} \, \mathcal{E}_{\kappa}(0) = 0 , \; \; \forall \kappa \neq \kappa_0 ,
\end{align*}
for all $m \geq 0$, so that \eqref{InfCascade} gives an analytical description of an energy transfer between lower-frequency and higher-frequency modes for the KG lattice. Such a phenomenon is called \emph{energy cascade} (see \cite{ponno2005energy}), in analogy with a similar phenomenon arising in weak turbulence theory.

We mention that the quantity
\begin{equation*}
\sum_{\kappa \in h \, \mathbb{Z}^d_N } |\kappa|^{2m} \, \cE_\kappa
\end{equation*}
is equivalent to a discrete $H^m$-Sobolev norm of the solution $(Q_j(t))_{j \in \mathbb{Z}^d_N}$, hence it is a good candidate in order to observe the evolution of the higher-frequency modes of the solution.

\end{remark}

\begin{remark} \label{InDataRem}

In Theorem \ref{InfRem} we considered an initial datum supported on the mode with wave vector $v_0$ just for simplicity. With simple modifications (see Remark 3.7 and Sec. 5 of \cite{bambusi2006metastability}), the result can be generalized to initial data supported on the mode with wave vector $v \in \mathbb{Z}^d_N \setminus \{ 0 \}$ and on its higher harmonics, under the assumption that $\left| \frac{v}{N+\frac{1}{2}} \right| \ll 1$ .
\end{remark}

As mentioned in comment $\mathrm{i.}$ of the introduction, the regimes of specific energy $\epsilon \sim h^{2\alpha}$ described in Theorem \ref{dDNLSrInf} are above the regime for which one can prove the existence of metastability phenomena (see Theorem 2.7 in \cite{gallone2021metastability}), and they are also below the so-called thermodynamic limit (since $\alpha >0$). 

We point out that in \eqref{InfCascade} the sum on the left-hand side only involves the wave vectors $k$ such that $|\kappa|  \leq h^{\delta}$, and that the index $m$ is very large; these two facts imply that the energy transfer exhibited in \eqref{InfCascade} is weaker compared to the energy equipartition observed numerically both for the KG lattice \cite{pistone2019universal} and for the FPU lattice \cite{benettin2011time,benettin2008study}.

\section{Galerkin Averaging} \label{galavsec}

In this Section we prove an abstract averaging theorem, which we will later apply in order to show that for large $N$ the dynamics of the KG lattice can be approximated by the dynamics of a NLS equation with small dispersion.  

The idea of its proof (following \cite{bambusi2005galerkin,bambusi2006metastability,faou2011hamiltonian,pasquali2018dynamics,gallone2021metastability}) is to make a Galerkin cut-off, namely to approximate the original infinite dimensional system by a finite dimensional one, to put in normal form the truncated system, and then to choose the dimension of the truncated system in such a way that the error due to the Galerkin cut-off and the error due to the truncation in the normalization procedure are of the same order of magnitude. The system one gets is composed of a part which is in normal form, and of a remainder.

If we neglect the remainder, we obtain a system whose solutions are approximate solutions of the original system; in Sec. \ref{ApprSec} we will show how to control the error with respect to a true solution of the original system.

This Section is divided into four parts. In the first part we introduce the analytic setting we are working with. In the second part we state the averaging Theorem \ref{gavthm}. In the third part we give a concise proof of the averaging Theorem, deferring the proof of the technical Lemma \ref{NFest} to Appendix \ref{BNFest}. In the last part we show an application of the abstract result to study the dynamics of the Klein-Gordon lattice.

\subsection{Functional Setting} \label{FuncSubsec}

\begin{definition} \label{defsigk}
Fix a constant $s \geq 0$, and let $\cN = \Z^d$ or $\N^d$, and $i=1,2$.
We will denote by $\ell^i_{s}$ the Banach space of complex sequences $v = (v_n)_{n \in \cN }$ with obvious vector space structure and with norm given by
\begin{align} \label{normsigk}
\|v\|_{\ell^i_s} &\coloneqq \left( \sum_{n \in \cN } |n|^{is} \, |v_n|^i  \right)^{1/i} , \; \; |n|^2 \coloneqq \max(1,n_1^2+\cdots+n_d^2), \; \; \forall n=(n_1,\cdots,n_d).
\end{align}
We will denote by $\ell^i$ the space $\ell^i_{0}$.
\end{definition}

We define the set $\cZ \coloneqq \cN \times \{\pm 1\}$; for $\mathtt{j}=(a,\sigma) \in \cZ$, we define $| \mathtt{j} | \coloneqq |a|$ and we denote $\bar{ \mathtt{j} }\coloneqq(a,-\sigma)$. We identify $(\xi,\eta) \in \C^{\cN} \times \C^{\cN}$ with $( \zeta_{\mathtt{j}} )_{ \mathtt{j} \in \cZ} \in \C^{\cZ}$ through the formula
\begin{equation*}
 \mathtt{j} = (a,\sigma) \in \cZ \; \Longrightarrow \; 
\begin{cases}
\zeta_{\mathtt{j}} = \xi_a & \text{if} \; \; \;  \sigma=1, \\
\zeta_{\mathtt{j}} = \eta_a & \text{if} \; \; \; \sigma=-1.
\end{cases}
\end{equation*}
Now fix $s \geq 1$, and consider the scale of Banach spaces 
$\cW^{s}\coloneqq \ell^1_{s} \times \ell^1_{s} \ni \zeta=(\xi,\eta)$, 
endowed with the symplectic form
\begin{align}
\Omega_1 &\coloneqq -\mathrm{i} \sum_{a \in \cN} \mathrm{d}\xi_a \wedge \mathrm{d}\eta_a ,\label{eq:SymplecticForm}
\end{align}
and the norm
\begin{align*}
\|\zeta\|_{\cW^s} &\coloneqq \|\xi\|_{\ell^1_s} + \|\eta\|_{\ell^1_s} , \; \; \forall \zeta=(\xi,\eta) \in \cW^s .
\end{align*}

If we fix $s$ and $\cU_{s} \subset \cW^{s}$ open, the Hamiltonian vector field of the Hamiltonian function $\Ham H \in C^1(\cU_{s},\C)$ is given by 
\begin{align} \label{HamVF}
X_{\Ham H}(\zeta) &\coloneqq \Omega^{-1}_{1} \nabla_{\zeta} \Ham H(\zeta),
\end{align}
where
\begin{align*}
\nabla_{\zeta} \Ham H(\zeta) &\coloneqq \left( \frac{\partial \Ham{H} }{\partial \zeta_{\mathtt{j}} } \right)_{ \mathtt{j} \in \cZ}, \\
\end{align*}
and for $\mathtt{j} =(a,\sigma) \in \cZ$
\begin{equation*}
\frac{\partial \Ham{H} }{\partial \zeta_{\mathtt{j}} } \; \coloneqq \; 
\begin{cases}
\frac{ \partial \Ham{H} }{ \partial \xi_a } & \text{if} \; \; \; \sigma=1, \\
\frac{ \partial \Ham{H} }{ \partial \eta_a } & \text{if} \; \; \; \sigma=-1.
\end{cases}
\end{equation*}
If we have $\Ham{F},\Ham{G} \in C^1(\cU_{s},\C)$ the Poisson bracket $\{ \Ham{F},\Ham{G} \}$ is defined by
\begin{align*}
\{ \Ham{F},\Ham{G} \} &\coloneqq \nabla_{\zeta}\Ham{F}^T \; \Omega^{-1}_1 \nabla_{\zeta}\Ham{G} \; = \; \mathrm{i} \; \sum_{a \in \cN} \frac{ \partial \Ham{F} }{ \partial \xi_a } \frac{ \partial \Ham{G} }{ \partial \eta_a } - \frac{ \partial \Ham{F} }{ \partial \eta_a } \frac{ \partial \Ham{G} }{ \partial \xi_a }  .
\end{align*}
We denote the open ball of radius $R$ and center $0$ in $\ell^1_{s}$ by $B_{s}(R)$; we write $\cB_{s}(R)\coloneqq B_{s}(R) \times B_{s}(R) \subset \cW^{s}$. 

Now, we introduce the Fourier projection operators $\hat{\pi}_j: \ell^1_{s} \to \ell^1_{s}$
\begin{align*} 
\hat{\pi}_j((v_n)_{ n\in\Z^d \setminus\{0\} }) &\coloneqq 
\begin{cases}
v_n & \text{if} \qquad j-1 \leq |n| < j \\
0 & \text{otherwise}
\end{cases} \qquad , \qquad j \geq 1,
\end{align*}
the operators $\pi_j: \cW^{s} \to \cW^{s}$
\begin{align} \label{smallpi}
\pi_j((\zeta_n)_{ n\in\Z^d \setminus\{0\} }) &\coloneqq 
\begin{cases}
\zeta_n & \text{if} \qquad j-1 \leq |n| < j \\
0 & \text{otherwise}
\end{cases} \qquad , \qquad j \geq 1,
\end{align}
and the operators $\Pi_M: \cW^{s} \to \cW^{s}$
\begin{align} \label{bigpi}
\Pi_M((\zeta_n)_{ n\in\Z^d \setminus\{0\} })  &\coloneqq 
\begin{cases}
\zeta_n & \text{if} \qquad |n| \leq M \\
0 & \text{otherwise}
\end{cases} \qquad , \qquad M \geq 0.
\end{align}
Last, we define the operator $\overline{\pi_j}\coloneqq\mathrm{id}-\pi_j$. The projection operators defined in \eqref{smallpi}-\eqref{bigpi} satisfy:
\begin{enumerate}
\item[$i.$] 
\begin{align*}
\zeta &= \sum_{j \geq 0} \pi_j\zeta , \;\; \forall \zeta \in \mathcal{W}^{s};
\end{align*}
\item[$ii.$] for any $M \geq 0$ 
\begin{align*}
\|\Pi_M\zeta\|_{\cW^{s}}&\leq  \, \|\zeta\|_{\cW^{s}} , \;\; \forall\zeta \in \mathcal{W}^{s}.
\end{align*}
\end{enumerate}

We defer the reader to Sec. III.1-III.2 of \cite{faou2012geometric} for more details.

In particular, we mention the following result about $\ell^p$ spaces (see Proposition III.2 in \cite{faou2012geometric}).

\begin{lemma} \label{lpRem}
Let $m,s \in \mathbb{R}$ be such that $m > s + d/2$. Then we have 
\begin{equation*}
\ell^2_{m} \subset \ell^1_s \subset \ell^2_s
\end{equation*}
and there exists $C=C(s,m)>0$ such that 
\begin{equation} \label{lpNormsIneq}
\| z \|_{\ell^2_s} \; \leq \; \| z \|_{\ell^1_s} \; \leq \; C \; \| z \|_{\ell^2_{m}}, \; \; \forall z \in \ell^2_{m} .
\end{equation}

\end{lemma}

We now describe the class of nonlinearities we are going to consider.
Let $l \geq 2$, we consider $\bfj =( \mathtt{j}_1,\ldots, \mathtt{j}_l) \in \cZ^{l}$, and where $\mathtt{j}_i=(a_i,\sigma_i)$, with $a_i \in \cN$ and $\sigma_i \in \{\pm 1\}$, for all $i \in \{1,\ldots,l\}$. We define
\begin{align*}
\overline{\bfj} \coloneqq (\overline{ \mathtt{j}_1},\ldots,\overline{ \mathtt{j}_{l}}), &\; \; \overline{ \mathtt{j}_i} \coloneqq (a_i,-\sigma_i), \; \; i = 1,\ldots,l.
\end{align*}
We also use the notation $\zeta_\bfj \coloneqq \zeta_{ \mathtt{j}_1} \cdots \zeta_{ \mathtt{j}_l}$, and we define the momentum $\cM(\bfj)$ of a multi-index $\bfj$ and the set $\cI_l$ of indices with zero momentum as following,
\begin{align*}
\cM(\bfj) &\coloneqq \sum_{i=1}^{l} a_i\sigma_i, \;\; \cI_{l} \coloneqq \{ \bfj =(\mathtt{j}_1,\ldots, \mathtt{j}_{l}) \in \cZ^{l} | \cM(\bfj)=0 \}. 
\end{align*}

\begin{definition} \label{defPkClass}

We say that a polynomial $P$ of degree $k$ belongs to the class $\cP_k$ if $P$ is real, it has a zero of order at least $2$ in $\zeta=0$, and if $P$ contains only monomials of zero momentum, namely it is of the form
\begin{align*}
P(\zeta) &= \sum_{l=2}^k \sum_{\bfj \in \cI_{l}} a_{\bfj} \zeta_{\bfj},
\end{align*}
where $a_{\overline{\bfj}} = \overline{a_{\bfj}}$, and the coefficients are bounded,
\begin{align*}
\forall \ell=2,\ldots,k, &\; \; \forall \bfj =( \mathtt{j}_1,\ldots, \mathtt{j}_l) \in \cI_{l}, \; \; |a_{\bfj}| \leq C.
\end{align*}
\end{definition}

As an example, polynomials belonging to the class $\cP_k$ are the ones corresponding in Fourier variables to nonlinearities of the form $\int_{\T^d} p(\psi(y),\phi(y)) \di y$, where $p \in \C[X,Y]$ is a polynomial of degree $k$ having a zero of order at least $2$ at the origin and such that $p(\psi,\overline{\psi}) \in \R$.

We define the following norm
\begin{align}
\| P \| &= \sum_{\ell=2}^k \; \sup_{\bfj \in \cI_{\ell}} |a_{\bfj}| , \; \; \forall P \in \cP_k . \label{defSupNorm}
\end{align}

We mention the following result (see Proposition III.6 of \cite{faou2012geometric} for more details; see also Proposition 2.8 of \cite{faou2011hamiltonian}), which exploits the property of zero momentum.

\begin{proposition} \label{PropFG}

Let $k \geq 2$ and $s \geq 0$, then $\cP_k \subset \cW_s$, and for $P \in \cP_k$ we have
\begin{align}
\| X_P(\zeta) \|_{\cW^s} &\leq 2k(k-1)^s \; \|P\| \; \|\zeta\|_{\cW^s} \; \max_{n=1,\ldots,k-2} \|\zeta\|_{\cW^s}^n, \; \; \forall \zeta \in \cW^s. \label{EstFG2} 
\end{align}
Moreover, if $P \in \cP_k$ and $Q \in \cP_h$, then $\{P,Q\} \in \cP_{k+h-2}$, and
\begin{align}
\| \{P,Q\} \| &\leq 2kh \, \|P\| \, \|Q\|.  \label{EstFG4}
\end{align}
\end{proposition}

\subsection{An Averaging Theorem} \label{BNFsubsec}

Now we let $\mathcal{N} = \mathbb{Z}^d$, and we consider a Hamiltonian system of the form
\begin{align} \label{Hamdecomp}
\Ham H &= \Ham h_0 + \delta \Ham F, 
\end{align}
where we assume that
\begin{itemize}
\item[(PER)]  $\Ham{h}_0$ generates a linear periodic flow $\Phi^\tau_{\Ham{h}_0}$ with period $T$, 
\begin{align*}
\Phi^{\tau+T}_{\Ham{h}_0} &= \Phi^\tau_{\Ham{h}_0} \qquad \forall \tau,
\end{align*}
which is analytic as a map from $\cW^{s}$ into itself for any $s \geq 1$. Furthermore, the flow is an isometry for any $s \geq 1 $. 

\item[(INV)] for any $s \geq 1$, 
$\Phi^\tau_{\Ham{h}_0}$ leaves invariant the space $\Pi_j\cW^{s}$ for any $j\geq 0$. 
Furthermore, for any $j \geq 0$ 
\[ \pi_j \circ \Phi^\tau_{\Ham{h}_0} = \Phi^\tau_{\Ham{h}_0} \circ \pi_j. \]
\end{itemize}

Next, we assume that there exist $\nu \in \mathbb{R}$ and $k \in \mathbb{N}$, where 
\begin{equation} \label{nuAss}
- \frac{1}{2} < \nu \leq 0
\end{equation}
and $k \geq 3$, such that the vector field of $\Ham F \in \cP_k$ admits an asymptotic expansion in $\delta$ of the form
\begin{align}
\Ham F &\sim \sum_{j \geq 1} \delta^{j-1} (\Ham F_{1j} + \delta^{\nu} \Ham F_{2j}) , \label{expF} \\
X_{\Ham F} &\sim \sum_{j \geq 1}  \delta^{j-1} (X_{\Ham F_{1j}} + \delta^{\nu} X_{\Ham F_{2j}}) , \label{expXF}
\end{align}
where we used the notation \eqref{HamVF} in \eqref{expXF}, and that the following property is satisfied
\begin{itemize}

\item[(HVF)] There exist $R^\ast>0$ and a non-increasing sequence of  integers $(\mathfrak{h}_j)_{j \geq 1}$ with $\mathfrak{h}_1=k$ such that for any $j \geq 1$ 
\begin{itemize}
\item[$\cdot$] $X_{\Ham F_{1j}} \in \cP_{1}$ is analytic from $\cB_{s+2j}(R^\ast)$ to $\cW^{s}$, and $X_{\Ham F_{2j}} \in \cP_{\mathfrak{h}_j-1}$ is analytic from $\cB_{s}(R^\ast)$ to $\cW^{s}$.
\end{itemize}
Moreover, there exists a non-increasing sequence of integers $(\mathfrak{j}_r)_{r \geq 1}$ with $\mathfrak{j}_r \leq k$ for all $r \geq 1$ such that for any $r \geq 1$ 
\begin{itemize}
\item[$\cdot$] $X_{\Ham F - \sum_{j=1}^r \delta^{j-1} ( \Ham F_{1j} + \delta^\nu \Ham F_{2j} ) } \in \cP_{k-1}$ is analytic from $\cB_{s+2(r+1)}(R^\ast)$ to $\cW^{s}$.
\end{itemize}

\end{itemize}

The main result of this Section is the following theorem.

\begin{theorem} \label{gavthm}
Fix $R>0$, $s_1\gg 1$. 
Consider \eqref{Hamdecomp}, and assume (PER), (INV) and (HVF), with $\nu$ as in \eqref{nuAss}.
Then $\exists$ $s_0>0$ with the following properties: 
for any $s \geq s_1$ there exists $\delta_{s} \ll 1$ such that 
for any $\delta<\delta_{s}$ there exists an analytic canonical transformation 
$\cT_\delta:\cB_{s}(R) \to \cB_{s}(R)$  such that
\begin{align} \label{TransfHam}
\Ham H_1 \coloneqq \Ham H \circ \cT_\delta &= \Ham{h}_0 + \delta \cZ_{11} + \delta^{1+\nu} \cZ_{21}  + \delta^{2(1+\nu)} \; \mathcal{R}^{(1)},
\end{align}
where $\cZ_{11} \in \cP_{2}$ and $\cZ_{21} \in \cP_{k}$ are in normal form, namely
\begin{align} \label{NFthm}
\{\cZ_{i1},\Ham{h}_0\} &= 0, \; \; i=1,2,
\end{align} 
and there exists a positive constant $C'_{s}$ (that depends on $s$) such that
\begin{align*} 
\sup_{\cB_{s+s_0}(R)} \|X_{\cZ_{i1}}\|_{\cW^{s}} &\leq C'_{s}, \; \; i=1,2,
\end{align*} 
\begin{align} \label{Remthm}
\sup_{\cB_{s+s_0}(R)} \|X_{\mathcal{R}^{(1)}}\|_{\cW^{s}} &\leq C'_{s},
\end{align} 
\begin{align} \label{CTthm}
\sup_{\cB_{s}(R)} \|\cT_\delta-\mathrm{id}\|_{\cW^{s}} &\leq C'_{s} \, \delta^{1+\nu}.
\end{align}
In particular, 
\begin{align} \label{average}
\cZ_{i1}(\zeta) = \la \Ham F_{i1} \ra(\zeta), \; \; i=1,2,
\end{align}
where $\la \Ham F_{i1} \ra(\zeta) \coloneqq \int_0^{T} \Ham F_{i1}\circ\Phi^\tau_{\Ham{h}_0}(\zeta) \frac{\di\tau}{T}$.
\end{theorem}

\begin{remark} \label{GalerkinRem}

We observe that the terms $\mathcal{Z}_{11}$ and $\mathcal{Z}_{21}$ appearing in Theorem \ref{gavthm} come from the terms $\Ham F_{11}$ and $\Ham F_{21}$, respectively. 

Moreover, we point out that the operator $\delta \Ham F - \delta (\Ham F_{11}+\delta^{1+\nu} \Ham F_{21})$ has order $\mathcal{O}( \delta^{2+\nu} )$, whereas the remainder in \eqref{TransfHam} has order $\mathcal{O}( \delta^{2(1+\nu)} )$, so it is larger (in terms of order of magnitude in $\delta$) for $\nu <0$; this loss comes from the normal form procedure performed in the proof of Lemma \ref{itlemma}. The assumption $\nu > - \frac{1}{2}$ ensures that the remainder of the normal form procedure is smaller than the terms which are in normal form with the unperturbed Hamiltonian $\Ham h_0$, so that it can be regarded as a lower order correction. 

On the other hand, in Sec. \ref{ApplSubsec} and Sec. \ref{NLSdyn} we will exploit the fact that the Hamiltonian in normal form that we obtain has a very simple structure, and is associated to a semilinear Hamiltonian PDE whose solutions, under suitable assumptions, exhibit growth of their Sobolev norms.

\end{remark}

\subsection{Proof of the Averaging Theorem} \label{BNFprsubsec}

The proof of Theorem \ref{gavthm} is actually an application of the techniques used in \cite{pasquali2018dynamics,bambusi2006metastability}.

First notice that by assumption (INV) the Hamiltonian vector field of $\Ham{h}_0$ 
generates a continuous flow $\Phi^\tau_{{\Ham h_0}}$ which leaves $\Pi_M\cW^{s}$ invariant. \\

Now we set $\Ham H = \Ham H_{1,M} + \cR_{1,M} + \cR_1$, where \\
\begin{align} \label{truncsys}
\Ham H_{1,M} &\coloneqq \Ham h_{0} + \delta \, \Ham F_{1,M}, \\ 
\Ham F_{1,M} &\coloneqq (\Ham F_{11}+\delta^\nu \, \Ham F_{21}) \circ \Pi_M,
\end{align}
and
\begin{align} \label{remsys}
\cR_{1,M} &\coloneqq \Ham{h}_0 + \delta (\Ham F_{11} + \delta^\nu \Ham F_{21}) - \Ham H_{1,M}, \\
\cR_1 &\coloneqq \delta \left( \Ham F - \Ham F_{11} - \delta^\nu \Ham F_{21} \right).
\end{align}

The system described by the Hamiltonian \eqref{truncsys} is the one that 
we will put in normal form.  \\
In the following we will use the notation $c \sleq b$ to mean: 
there exists a positive constant $K$ independent of $M$ and $R$ 
(but eventually on $s$), such that $c \leq Kb$. 
We exploit the following intermediate results:

\begin{lemma} \label{truncest}
For any $s \geq s_1$ there exists $R>0$ such that $\forall \, \sigma >0$, $M>0$ 
\begin{align} \label{truncremt}
\sup_{\zeta \in  \cB_{s +\sigma+2 }(R) } \|X_{\cR_{1,M}}(\zeta)\|_{\cW^{s}} &\sleq \; \frac{ \delta^{1+\nu} }{M^{ \sigma }}, 
\end{align}
\begin{align} \label{expremest} 
\sup_{\zeta \in  \cB_{s+4}(R) } \|X_{\cR_1}(\zeta)\|_{\cW^{s}} &\sleq \delta^{ {2+\nu} }.
\end{align}
\end{lemma}

\begin{proof}
We set $\Ham F_{i1,M} \coloneqq \Ham F_{i1} \circ \Pi_{M}$, for $i=1,2$.

We recall that $\cR_{1,M} = \Ham{h}_0 + \delta \Ham F_1 - \Ham H_{1,M} =\delta (F_{11}-F_{11,M}) + \delta^{1+\nu}(F_{21}-F_{21,M}) $. \\
We first notice that  $\|\mathrm{id}-\Pi_M\|_{ \cW^{s+\sigma} \to \cW^{s} } \leq M^{-\sigma}$: indeed, let $f=(f_1,f_2) \in \cW^{s+\sigma}$, then
\begin{align*}
\left\| \sum_{j \geq M+1} \pi_j f \right\|_{\cW^{s}} &= \sum_{k \in \mathbb{Z}^d} |k|^s \left[ \left( \sum_{j \geq M+1} \pi_j f_1 \right)_k  +\left( \sum_{j \geq M+1} \pi_j f_2 \right)_k \right] \\
&= \sum_{k \in \mathbb{Z}^d} |k|^s \left[  \sum_{j \geq M+1} \pi_j f_{1k}   + \sum_{j \geq M+1} \pi_j f_{2k}  \right] \\
&\stackrel{ f \in \cW^s }{=} \sum_{j \geq M+1}  \; \sum_{k \in \mathbb{Z}^d} |k|^s ( \pi_j f_{1k}   + \pi_j f_{2k} ) \\
&= \sum_{j \geq M+1}  \; \sum_{j-1 \leq |k| < j} |k|^s (  f_{1k}   +  f_{2k} ) \\
&\leq \sum_{j \geq M+1} |j-1|^{-\sigma} \sum_{j-1 \leq |k| < j} |k|^{s+\sigma} (  f_{1k}   +  f_{2k} ) \\
&\leq M^{-\sigma} \|f\|_{\cW^{s+\sigma}} .
\end{align*}

Inequality \eqref{truncremt} follows from the following estimate,
\begin{align*}
&\sup_{\zeta \in \cB_{s+2+\sigma}(R)} \; \| X_{\cR_{1,M}}(\zeta)\|_{\cW^{s}} \\
& \sleq \; \|\mathrm{d} X_{ \delta \Ham F_{11} + \delta^{1+\nu} \Ham F_{21} }\|_{ L^\infty(\cB_{s+2}(R),\cW^{s}) } \|\mathrm{id}-\Pi_M\|_{ L^\infty(\cB_{s+2+\sigma}(R),\cB_{s+2}(R) ) } \\
&\sleq \delta^{1+\nu} \, M^{-\sigma},
\end{align*}
where we denoted by $\mathrm{d} X_{ \delta \Ham F_{11} + \delta^{1+\nu} \Ham F_{21} }$ the differential of $X_{ \delta \Ham F_{11} + \delta^{1+\nu} \Ham F_{21} }$; the fact that the differential acts from $\cB_{s+2}(R)$ to $\cW^{s}$ follows from (HVF). Estimate \eqref{expremest} is an immediate consequence of (HVF).
\end{proof}

\begin{lemma} \label{pertestlemma}
For any $s \geq s_1$ 
\begin{align*}
\sup_{\zeta \in  \cB_{s}(R^*) } \|X_{\Ham F_{1,M}}(\zeta)\|_{\cW^{s}} &\leq ( K^{(\Ham F_{11} )}_{1,s} + \delta^\nu K^{(\Ham F_{21} )}_{1,s} ) M^{2} , 
\end{align*}
where 
\begin{align*}
K^{( \Ham F_{11} )}_{1,s} \coloneqq \sup_{\zeta \in \cB_{s}(R^*) } \|X_{\Ham F_{11}}(\zeta)\|_{\cW^{s-2}} \; < \; + \infty, &\;\;
K^{( \Ham F_{21} )}_{1,s} \coloneqq \sup_{\zeta \in \cB_{s}(R^*) } \|X_{\Ham F_{21}}(\zeta)\|_{\cW^{s}} \; < \; + \infty.
\end{align*}
\end{lemma}

\begin{proof}
We have
\begin{align*}
& \sup_{\zeta \in \cB_{s}(R)} \left\| \sum_{j \leq M} \pi_j X_{\Ham F_{1,M}} (\zeta) \right\|_{\cW^{s}} \\
&\leq M^{2} \sup_{\zeta \in \cB_{s}(R)} \|X_{ \Ham F_{11} + \delta^\nu \Ham F_{21} }(\zeta)\|_{\cW^{s-2}} \leq ( K^{(\Ham F_{11} )}_{1,s} + \delta^\nu K^{(\Ham F_{21} )}_{1,s} ) M^{2},
\end{align*}
where the last quantity is finite for $R \leq R^\ast$ by property (HVF).
\end{proof}

To normalize \eqref{truncsys} we need to prove a reformulation of Theorem 4.4 in \cite{bambusi1999nekhoroshev}. Here we report a statement of the result adapted to our context which is proved in Appendix  \ref{BNFest}.

\begin{lemma} \label{NFest}
Let $s \geq s_1+2$, $R>0$, and consider the 
system \eqref{truncsys}. 

Then there exists $\delta_0=\delta_0(T,\|\Ham F_{11}\|,\|\Ham F_{21} \|,R,k,\nu)>0$ such that if $\delta \leq \delta_0$, then there exists an analytic canonical transformation 
$\cT^{(0)}_{\delta,M}: \cB_{s}(R) \to \cB_{s}(R)$ such that
\begin{align}
&\sup_{\zeta \in \cB_{s}(R)} \|\cT^{(0)}_{\delta,M}(\zeta)-\zeta\|_{\cW^{s}} \nonumber \\
&\leq 4T \, k (k-1)^s \, M^2 \left( \delta^{1+\nu}  \|\Ham F_{21}\| \, (1+\delta M^2 \|\Ham F_{11}\|) + \delta \|\Ham F_{11}\| \right) \, R \;\sup_{\zeta \in \cB_{s}(R)} \max_{n=1,\ldots,k-2} \|\zeta\|_{\cW^s}^n , \label{CTlemma}
\end{align}
and that puts \eqref{truncsys} in normal form up to a small remainder, 
\begin{align} \label{stepr}
\Ham H_{1,M} \circ \cT^{(0)}_{\delta,M} &= \Ham h_{0} + \delta \Ham Z^{(1)}_M + \delta^{1+\nu} \Ham Z^{(2)}_M + \delta^{ 2(1+\nu) } \cR^{(1)}_M, 
\end{align}
with $Z^{(i)}_M$ in normal form, namely $\{h_{0,M},Z^{(i)}_M\}=0$, for $i=1,2$.  Moreover, we have
\begin{align}
&\sup_{\zeta \in \cB_{s}(R)} \|X_{\Ham Z_M^{(i)}}(\zeta)\|_{\cW^{s}} \nonumber \\
&\leq 2k (k-1)^s \; \| \Ham F_{i1} \| M^2 \; R \; \sup_{\zeta \in \cB_{s}(R)} \max_{n=1,\ldots,k-2} \|\zeta\|_{\cW^s}^n,  \; \; \forall \zeta \in \cW^s, \; \; i=1,2, \label{itervfj1} 
\end{align}
and 
\begin{align}
&\sup_{\zeta \in \cB_{s}(R)} \Vert X_{ \mathcal{R}_M^{(1)}  } (\zeta) \Vert_{\cW^s} \nonumber \\
&\leq 2^4 \, T \, M^4 \, ( k_1 -1)(2 k_1 -3)^s \, R \, \sup_{\zeta \in \cB_s(R)} \max_{n=1,\ldots,2(k_1-2)} \|\zeta\|_{\cW^s}^n \, \times \nonumber \\
&\;\; \times  \bigg\{ 4 \, \| \Ham F_{11} \| \, \left[ \left( \left( 2T k_1(k_1-1)^s M^2 \, \| \Ham F_{11} \| \, R \, \sup_{\zeta \in \cB_{s}(R)}  \max_{n=1,\ldots,k_1-2} \|\zeta\|_{\cW^s}^n \right) + 4 \right) \, \| \Ham F_{11} \| +  \| \Ham F_{21} \|  \right] \nonumber \\
&\;\;\;\;  + \, \| \Ham F_{21} \|^2 \left[ \left( 4T k_1(k_1-1)^s \, M^2 \, \| \Ham F_{21} \| \, R \, \sup_{\zeta \in \cB_{s}(R)}  \max_{n=1,\ldots,k_1-2} \|\zeta\|_{\cW^s}^n \right) + 6 \right] + 2 \, \| \Ham F_{11} \|^2   \bigg\} \nonumber \\
&\;\;\;\; + 2k_1(k_1-1)^s \, ( \rho_0+ 4T \, M^4 \, \| \Ham F_{11} \| \, \| \Ham F_{21} \| )  \, R \, \sup_{\zeta \in \cB_{s}(R)} \max_{n=1,\ldots,k_1-2} \|\zeta\|_{\cW^s}^n \label{vecfrem0}
\end{align}
for $\nu=0$, and
\begin{align}
&\sup_{\zeta \in \cB_{s}(R)} \Vert X_{ \mathcal{R}_M^{(1)}  } (\zeta) \Vert_{\cW^s} \nonumber \\
&\leq \left[ \left( 4T k_1(k_1-1)^s M^2 \, \| \Ham F_{21} \| \, R \, \sup_{\zeta \in \cB_{s}(R)}  \max_{n=1,\ldots,k_1-2} \|\zeta\|_{\cW^s}^n \right) + 8 \right] \times \nonumber \\
&\; \; \; \; \; \; \; \; \times \,  2^4 (k_1 -1)(2 k_1 -3)^s \; T M^4 \, \| \Ham F_{21} \|^2 \,  R \, \sup_{\zeta \in \cB_s(R)} \max_{n=1,\ldots,2(k_1-2)} \|\zeta\|_{\cW^s}^n  \label{vecfrem}
\end{align}
for $- \frac{1}{2} < \nu < 0$. 

\end{lemma}

Now we conclude with the proof of  Theorem \ref{gavthm}. 

\begin{proof}[Proof of  Theorem \ref{gavthm}]
If we define $\delta_s \coloneqq \delta_0$, where $\delta_0$ is the quantity appearing in the statement of Lemma \ref{NFest}, and we choose
\begin{align*}
s_0 &= \sigma+2, \\
\sigma &\geq 2,
\end{align*}
then the transformation $\cT_\delta\coloneqq\cT^{(0)}_{\delta,M}$ defined by 
Lemma \ref{NFest} satisfies \eqref{TransfHam} because of \eqref{stepr}. 

Next, Eq. \eqref{NFthm} follows from Lemma \ref{NFest}, Eq. \eqref{Remthm} follows from \eqref{vecfrem}, while \eqref{CTthm} is precisely \eqref{CTlemma}.
Finally, \eqref{average} can be deduced by applying Lemma \ref{homeqlemma} to $\Ham G=\Ham F_{i1}$, for $i=1,2$.
\end{proof}

\subsection{Application to the KG lattice} \label{ApplSubsec}

We want to study small amplitude solutions of \eqref{dDKGseq}, with initial data in which only one low-frequency Fourier mode is excited.

The first step is to introduce an interpolating function $\mathsfit{Q}=\mathsfit{Q}(t,x)$ such that $ $
\begin{itemize}
\item[(B1)] $\mathsfit{Q}(t,j)=Q_j(t)$, for all $j \in \Z^d_{N}$;
\item[(B2)] $\mathsfit{Q}$ is periodic with period $2N+1$ and odd in the $x_k$-variable, for all $k=1,\ldots,d$;
\item[(B3)] $\mathsfit{Q}$ fulfills 
\begin{align}
\ddot{\mathsfit{Q}} &= -\mathsfit{Q} + \Delta_1 \mathsfit{Q} - \beta \mathsfit{Q}^{2\ell+1}, \label{KGseqc}
\end{align}
where $\Delta_1$ is the operator defined by
\begin{align}  \label{Delta1c}
\Delta_1 &\coloneqq 4 \sum_{k=1}^d \sinh^2\left(\frac{\partial_{x_k}}{2}\right) .
\end{align}
\end{itemize}
One can check that the operator \eqref{Delta1c} is equivalent to \eqref{Delta1Op} by using functional calculus. We refer to Sec. 7 of \cite{rink2003symmetric} for a study of invariant manifolds based on symmetries for Hamiltonian lattices and their continuous approximation.

It is easy to verify that \eqref{KGseqc} is Hamiltonian with Hamiltonian function
\begin{align} \label{HamKGsc}
\Ham H_{KG}(\mathsfit{Q},\mathsfit{P}) &= \int_{ [-\frac{1}{h},\frac{1}{h}]^d } \left( \frac{\mathsfit{P}^2}{2} +  \frac{\mathsfit{Q}^2}{2} - \frac{\mathsfit{Q} \; \Delta_1 \mathsfit{Q}}{2} + \beta \frac{\mathsfit{Q}^{2\ell+2}}{2\ell+2}\right) \, \di x,
\end{align}
where $\mathsfit{P}$ is a periodic function and is canonically conjugated to $\mathsfit{Q}$.

Starting from the Hamiltonian \eqref{HamKGsc}, we look for small amplitude solutions of the form
\begin{align} 
\mathsfit{Q}(t,x) &= h^{\alpha} \, q( t, h x), \; \; 0 < \alpha \leq \frac{1}{\ell}, \label{NLSr}
\end{align} 
where $q : \mathbb{R} \times [-1,1]^d \to \mathbb{R}$ is a space-periodic function and $h$ is defined in \eqref{small}.

\begin{remark} \label{RegimeRem}

The ansatz \eqref{NLSr} introduces a relation between the number of particles of the lattices (through the quantity $h$) and the specific energy, which we denote by $\epsilon$. Indeed, we have
\begin{align*}
\epsilon &\sim \frac{1}{(2N+1)^d} \sum_{k \in \Z_N^d} h^{2\alpha} q^2,
\end{align*}
which implies
\begin{align}
\epsilon &\sim h^{2\alpha}. \label{RegimeRel}
\end{align}
\end{remark}

We introduce the rescaled space variable $y = h x$, and we define $I$ as 
\begin{align} 
I &\coloneqq [-1,1]^d. \label{I}
\end{align}
Recalling \eqref{KGseqc}, the rescaled equation of motion is given by
\begin{align}
q_{tt} &= - q + \Delta_{h} q - h^{2 \ell \alpha} \, \beta \, q^{2\ell+1}, \label{KGsc2} 
\end{align}
which is a Hamiltonian PDE with associated Hamiltonian given by 
\begin{align}
\Ham K(q,p) &= \int_{I} \, \frac{p^2}{2} +  \left( \frac{q^2}{2} - \frac{q \; \Delta_{h} q}{2} \right) + h^{2\ell \alpha} \, \beta \, \frac{q^{2\ell+2}}{2\ell+2} \, \di y, \; \; \Delta_h \coloneqq 4 \sum_{k=1}^d \sinh^2\left(\frac{h \partial_{y_k}}{2}\right),  \label{Ham2KGcNLS} 
\end{align}
where $p$ is the variable canonically conjugated to $q$. 

We now exploit the change of coordinates $(q,p) \mapsto (\psi,\bar\psi)$ given by
\begin{align} \label{psi}
\psi &= \frac{1}{\sqrt{2}} (q- \mathrm{i} p),
\end{align}
so that the inverse change of coordinates is given by
\begin{align}
q = \frac{1}{\sqrt{2}} (\psi+\bar\psi), \; \; &\; \; \; p = \frac{\mathrm{i}}{\sqrt{2}} (\psi-\bar\psi) , \label{qp} 
\end{align}
while the symplectic form is given by $-\ii \, \di \psi \wedge \di\bar\psi$. 
With this change of variables the Hamiltonian takes the form
\begin{align}
\Ham K(\psi,\bar\psi) &= \int_{I}  \, \psi\bar{\psi} + \, \frac{ (\psi+\bar\psi)(-\Delta_h (\psi+\bar\psi)) }{4} + h^{2\ell \alpha} \, \beta \, \frac{ (\psi+\bar\psi)^{2\ell+2} }{ 2^{\ell+1}(2\ell+2) } \mathrm{d}y \nonumber \\
&\sim  \Ham{h}(\psi,\bar\psi) + h^{2} \Ham F(\psi,\bar\psi) + h^{4} \cR(\psi,\bar\psi), \label{AsExpPsi} \\
\Ham{h}(\psi,\bar\psi) &= \int_{I} \psi \, \bar\psi \, \di y, \nonumber \\
\Ham F(\psi,\bar\psi) &= \int_{I}  -\frac{(\psi+\bar\psi) \, [-\Delta(\psi+\bar\psi)] }{4} + h^{2(\ell\alpha-1)} \, \beta \, \frac{(\psi+\bar\psi)^{2\ell+2}}{ 2^{\ell+1}(2\ell+2) } \, \di y. \label{hFNLSpsi}
\end{align}
Observe that the remainder is of higher-order with respect to the leading part of the Hamiltonian, since
\begin{align*}
4 > 2\ell \alpha ; &\; \; \alpha < \frac{2}{\ell},
\end{align*}
which is true due to \eqref{NLSr}.

If we denote by $e_{n}(y)\coloneqq e^{\pi \mathrm{i} \; n \cdot y}$, $n \in \mathbb{Z}^d$, the Fourier basis on $L^2(I)$, by the decomposition
\begin{align} \label{xieta}
\psi(y) = 2^{-d/2} \sum_{n \in \mathbb{Z}^d} \xi_n e_n(y), \; \; &\; \; \; \overline{\psi(y)} = 2^{-d/2} \sum_{n \in \mathbb{Z}^d} \eta_n e_{-n}(y),
\end{align}
we get that 
\begin{align}
\Ham{h}(\xi,\eta) &= \sum_{n \in \mathbb{Z}^d} \xi_n \, \eta_n , \nonumber \\
\Ham F(\xi,\eta) &=  \sum_{n \in \mathbb{Z}^d} \frac{1}{4} (\xi_{n}+\eta_{-n}) \; \; |h|^2(\xi_{-n} +\eta_{n})  \nonumber \\
&\; \; \; \;  + h^{2(\ell\alpha-1)} \, \frac{\beta}{ 2^{\ell+1}(2\ell+2) } \sum_{\substack{n_1,\ldots,n_{2\ell+2} \in \mathbb{Z}^d\\ \sum_j  n_j = 0 } } (\xi_{n_1}+\eta_{-n_1}) \cdots (\xi_{n_{2\ell+2}}+\eta_{-n_{2\ell+2}}) . \label{hFNLSxieta}
\end{align}

Now we apply the averaging Theorem \ref{gavthm} to the Hamiltonian \eqref{AsExpPsi}, with $\mathcal{N} = \mathbb{Z}^d$, $\delta=h^2$, $\nu=\ell\alpha-1$; the assumption \eqref{nuAss} corresponds to
\begin{equation} \label{alphaBound0}
\frac{1}{2\ell} < \alpha \leq \frac{1}{\ell}.
\end{equation}
Observe that $\Ham{h}$ generates a periodic flow,
\begin{align}
-\mathrm{i} \partial_{t} \xi_n = \xi_n; \; &\; \; \xi_n(t ,y) = e^{\mathrm{i} t}\xi_{0,n}(y) , \; \; n \in \mathbb{Z}^d, \label{h0flowNLS}
\end{align}
and similarly for $\eta$.

\begin{proposition} \label{FavNLSprop}
Assume that the condition \eqref{alphaBound} holds true. Then the average of $\Ham F$ in \eqref{hFNLSxieta} with respect to the flow of $\Ham{h}$ \eqref{hFNLSxieta} is given by
\begin{align} 
& \la \Ham F \ra(\xi,\eta) \nonumber \\
&=  \sum_{n \in \mathbb{Z}^d}  \frac{1}{2} \, |n|^2 \; \xi_n \eta_n + h^{2(\ell\alpha-1)} \, \frac{\beta}{ 2^{\ell+1}(2\ell+2) } \sum_{ \substack{n_1,\ldots,n_{2\ell+2} \in \mathbb{Z}^d \\ \sum_j (-1)^{j-1} n_j = 0   } } \xi_{n_1} \eta_{n_2} \cdots \xi_{n_{2\ell+1}} \eta_{n_{2\ell+2}} . \label{FavNLS}
\end{align}
\end{proposition}

\begin{corollary} \label{NLScor}
Assume that the condition \eqref{alphaBound0} holds true. Then the equations of motion associated to $\Ham{h}(\xi,\eta)+h^2\la \Ham F \ra(\xi,\eta)$ are given by the following nonlinear Schr\"odinger equation (NLS),
\begin{align} \label{NLSeqFT} 
-\mathrm{i} \, \partial_t \xi_{n} &= \xi_n - \frac{h^2}{2} \, |n|^2 \xi_n + h^{2 \ell\alpha} \,  \frac{\beta}{2^{\ell+2}} \, \binom{2\ell+2}{\ell+1}  \, \sum_{ \substack{n_1,\ldots,n_{2\ell+1} \in \mathbb{Z}^d \\ \sum_j (-1)^{j-1} n_j = n    } }  \xi_{n_{1}} \eta_{n_{2}} \cdots \eta_{n_{2\ell}} \xi_{n_{2\ell+1}} , 
\end{align}
for all $n \in \mathbb{Z}^d$, and similarly for $\eta$.  Moreover, if
\begin{equation} \label{alphaBound}
\frac{1}{2\ell} < \alpha < \frac{1}{\ell}
\end{equation}
is satisfied, then Eq. \eqref{NLSeqFT} is a small-dispersion NLS.
\end{corollary}

\begin{remark} \label{rem:scaling}

By Corollary \ref{NLScor} and \eqref{xieta} we obtain that the function $\psi$ introduced in \eqref{psi} satisfies
\begin{align} \label{NLSeq}
-\mathrm{i} \, \psi_{t} &= \psi - \frac{1}{2} \; h^2 \; \Delta\psi + h^{2 \ell\alpha} \, \beta \, \binom{2\ell+2}{\ell+1} \frac{1}{2^{\ell+2}} \, |\psi|^{2\ell}\psi.
\end{align}

Observe that if we introduce the second time scale $\tau = \frac{1}{2} \, h^{2\ell\alpha} t$ and we consider $\varphi = e^{-\mathrm{i} \, 2\, \tau/{h^{2\ell\alpha}} } \, \psi$, then the rescaled Hamiltonian $\Ham K_2 = 2 h^{-2\ell\alpha} \, \Ham K_1$ has associated equations of motions equivalent to
\begin{align} \label{NLSeqResc}
-\mathrm{i} \, \varphi_{\tau} &= -  h^{2 (1-\ell\alpha)} \; \Delta\varphi + \beta \, \binom{2\ell+2}{\ell+1} \frac{1}{2^{\ell+1}} \, |\varphi|^{2\ell}\varphi ,
\end{align}
 which in Fourier reads
\begin{align}
-\mathrm{i} \, \partial_{\tau} \xi_{n} &= - h^{2(1-\ell\alpha)} \, |n|^2 \xi_n +  \beta \, 2^{-(\ell+1)} \, \binom{2\ell+2}{\ell+1}  \, \sum_{ \substack{n_1,\ldots,n_{2\ell+1} \in \mathbb{Z}^d \\ \sum_j (-1)^{j-1} n_j = n    } }  \xi_{n_{1}} \eta_{n_{2}} \cdots \eta_{n_{2\ell}} \xi_{n_{2\ell+1}} ,  \nonumber \\
& \nonumber \\
\mathrm{i} \, \partial_{\tau} \eta_{n} &= - h^{2(1-\ell\alpha)} \, |n|^2 \eta_n +  \beta \, 2^{-(\ell+1)} \, \binom{2\ell+2}{\ell+1}  \, \sum_{ \substack{n_2,\ldots,n_{2\ell+2} \in \mathbb{Z}^d \\ \sum_j (-1)^{j-1} n_j = n     } }  \eta_{n_{2}} \cdots \xi_{n_{2\ell+1}} \eta_{n_{2\ell+2}} , \label{NLSeqFTResc}
\end{align}
for all $n \in \mathbb{Z}^d$.
\end{remark}

\section{Dynamics of the NLS equation with small dispersion} \label{NLSdyn}

In this Section we recall some known facts on the dynamics of the defocusing small-dispersion NLS equation with periodic boundary conditions. The interested reader can find more detailed explanations and proofs in \cite{kuksin1995squeezing,kuksin1996growth,kuksin1997oscillations}.

Consider the defocusing NLS equation 
\begin{align} \label{deltaDNLSeq}
-\mathrm{i} \varphi_\tau &= -\varepsilon \; \Delta_{y} \varphi +  |\varphi|^{2\ell}\varphi, \; \; \varphi(0)=\varphi_0, \; \; \varepsilon >0, \; \; y \in \T^d\coloneqq(\R/(2\pi\Z))^d.
\end{align}

Eq. \eqref{deltaDNLSeq} is a PDE admitting a Hamiltonian structure: indeed, the associated Hamiltonian is given by
\begin{align}
\Ham H_{\mathrm{NLS}}(\varphi,\bar\varphi) &\coloneqq \int_{\T^d} \left( \frac{\varepsilon}{2} \, |\nabla\varphi|^2 + \frac{1}{2\ell+2} \, |\varphi|^{2\ell+2} \right) \; \di y. \label{1DNLSHam}
\end{align}

We first mention the following result, which states that if the dispersion in \eqref{deltaDNLSeq} is sufficiently small, then there exists a time at which the solution doubles its Sobolev norm (see Theorem 1 of \cite{kuksin1997oscillations} and its rescaled version in Sec. 3 of \cite{kuksin1997oscillations}).

\begin{proposition} \label{KuksinProp}

Consider \eqref{deltaDNLSeq} for $d \in \{1,2,3\}$, under odd periodic boundary conditions.

Let $m \geq 3$, and set $\nu \coloneqq 2\ell + 1/m$. Assume that the initial datum $\varphi_0$ of \eqref{deltaDNLSeq} satisfies
\begin{align*}
\| \varphi_0 \|_{L^2} = \varrho>0 , &\; \; \| \varphi_0 \|_{H^m} = r >0 .
\end{align*}
Then there exists $\lambda^{\star}=\lambda^{\star}(m) \in (0,1/\nu)$ such that for all $\lambda \in (0,\lambda^{\star})$, if the condition
\begin{align} \label{smallDispAss}
\varrho &\geq K \, \varepsilon^\lambda \, r^{1-2\ell\lambda}
\end{align}
holds true, then there exist $C>0$ independent from $m$ and $t_1>0$ satisfying
\begin{align} \label{t1EstKuk}
t_1 &\leq \, K_1 \, r^{\nu - 2 \ell} \, \varrho^{-\nu} \leq \, K_2 \, \varepsilon^{-\nu \lambda} \, r^{-2\ell(1-\nu\lambda)} , \; \; K_2 = K_1 \, K^{-\nu},
\end{align}
such that $\| \varphi(t_1) \|_{H^m} =2r$, and $\| \varphi(t_1) \|_{H^m} < 2r$ for $0 \leq t < t_1$.
\end{proposition}

A consequence of Proposition \ref{KuksinProp} is the following result, which describes the dynamics of solutions of small-dispersion NLS up to a time scale proportional to the inverse of the dispersion, showing that such solutions display a growth of their Sobolev norm (see Theorem 2 of \cite{kuksin1997oscillations}).

\begin{theorem} \label{KuksinThm}

Consider \eqref{deltaDNLSeq} for $d \in \{1,2,3\}$, under odd periodic boundary conditions.

Let $m \geq 3$, and set $\nu \coloneqq 2\ell + 1/m$. Then there exists $\lambda^{\star}=\lambda^{\star}(m) \in (0,1/\nu)$ such that for all $\lambda \in (0,\lambda^{\star})$, if we define the set
\begin{align} \label{AKuksin}
\mathcal{A} \, \coloneqq \, \mathcal{A}(m,\ell,\lambda) &\coloneqq \left\{ \varphi \in H^m(\mathbb{T}^d) \, | \, \|\varphi\|_{H^m} > K^{- \frac{1}{1-2\ell\lambda}} \, \varepsilon^{- \frac{\lambda}{1-2\ell\lambda}} \, \|\varphi\|_{L^2}^{\frac{1}{1-2\ell\lambda}} \right\},
\end{align}
the following holds true: if $\varphi$ is a solution of \eqref{deltaDNLSeq} such that $\varphi_0 \notin \mathcal{A}$ and $\|\varphi_0\|_{H^m}=r$, then there exists $T>0$ satisfying
\begin{align} \label{TEstKuk}
T &\leq \frac{1}{ 1-2^{-2\ell(1-\nu\lambda)} } \, K^{-\nu} \, \varepsilon^{-\nu \lambda} \, r^{-2\ell(1-\nu\lambda)} ,
\end{align}
such that $\varphi(T) \in \mathcal{A}$.
\end{theorem}

\begin{remark} \label{KukRem}

The fact that $\varphi_0 \notin \mathcal{A}$ can be seen as a condition on the dispersion $\varepsilon$, since it is equivalent to
\begin{align}
\varepsilon & \leq \varepsilon_0 \coloneqq \left( \frac{ \|\varphi_0\|_{L^2} }{ K \, \|\varphi_0\|_{H^m}^{1-2\ell\lambda} } \right)^{1/\lambda} , \label{CondSmallDisp}
\end{align}
see also the corresponding condition \eqref{smallDispAss} in Proposition \ref{KuksinProp}.

Kuksin also proved that for every $\mathtt{a}>0$ one can choose $K$ such that
\begin{align} 
& \lambda^{\star}(m) = \frac{1}{2\ell} \; \left( 1 - \frac{B_0}{2\ell m} \right) + \mathcal{O}\left(\frac{1}{m^2}\right), \; \; B_0=B_0(d,\ell,\mathtt{a}) = d\ell+4+\mathtt{a}, \; \; \text{as} \; \; m \to +\infty , \label{lambdaAsymp} \\
& \lim_{m\to\infty} \frac{1}{m} \; \frac{\lambda^\star(m)}{1-2\ell\lambda^\star(m)} =\frac{1}{B_0},  \label{lambdaAsymp2}
\end{align}
see Proposition 1 and pag. 346 in \cite{kuksin1997oscillations}. 

We point out that the time $T$ in \eqref{TEstKuk} is obtained by looking at the proof of Theorem 2 in \cite{kuksin1997oscillations}, which is in turn based on the estimate \eqref{t1EstKuk} in Proposition \ref{KuksinProp}. In the actual statement of the Theorem, the time $T$ is bounded by using that $\frac{1}{1-2^{-2\ell(1-\nu\lambda)} } \leq C_1 2\ell(1-\nu\lambda)$ and the asymptotics \eqref{lambdaAsymp}, leading to
\begin{align} \label{TEst2Kuk}
T &\leq C \, m \, K^{-\nu} \, \varepsilon^{-\nu \lambda} \, r^{-2\ell(1-\nu\lambda)} \, \leq \, C \, m \, K^{-1/\lambda} \, \varepsilon^{-1},
\end{align}
where the second estimate in \eqref{TEst2Kuk} follows from the fact that if $\varphi_0 \notin \mathcal{A}$, we have $r \geq \|\varphi_0\|_{L^2} \geq K \, \varepsilon^{\lambda} \, r^{1-2\ell\lambda}$, hence $r^{2\ell} \geq K^{1/\lambda} \, \varepsilon$.

\end{remark}

\begin{remark} \label{KukRemNorm}

Kuksin proved also that under the assumptions of Theorem \ref{KuksinThm}, if moreover $\|\varphi_0\|_{L^2} = 1$, then
\begin{align} \label{lambdaEstKuk}
\lambda^{\star}(m) &\leq \frac{1}{2} \, \frac{m}{1+2\ell m} ,
\end{align}
see Theorem 5 in \cite{kuksin1997oscillations}.

\end{remark}

\begin{remark} \label{applKGrem}

For our application to the KG lattice we have from \eqref{NLSeqResc} that $\varepsilon = h^{2(1-\ell\alpha)}$, which by \eqref{TEst2Kuk} it implies that in the original time scale $t$ 
\begin{align*}
t &\leq \mathcal{O}\left( \frac{1}{ h^{2 \left( \nu\lambda+\ell\alpha(1-\nu\lambda) \right)  } }\right) \; \leq \; \mathcal{O}\left( \frac{1}{h^2} \right) . 
\end{align*}
\end{remark}

\section{Approximation results} \label{ApprSec}

In this Section we show how to use the normal form equations in order to construct approximate solutions of \eqref{dDKGseq}, and we estimate the difference with respect to the true solutions with corresponding initial data. 

Recall that for \eqref{dDKGseq} we defined the energy of normal mode $E_k$, $k \in \Z^d_{N}$, in \eqref{EnNormModeKG}; also recall that for any $k \in \mathbb{Z}^d_N$ we defined the corresponding specific wave vector $\kappa \coloneqq \kappa(k)$ in \eqref{kappa}. 

In this section, we first point out a relation between the specific energy of normal mode $\cE_\kappa=2^d E_k/N_d$, and the Fourier coefficients of the solutions of the normal form equations, see Lemma \ref{EnSpecNLSLemma} and Proposition \ref{NLSphiProp}; this procedure has to be done carefully, since all wavevectors of the form $k+\left( (2N+1) \, n_1,\ldots, (2N+1) \, n_d \right)$ ($n = (n_1,\ldots,n_d) \in \mathbb{Z}^d$) contribute to the specific energy $\cE_\kappa$ of the discrete system. Then we prove that the approximate solutions approximate the specific energy of the true normal mode $\cE_\kappa$ up to the time-scale for which the continuous approximation is valid, see Proposition \ref{ApprPropdNLS}. For simplicity we present in this Section only the main part of the argument, and we defer the proof of Proposition \ref{NLSphiProp} to Appendix \ref{ApprEstSec11}. \\

The component $k$ of the Fourier transform of the continuous approximation at time $t$ will be denoted by $\hat{\mathsfit{Q}}(t,k)$ to distinguish it from the $k$-th component of the Fourier transform of the discrete variable $\hat{Q}_k(t)$. Since for the rescaled functions such as $q$ and $p$ this ambiguity does not hold, we prefer the short notation $\hat{q}_k(\tau)$ and $\hat{p}_k(\tau)$ for their Fourier transform. \\

Let $\beta >0$ and let $I$ be as in \eqref{I}, we define the Fourier coefficients of the function $q:I \to \R$ by

\begin{align} \label{FourierqcontKG}
\hat{q}_{k} \;&\coloneqq\; \frac{1}{2^{d/2}} \int_I q(y) \, e^{-\mathrm{i} \pi \, k \cdot y } \di y , \;\; k \in \mathbb{Z}^d ,
\end{align}
and similarly for the Fourier coefficients of the function $p$. 

\begin{lemma} \label{EnSpecNLSLemma}
Consider the lattice \eqref{Ham2KGs} in the regime \eqref{RegimeRel} and with interpolating function \eqref{NLSr}. Then for a state corresponding to $(q,p)$ one has
\begin{align} \label{SpecEnNormModeNLS}
\cE_\kappa &= \frac{ h^{2\alpha} }{2}  \sum_{ L \in \Z^d:h L \in (2\Z)^d } |\hat{p}_{K+L}|^2 + \omega_K^2 \left| \hat{q}_{K+L} \right|^2, \; \; \kappa = h \, K, \; \; K = ( K_1,\ldots,  K_d) \in \mathbb{Z}^d_N ,
\end{align}
where the $\omega_K$ are defined as in \eqref{FreqNormModeKG}.

\end{lemma}

\begin{proof}
Let us consider the smooth $(2N+1)^d$-periodic interpolating function $\mathsfit{Q}$ for $(Q_j)_{j \in \mathbb{Z}^d_N}$, and similarly $\mathsfit{P}$ for $(P_j)_{j \in \mathbb{Z}^d_N}$. We denote by
\begin{align} \label{FourierQcont}
\hat{Q}(j) &\coloneqq \frac{1}{ (2N+1)^{d/2}  } \int_{\left[-\left( N+\frac{1}{2} \right),\left( N+\frac{1}{2} \right) \right]^d } \mathsfit{Q}(x) e^{-\mathrm{i}  \frac{2\pi}{ 2N+1  } \, j \cdot x } \; \di x,
\end{align}
so that by the interpolation property we obtain
\begin{align*}
Q_j(t) =  \mathsfit{Q}(t,j) &=  \frac{1}{ (2N+1)^{d/2}  } \sum_{k \in \Z^d} \hat{\mathsfit{Q}}(t,k) e^{-\mathrm{i}  \frac{2\pi}{ 2N+1  }  \, j \cdot k }  \nonumber \\
&= \frac{1}{ (2N+1)^{d/2}  } \sum_{k=(k_1,\ldots,k_d) \in \Z^d_{N} } \, e^{-\mathrm{i} \frac{2\pi}{ 2N+1  } \, j \cdot k } \times \nonumber \\
&\; \; \times \left[ \sum_{n=(n_1,\ldots,n_d) \in \Z^d} \hat{ \mathsfit{Q} }( t, k_1+(2N+1)n_1,\ldots, k_d+(2N+1)n_d ) \right] ,
\end{align*}
hence
\begin{align} \label{FourierRel}
\hat{Q}_k(t) &= \sum_{n \in \Z^d} \hat{ \mathsfit{Q} }( t, k_1+(2N+1)n_1, \ldots, k_d+(2N+1)n_d ), \;\; \forall \, k \in \mathbb{Z}^d_N .
\end{align}
The relation between $\hat{\mathsfit{Q}}(t,k)$ and $\hat{q}_k(t)$ can be deduced from \eqref{NLSr}, which gives $\mathsfit{Q}(t,j) = h^{\alpha} q(t,h j)$ for all $j \in \mathbb{Z}^d_N$ . This implies that
\begin{align}
\hat{\mathsfit{Q}}(t,k) &=  \frac{1}{2^{d/2}} h^{d/2} \, \int_{ \left[ -\frac{1}{h},\frac{1}{h} \right]^d } \mathsfit{Q}(t,x) e^{-i\pi (h \, k \cdot x )} \di x \nonumber \\
&= \frac{1}{2^{d/2}} \, h^{\alpha + d/2} \, \int_{ \left[ -\frac{1}{h},\frac{1}{h} \right]^d } q(t, h x) e^{-i\pi (h \, k \cdot x )} \di x \nonumber \\
&\stackrel{\eqref{NLSr}}{=} \frac{1}{2^{d/2}} \, h^{\alpha - d/2} \int_{I} q(t,y) e^{-i \pi \, k \cdot y } \di y \; = \; h^{\alpha - d/2} \, \hat{q}_k(t) , \; \; \forall \, k \in \mathbb{Z}^d_N, \label{FourierRelQq}
\end{align}
and similarly
\begin{align} \label{FourierRelPp}
\hat{\mathsfit{P}}(t,k) &= h^{\alpha - d/2} \, \hat{p}_k(t) , \;\; \forall \, k \in \mathbb{Z}^d_N.
\end{align}
By using \eqref{enkappa} and \eqref{FourierRel}-\eqref{FourierRelPp} we have that for $\kappa = h \, K$, $K = ( K_1,\ldots, K_d) \in \mathbb{Z}^d_N$ (we omit the time dependence for simplicity; see also the proof of Lemma B.1 in \cite{bambusi2006metastability})
\begin{align*}
\cE_\kappa &\stackrel{\eqref{enkappa},\eqref{FourierRelQq}  }{=} h^{d} \, h^{2\alpha-d} \, \frac{1}{2} \,  \sum_{ L \in \Z^d:h L \in (2\Z)^d } |\hat{p}_{K+L}|^2 + \omega_K^2 \left| \hat{q}_{K+L} \right|^2 ,
\end{align*}
which leads to the thesis.
\end{proof}

\begin{proposition} \label{NLSphiProp}
Fix $s > \frac{d}{2} + 1$ and $0 < \delta < 1$. 

Consider the normal form system \eqref{NLSeqFTResc}, with $(\xi,\eta) \in \cW^{s}$, and denote by $\cE_\kappa$ the specific energy of the normal mode with index $\kappa$ as defined in \eqref{kappa}-\eqref{enkappa}. Then for any sufficiently small $h>0$  there exist positive constants $C$, $c_1$, $c_2$ and $M$ such that
\begin{align} \label{NLScoefvarphi}
\left| \frac{ \cE_\kappa }{ h^{2\alpha} } - \frac{\xi_K \; \eta_K}{2} \right| &\leq C h^{ c_1 } \|(\xi,\eta)\|_{\cW^{s}}^2 , \; \; c_1 = c_1(\delta) , 
\end{align}
for $\kappa = h \, K$, $K = ( K_1,\ldots, K_d) \in \mathbb{Z}^d_N$ and $|K_1|+\cdots+|K_d| \leq M=M(\delta)$, and
\begin{align} \label{NLSSpecEnEst}
|\cE_\kappa| &\leq C \, h^{ c_2 } \|(\xi,\eta)\|_{\cW^{s}}^2 , \; \; c_2=c_2(\delta,\alpha,s,d)  
\end{align}
for  $\kappa = h \, K$, $K = ( K_1,\ldots, K_d) \in \mathbb{Z}^d_N$ and $|K_1|+\cdots+|K_d| > M$. More explicitly, we have 
\begin{align}
c_1  = 2 \delta , \; \; M  &= h^{ - (1-\delta) } , \; \; c_2 = 2\alpha + (1-\delta) (2s-d-2) . \label{constAppr}
\end{align}

\end{proposition}

The proof of the above Proposition is deferred to Appendix \ref{ApprEstSec11}. \\

\begin{proposition} \label{ApprPropdNLS}
Let us consider \eqref{Ham2KGs} with $\beta>0$.

Let $m > s+s_0+\frac{d}{2} \gg 1$, where $s_0>0$ is the exponent for which Theorem \ref{gavthm} holds true. Let $\nu,\lambda^\star>0$ be defined as in Theorem \ref{KuksinThm}, and assume that $\lambda \in (0, \lambda^\star )$.

Let us assume that the initial datum satisfies \eqref{InDatumHyp} with $\alpha_0 < \alpha < 1/\ell$, where
\begin{align} 
\alpha_0 = \alpha_0(m,d,\ell,\lambda) &\coloneqq \frac{1}{\ell} \; \max\bigg( 1- (1+\nu\lambda)^{-1} \; \left( 1- \frac{d}{4} \right)  ,   \frac{d}{4} + \frac{1}{16} ,\frac{1}{2} \bigg) , \label{alpha0Bound}
\end{align}
and denote by $(Q_{j}(t),P_{j}(t))_{j \in \mathbb{Z}^d_{N}}$ the corresponding solution. Consider the approximate solution $\widetilde{\zeta}_a(t,x)=(\widetilde{\xi}_a(t,x), \widetilde{\eta}_a(t,x))$ with initial datum $\widetilde{\zeta}_{a,0}=(\widetilde{\xi}_{a,0}, \widetilde{\eta}_{a,0})\in \cW^{s+s_0} \cap (\ell^2_m \times \ell^2_m)$, and assume that $\widetilde{\zeta}_a \in \cW^{s+s_0}$ for all times. Moreover, let $0 < \delta < 1$, let $c_1=c_1(\delta)>0$ and $M=M(\delta)>0$ be as in Proposition \ref{NLSphiProp}, and assume that $\frac{1}{16} \, c_1 \, \in \left( 0 , 2\ell\alpha - \frac{d}{2} \right)$. \\

Then there exist $T_0$, $h_0=h_0( \alpha, \|( \widetilde{\xi}_{a,0},\widetilde{\eta}_{a,0} )\|_{\cW^{s+s_0}} )$ such that, if $h < h_0$, we have that there exists $C>0$ such that
\begin{align} \label{apprDiscrCont3}
\sup_j |Q_j(t) - \mathsfit{Q}_a(t,j)| + |P_j(t) - \mathsfit{P}_a(t,j)| &\leq C h^{\alpha+c_1/16}, \; \; \forall  \, |t| \leq \frac{T_0}{ h^{2 \left( \nu\lambda+\ell\alpha(1-\nu\lambda) \right)  } } ,
\end{align}
where $(\mathsfit{Q}_a,\mathsfit{P}_a)$ are given by considering $(\widetilde{\xi}_{a,0},\widetilde{\eta}_{a,0})$ and applying \eqref{xieta}, \eqref{qp} and \eqref{NLSr}. Moreover,
\begin{align} \label{LowModesApprNLS}
\left| \frac{ \cE_\kappa{(t)} }{ h^{2\alpha} } - \frac{\widetilde{\xi}_{a,K}(t) \widetilde{\eta}_{a,K}(t)}{2} \right| &\leq C \, h^{  c_1/16 } \, , \; \; \forall |t| \leq \frac{T_0}{ h^{2 \left( \nu\lambda+\ell\alpha(1-\nu\lambda) \right)  } }
\end{align}
for  $\kappa = h \, K$, $K = ( K_1,\ldots, K_d) \in \mathbb{Z}^d_N$ and $|K_1|+\ldots+|K_d| \leq M$. Moreover,  there exists $\varrho \coloneqq \varrho(\delta,\alpha,s,d)$ such that
\begin{align} \label{HighModesApprNLS}
|\cE_\kappa{(t)}| &\leq h^{ 2\alpha + \varrho } \, , \; \; \forall |t| \leq \frac{T_0}{ h^{2 \left( \nu\lambda+\ell\alpha(1-\nu\lambda) \right)  } } ,
\end{align}
for  $\kappa = h \, K$, $K = ( K_1,\ldots, K_d) \in \mathbb{Z}^d_N$ and $|K_1|+\ldots+|K_d| > M$.
\end{proposition}

\begin{remark} \label{ApproxAssRem}

Observe that the assumption $ c_1 < 16 \, \left( 2\ell\alpha - \frac{d}{2} \right) $ is equivalent to $\delta < 8 \left( 2\ell\alpha-\frac{d}{2} \right)$, which is automatically satisfied since $1 \leq 8 \left( 2\ell\alpha_0 - \frac{d}{2} \right)$; $\ell \alpha_0 \geq \frac{d}{4} + \frac{1}{16}$.

One can check from \eqref{alpha0Bound} that 
\begin{align*}
\frac{1}{2\ell} \leq \alpha_0 &< \frac{7}{8 \, \ell}  , \; \; \forall d \in \{ 1,2,3 \} , \; \; \forall \ell \geq 1, \; \; \forall m>1, \; \; \forall \lambda \in (0,\lambda^{\star}) .
\end{align*}

\end{remark}

\begin{proof}

The proof combines the argument in Appendix E of \cite{gallone2021metastability} (see also Appendix C of \cite{bambusi2006metastability} for the one-dimensional case) together with the argument in Appendix B of \cite{guardia2015growth} .

Exploiting the canonical transformation found in Theorem \ref{gavthm}, we also define
\begin{align}
\zeta_a \, &\coloneqq  \cT_{h^2}( \tilde\zeta_a ) \, = \, \tilde\zeta_a + \phi_a(\tilde\zeta_a),
\end{align}
where $\phi_a(\tilde\zeta_a)\coloneqq ( \phi_\xi(\tilde\zeta_a),\phi_\eta(\tilde\zeta_a) )$; by \eqref{CTthm} we have
\begin{align} \label{estRemThm3}
\sup_{\zeta \in \cB_{s}(R)} \|\phi_a(\zeta)\|_{\cW^{s}} &\leq C'_s h^{2\ell\alpha} \, R.
\end{align}

For convenience we define
\begin{align}
q_a(\tau,y) \coloneqq \frac{1}{\sqrt{2}} \left[ e^{\mathrm{i} \tau} \tilde{\xi}_a(\tau,y) + e^{-\mathrm{i}\tau} \tilde{\eta}_a(\tau,y)  \right] , &\;\; p_a(\tau,y) \coloneqq \frac{1}{\sqrt{2}\mathrm{i}} \left[ e^{\mathrm{i}\tau} \tilde{\xi}_a(\tau,y) - e^{-\mathrm{i}\tau} \tilde{\eta}_a(\tau,y)  \right] .   \label{qpappr3} 
\end{align}

We observe that the pair $(q_a,p_a)$ satisfies
\begin{align*}
h (q_a)_t = h p_a + h^{2\ell\alpha+1} \cR_q, \; &\; \; h (p_a)_t = - h q_a + h \Delta_1 q_a- h^{2\ell\alpha+1} \, \beta \, \overline{\pi_{0}}q_a^{2\ell+1} + h^{4\ell\alpha+1} \cR_p,
\end{align*}
where the operator $\Delta_1$ acts on the variable $x$, $\overline{\pi_{0}}$ is the projector on the space of the functions with zero average, and the remainders are functions of the rescaled variables $\tau$ and $y$ which satisfy
\begin{align} \label{remWs}
\sup_{\cB_{s+s_0}(R)} \| (\cR_q,\cR_p) \|_{\cW^s} &\leq C ,
\end{align}
due to estimate \eqref{Remthm} in Theorem \ref{gavthm}. \\

We now restrict the space variables to integer values; keeping in mind that $q_a$ and $p_a$ are periodic, we restrict to $j \in \Z^d_{N}$. Now let $Q=(Q_j)_{j \in  \Z^d_{N} }$ a finite sequence, we define the following norms 
\begin{align*}
\|Q\|_k &\coloneqq \left( \, \sum_{j \in \mathbb{Z}^d_{N}}  \, |Q_j|^k  \, \right)^{1/k} , \;\; k \geq 1.
\end{align*}

Now we consider the discrete model \eqref{Ham2KGs}: we rewrite it in the following form,
\begin{align}
\dot{Q}_j &= P_j , \label{DiscrEq31} \\
\dot{P}_j &= - Q_j + (\Delta_1Q)_j - \, \beta \, \overline{\pi_{0}}Q_j^{2\ell+1} , \label{DiscrEq32}
\end{align}
and we want to show that there exist two sequences $E=(E_j)_{j \in \Z^d_{N} }$ and $F=(F_j)_{j \in \Z^d_{N} }$ such that
\begin{align*}
Q \, = \, h^{\alpha} \, q_a + h^{\alpha+\gamma} E, &\; \; P \, = \, h^{\alpha} p_a + h^{\alpha+\gamma} F , \; \; \gamma \coloneqq c_1/16 ,
\end{align*}
($\gamma>0$ will be specified later) fulfills \eqref{DiscrEq31}-\eqref{DiscrEq32}. Therefore, we have that
\begin{align}
\dot{E} &= F - h^{2\ell\alpha-\gamma} \cR_q  \label{EqSeq31} \\
\dot{F} &= - E + \Delta_1E  - \beta \pi_0\, h^{2\ell\alpha} \left( \sum_{k=0}^{2\ell}  h^{(2\ell-k) \gamma} \,\binom{2\ell+1-k}{k} \, q_a^{k} E^{2\ell+1-k} \right) - h^{4\ell\alpha-\gamma}\cR_p , \label{EqSeq32}
\end{align}
where we impose initial conditions on $(E,F)$ such that $(q_a,p_a)$ has initial conditions corresponding to the ones of the true initial datum,
\begin{align*}
h^{\alpha} q_a(0,h j) + h^{\alpha+\gamma} E_{0,j} = Q_{0,j}, &\;\; h^{\alpha} p_a(0,h j) + h^{\alpha+\gamma} F_{0,j} = P_{0,j} . 
\end{align*}

We now define the operator $\partial_{i}$, $i=1,\ldots,d$,  by $(\partial_i f)_j \coloneqq f_j - f_{j-e_i^{(d)}}$ for all $f \in \ell^2_{0} \cap \ell^1_{0}$ such that $\| f \|_k < +\infty$ , $k=1,2$. 

\begin{itemize}
\item Claim 1: Let $\gamma >0$, we have
\begin{align*}
\|E_0\|_2 \leq C' h^{2\ell\alpha-\gamma-d/2}, \; \|F_0\|_2 &\leq C' h^{2\ell\alpha-\gamma-d/2}, \\ 
\|\partial_j E_0\|_2 \leq C' h^{2\ell\alpha+1-\gamma-d/2}, \; \|\partial_j F_0\|_2 &\leq C' h^{2\ell\alpha+1-\gamma-d/2}, \; j=1,\ldots,d.
\end{align*}
\end{itemize}
In order to prove Claim 1 we observe that
\begin{align*}
E_0 &= h^{\alpha} \frac{ \psi_a+\bar\psi_a-(\widetilde{\psi}_a+\bar{\widetilde{\psi}_a}) }{\sqrt{2} h^{\alpha+\gamma} }\, = \, h^{-\gamma} \, \frac{ \phi_\xi+\phi_\eta }{\sqrt{2}}, \\
F_0 &= h^{\alpha} \frac{ \psi_a-\bar\psi_a-(\widetilde{\psi}_a-\bar{\widetilde{\psi}_a}) ] }{\sqrt{2} \, \mathrm{i} \, h^{\alpha+\gamma} }\, = \, h^{-\gamma} \frac{ \phi_\xi-\phi_\eta }{\sqrt{2} \, \mathrm{i} }, \\
\end{align*}
and by \eqref{estRemThm3} we can deduce
\begin{align*}
\|E_0\|_2^2 &= \sum_{j \in \Z^d_{N} } |E_{0,j}|^2 \,\leq \, C \, N^{d} \, (h^{2\ell\alpha-\gamma})^2 \,  = \, C \, h^{4\ell\alpha-2\gamma-d}, \\
\|F_0\|_2^2 &= \sum_{j \in \Z^d_{N} }  |F_{0,j}|^2 \,\leq \, C \, N^{d} \, (h^{2\ell\alpha-\gamma})^2 \,  = \, C \, h^{4\ell\alpha-2\gamma-d}, \\
\|\partial_k E_0\|_2^2 &= \sum_{j \in \Z^d_{N} } |\partial_k E_{0,j}|^2 \, \leq \, \, C \, N^{d} \, (h^{2\ell\alpha+1-\gamma})^2 \, \leq \, C \, h^{4\ell\alpha+2-2\gamma-d}, \; \; k=1,\ldots,d, \\
\|\partial_k F_0\|_2^2 &= \sum_{j \in \Z^d_{N} } |\partial_k F_{0,j}|^2 \, \leq \, \, C \, N^{d} \, (h^{2\ell\alpha+1-\gamma})^2 \, \leq \, C \, h^{4\ell\alpha+2-2\gamma-d}, \; \; k=1,\ldots,d, 
\end{align*}
and this leads to the thesis.

\begin{itemize}
\item Claim 2: Fix $T_0 > 0$  and $K_\ast >0$. If \eqref{alpha0Bound} holds true, then for any $h < h_s$ and for any $\gamma>0$ satisfying 
\begin{align} \label{CondApprox}
\gamma &< 2\ell\alpha-\frac{d}{2}  
\end{align}
we have
\begin{align}
\|E\|_2^2 + \|F\|_2^2 +  \sum_{j=1}^d \|\partial_j E\|_2^2  &\leq K_\ast, \; \; |t| < \frac{T_0}{ h^{2 \left( \nu\lambda+\ell\alpha(1-\nu\lambda) \right)  } } .
\end{align}
\end{itemize}

To prove the claim, we define
\begin{align} 
\cG(E,F) &\coloneqq \sum_{j \in \Z^d_{N} }  \, \bigg[ \frac{F_j^2 + E_j^2 + E_j (-\Delta_1E)_j}{2}  + \beta \, h^{2\ell\alpha} \left( \sum_{k=0}^{2\ell}  h^{(2\ell-k) \gamma} \,\binom{2\ell+1-k}{k} \, q_{a,j}^{k} \frac{ E_j^{2(\ell+1)-k} }{2(\ell+1)-k} \right) \bigg] . \label{auxFuncClaim32}
\end{align}

Since there exists $C_1,C_2>0$ (depending on $\ell$) such that for $h>0$ sufficiently small 
\begin{align*}
 \left| \sum_{j \in \Z^d_{N} } \sum_{k=0}^{2\ell}  h^{(2\ell-k) \gamma} \,\binom{2\ell+1-k}{k} \, q_{a,j}^{k} \frac{ E_j^{2(\ell+1)-k} }{2(\ell+1)-k} \right| 
&\leq C_1 \,  \sum_{j \in \Z^d_{N} } \sum_{k=0}^{2\ell}  h^{(2\ell-k) \gamma}  \left|  E_j \right|^{2(\ell+1)-k}  \\
&=  C_1 \, \sum_{k=0}^{2\ell}  h^{(2\ell-k) \gamma}  \left\|  E \right\|_{2(\ell+1)-k}^{2(\ell+1)-k}  \\
&\leq C_1 \, \sum_{k=0}^{2\ell}  h^{(2\ell-k) \gamma}  \left\|  E \right\|_{2}^{2(\ell+1)-k} \\
&\leq C_2 \, h^{2\ell\alpha} \| E \|_2^2 ,
\end{align*}
we have that for $h>0$ sufficiently small 
\begin{align*}
\frac{1}{2} \cG(E,F) \leq \|E\|_2^2 + \|F\|_2^2 +  \sum_{j=1}^d \|\partial_j E\|_2^2   &\leq  2 \cG(E,F).
\end{align*}

Now we compute the time derivative of $\cG$. Exploiting \eqref{EqSeq31}-\eqref{EqSeq32}, we obtain
\begin{align}
\dot{\cG} &= \sum_{j} F_j \bigg[ -E_j + (\Delta_1E)_j - \beta \, h^{2\ell\alpha} \left( \sum_{k=0}^{2\ell}  h^{(2\ell-k) \gamma} \,\binom{2\ell+1-k}{k} \, q_{a,j}^{k} E_j^{2\ell+1-k} \right) - h^{4\ell\alpha-\gamma}(\cR_p)_j \bigg] \nonumber \\ 
&\;\; + \sum_{j}   (E_j - (\Delta_1E)_j) \left[ F_j - h^{2\ell\alpha-\gamma} (\cR_q)_j \right] \nonumber \\ 
&\; \; + \beta \, h^{2\ell\alpha} \, \sum_j   \, \sum_{k=0}^{2\ell}  h^{(2\ell-k) \gamma} \,\binom{2\ell+1-k}{k} \, q_{a,j}^{k} \, E_j^{2\ell+1-k} \, \left[ F_j - h^{2\ell\alpha-\gamma} (\cR_q)_j  \right] \nonumber \\ 
&\; \; + \beta \, h^{2\ell\alpha} \, \sum_j   \, \sum_{k=1}^{2\ell}  h^{(2\ell-k) \gamma} \,\binom{2\ell+1-k}{k} \, \frac{k}{2(\ell+1)-k} \, q_{a,j}^{k-1} \, h^{2\ell\alpha} \, \frac{\mathrm{d} q_{a,j}}{\mathrm{d}\tau} \, E^{2(\ell+1)-k}  , \nonumber 
\end{align}
which simplifies to
\begin{align}
\dot{\cG} &=  - h^{4\ell\alpha-\gamma} \sum_j  F_j (\cR_p)_j - h^{2\ell\alpha-\gamma} \sum_{j}  (E_j - (\Delta_1E)_j) (\cR_q)_j  \label{TimeDerAuxFunc42} \\
&\; \; - \beta \, h^{4\ell\alpha-\gamma} \, \sum_j   (\cR_q)_j \, \sum_{k=0}^{2\ell}  h^{(2\ell-k) \gamma} \,\binom{2\ell+1-k}{k} \, q_{a,j}^{k} \, E_j^{2\ell+1-k}  \label{TimeDerAuxFunc43} \\
&\; \; + \beta \, h^{4\ell\alpha} \, \sum_j   \, \sum_{k=1}^{2\ell}  h^{(2\ell-k) \gamma} \,\binom{2\ell+1-k}{k} \, \frac{k}{2(\ell+1)-k} \, q_{a,j}^{k-1} \,\frac{\mathrm{d} q_{a,j}}{\mathrm{d}\tau} \, E_j^{2(\ell+1)-k}   \label{TimeDerAuxFunc44} 
\end{align}
In order to estimate \eqref{TimeDerAuxFunc42}-\eqref{TimeDerAuxFunc44}, we notice that
\begin{align*}
\sup_j |(\Delta_1 E)_j| &\leq 2 \sup_j   \sum_{k=1}^d |(\partial_k E)_j| \leq 4d \sqrt{\cG}, \\
\|\cR_q\|_2^2 \, \leq \, \sum_j  |(\cR_q)_j|^2 & \leq C \, N^{d} \sup_y |\cR_q(y)|^2 \, \leq \, C h^{-d}, \\
\|\cR_p\|_2^2 &\leq \, C h^{-d},
\end{align*}
and that $ |(\partial_i\cR_q)_j| \leq h \sup_y \left| \frac{\partial \cR_q}{\partial y} (y) \right|$, which implies 
\begin{align*}
\|\partial_i\cR_q\|_2^2 & \leq  C h^{2-d}.
\end{align*}

Now, we can estimate \eqref{TimeDerAuxFunc42} by
\begin{align} 
C \, h^{4\ell\alpha-\gamma-d/2}  \, \mathcal{G}^{1/2} + C \, h^{2\ell\alpha-\gamma-d/2} \, \mathcal{G}^{1/2} &= C \, h^{2\ell\alpha-\gamma-d/2} (1+h^{2\ell\alpha})  \, \mathcal{G}^{1/2} , \label{EstTimeDerAuxFunc42}
\end{align}
so that this term is small if $\gamma < 2\ell\alpha-\frac{d}{2}$, which in turn implies that $\alpha > \frac{d}{4\ell}$.

Next, recalling that $q_{a}$ is bounded, we can estimate \eqref{TimeDerAuxFunc43} by
\begin{align} 
&C \, \beta \, h^{4\ell\alpha-\gamma-d/2} \left[ \mathcal{G}^{1/2} + h^{\gamma} \mathcal{G} +  \sum_{k=0}^{2\ell-2} h^{(2\ell-k)\gamma} \| E \|_{2}^{2\ell+1-k} \right] \nonumber \\
&\leq C \, \beta \, h^{4\ell\alpha-\gamma-d/2} \left[ \mathcal{G}^{1/2} + h^{\gamma} \mathcal{G} +   h^{2\ell \gamma } \mathcal{G}^{3/2} \right] \nonumber \\ 
&\leq C \, \beta \, \left[h^{4\ell\alpha-\gamma-d/2} \mathcal{G}^{1/2} + h^{4\ell\alpha-d/2} \mathcal{G} +  h^{4\ell\alpha+\gamma-d/2} \, \mathcal{G}^{3/2}   \right]  . \nonumber 
\end{align}

Finally, we can estimate \eqref{TimeDerAuxFunc44} by
\begin{align} 
& C \, \beta \, h^{4\ell\alpha} \left[ \mathcal{G} + \sum_{k=1}^{2\ell-1} h^{(2\ell-k)\gamma} \| E \|_{2}^{2\ell+2-k}  \right] \nonumber \\ 
&\leq C \, \beta \, h^{4\ell\alpha} \left[ \mathcal{G} +  h^{ \gamma } \mathcal{G}^{3/2} \right] \nonumber \\ 
&\leq C \, \beta \, \left[ h^{4\ell\alpha} \mathcal{G} + h^{4\ell\alpha+\gamma} \mathcal{G}^{3/2}  \right] . \nonumber 
\end{align}

Hence, as long as $\cG < 2 K_\ast$ and provided that condition \eqref{CondApprox} holds true,

\begin{align}
  \left| \dot{\cG} \right| &\leq C  K_\ast^{1/2} \, \left[  h^{2\ell\alpha-\gamma-d/2} (1+h^{2\ell\alpha}) +\beta \, h^{4\ell\alpha-\gamma-d/2}  \right] \nonumber \\ 
&\; \; \; \; + C \, \beta \,  h^{2\ell\alpha} \left[  h^{2\ell\alpha-d/2}  + K_\ast^{1/2} h^{2\ell\alpha+\gamma-d/2} +h^{2\ell\alpha} + K_\ast^{1/2} h^{2\ell\alpha+\gamma} \right] \mathcal{G} , \label{EstTimeDer32} 
\end{align}
where we can rewrite the term on the last line of Eq. \eqref{EstTimeDer32} in the following form,
\begin{align*}
& C \, \beta \,  h^{2\ell\alpha + 2\nu\lambda(1-\ell\alpha) } \, \mathcal{G} \,  \bigg[ h^{2\ell\alpha- 2\nu\lambda(1-\ell\alpha) - d/2}   (1 +  K_\ast^{1/2} h^{\gamma} )    (1 + h^{d/2}) \bigg] 
\end{align*}
where
\begin{align*}
2\ell\alpha- 2\nu\lambda(1-\ell\alpha) - d/2 &> 0
\end{align*}
for
\begin{align*}
\ell \alpha &> (1+\nu\lambda)^{-1} \; \left( \frac{d}{4} + \nu \lambda \right) = 1- (1+\nu\lambda)^{-1} \; \left( 1- \frac{d}{4} \right) ,
\end{align*}
(we also observe that $1- (1+x)^{-1} \; \left( 1- \frac{d}{4} \right) > \frac{d}{4}$ for all $0<x<1$ and $d=1,2,3$), from which we can derive the assumption \eqref{alpha0Bound}.

By applying Gronwall's lemma we get
\begin{align} 
\cG(t) &\leq \cG(0) e^{ C \, \beta \,  h^{2\ell\alpha + 2\nu\lambda(1-\ell\alpha) } \, t \,  \left[ h^{2\ell\alpha- 2\nu\lambda(1-\ell\alpha) - d/2}   (1 +  K_\ast^{1/2} h^{\gamma} )    (1 + h^{d/2}) \right]  } \nonumber \\
&\; \;  + e^{ C \, \beta \,  h^{2\ell\alpha + 2\nu\lambda(1-\ell\alpha) } \, t \,  \left[ h^{2\ell\alpha- 2\nu\lambda(1-\ell\alpha) - d/2}   (1 +  K_\ast^{1/2} h^{\gamma} )    (1 + h^{d/2}) \right]  }  \nonumber \\
&\; \; \; \; \times  C \, \beta \,  h^{2\ell\alpha + 2\nu\lambda(1-\ell\alpha) } \, t  \,  \bigg[ h^{2\ell\alpha- 2\nu\lambda(1-\ell\alpha) - d/2}   (1 +  K_\ast^{1/2} h^{\gamma} )    (1 + h^{d/2}) \bigg]   \nonumber \\
&\; \; \; \; \times  C  K_\ast^{1/2} \, \left[  h^{2\ell\alpha-\gamma-d/2} (1+h^{2\ell\alpha}) +\beta \, h^{4\ell\alpha-\gamma-d/2}  \right] , \label{Gronwall} 
\end{align}
and from \eqref{Gronwall} we can deduce the thesis. 

\end{proof}

We conclude with some important remarks about the proof of Proposition \ref{ApprPropdNLS}.

\begin{remark} \label{RemEstRem}

We stress the fact that \eqref{remWs}, unlike the bound (C.7) in \cite{bambusi2006metastability} and  the estimate after Eq. (260)-(261) in \cite{gallone2021metastability}, means that the remainder (which is evaluated on an approximate solution) is bounded in the space $\cW^s$, which is based on $\ell^1_s$; see also the estimate (3.19) in Theorem 4 of \cite{guardia2015growth}, where such bound is used in the context of proving growth of Sobolev norms for the NLS equation.
\end{remark}

\begin{remark} \label{Rem1Norm1}

Using the notations of Claim $1$ in the proof of Proposition \ref{ApprPropdNLS}, we can exploit the assumption on the regularity of the initial datum $\widetilde{\zeta}_{a,0}$ and Lemma \ref{lpRem} in order to prove the following estimates,
\begin{align*}
\|E_0\|_1 \leq C' h^{2\ell\alpha-\gamma-d}, \; \|F_0\|_1 &\leq C' h^{2\ell\alpha-\gamma-d},  \\
\|\partial_j E_0\|_1 \leq C' h^{2\ell\alpha+1-\gamma-d}, \; \|\partial_j F_0\|_1 &\leq C' h^{2\ell\alpha+1-\gamma-d}, \; j=1,\ldots,d.
\end{align*}

\end{remark}

\begin{remark} \label{Rem2Norm1}

Notice that the function $\mathcal{G}$ introduced in \eqref{auxFuncClaim32} is adapted to the norm $\| \cdot \|_2$. In principle, we could have introduced an analogous function adapted to the norm $\| \cdot \|_1$; in that case, however, using the above argument we would have obtained worse estimates with respect to $h$, and in turn we would have obtained worse estimates on the exponent $\alpha_0$ (see also Remark \ref{Rem1Norm1}, and compare it to the Claim 1).
\end{remark}

\section{Statement and Proof of the Main Result} \label{ThmProof}

Here we state and we prove our main result, Theorem \ref{dDNLSrThm}. The main issue in the proof of the theorem is to determine to which extent we can show that the discrete Klein-Gordon lattice \eqref{Ham2KGs} displays an energy cascade phenomenon as its approximating PDE, the small-dispersion NLS \eqref{NLSeq}. 

Now we state our main result.\\

\begin{theorem} \label{dDNLSrThm}
Consider \eqref{Ham2KGs} with $\beta > 0$, $\ell \geq 1$, and let $ \frac{7}{8} \, \ell^{-1} < \alpha_1 < \left( 1 - \frac{1}{2^6} \right) \, \ell^{-1} $.

Then there exists $m_0 \gg 1$ such that for all $m \geq m_0$ there exist positive $\lambda_0$ and $\lambda_1$ such that for any $\lambda \in (\lambda_0,\lambda_1)$ the following holds true.

There exist $\alpha_0 >0$ and $\delta_0>0$ such that for any $\delta \in ( \delta_0 , 1 )$ and for any $\alpha \in (\alpha_0,\alpha_1)$ there exists a sufficiently small $h_0>0$ such that for all positive $h < h_0$, if the initial datum satisfies
\begin{align} \label{InDatumHyp}
\mathcal{E}_{\kappa_0}(0) = C_0 \; h^{2\alpha} , &\; \; \mathcal{E}_{\kappa}(0) = 0 , \; \; \forall \kappa \neq \kappa_0 ,
\end{align}
then there exist $C_1>0$, $K > 0$, a sufficiently large $T_m >0$ satisfying
\begin{align} \label{TmCascade}
 T_m &\leq C_1 \, K \, h^{-2 \ell\alpha - 2(2\ell+1/m)\lambda(1-\ell\alpha) }  ,
\end{align}
and  $c>0$ such that
\begin{align} \label{Cascade}
\sum_{\kappa \in h \, \mathbb{Z}^d_N : |\kappa|  \leq h^{\delta} } |\kappa|^{2m} \, \cE_\kappa(T_m)  &\geq \frac{c}{8} \;  h^{2\alpha - 2(1-\ell\alpha) \frac{2\lambda}{1-2\ell\lambda} +2m } .
\end{align}

\end{theorem}

\begin{remark} \label{QuantDepRem}

From the proof of Theorem \ref{dDNLSrThm} one can highlight the dependence of the quantities $m_0$, $\lambda_0$, $\lambda_1$, $\alpha_0$, $\delta_0$, $h_0$, $K$, $T_m$ and $c$ from the parameters introduced in the statement of the Theorem \ref{dDNLSrThm}. 

More precisely, we have $m_0=m_0(d,\ell,s,s_0)  > s+s_0+\frac{d}{2}$, where the exponents $s_0$ and $s$ appear in the statement of Theorem \ref{gavthm}. We also have $\lambda_0=\lambda_0(m,d,\ell,\alpha_1)$, $\lambda_1=\lambda_1(m,d,\ell)$, $\alpha_0 = \alpha_0(m,d,\ell,\lambda)$, $\delta_0=\delta_0(m,d,\ell,\alpha_1,\lambda)$, $h_0=h_0(m,d,\delta,\ell,\alpha,\lambda)$, $K=K(\ell,m,\lambda)$, $T_m=T_m(h,\ell,\alpha,\lambda)$ and $c=c(\ell,\lambda)$.

\end{remark}

\smallskip

\begin{remark} \label{MainThmRem}

Theorem \ref{dDNLSrThm}, and in particular \eqref{Cascade}, gives an analytical description of an energy cascade for the KG lattice. As mentioned in comment $\mathrm{i.}$ of the introduction, the regimes of specific energy described in Theorem \ref{dDNLSrThm} are above the regime for which one can prove the existence of metastability phenomena, but they also are below the so-called thermodynamic limit. 

We observe that the fact that $d \leq 3$ is a consequence both of Theorem \ref{KuksinThm} about the growth of Sobolev norms for the small-dispersion NLS equation, and of Proposition \ref{ApprPropdNLS} regarding the approximation of the normal modes of the lattice with the Fourier modes of the PDE.
\end{remark}

\smallskip

\begin{remark} \label{ConstantsRem}

Here we mention where the different quantities mentioned in the statement of Theorem \ref{dDNLSrThm} appear in the paper, in order to highlight their role and to give a better explanation of the interplay between the different intermediate results.

The exponent $m_0$ is defined at the end of the proof of Theorem \ref{dDNLSrThm}.

The exponent $\lambda_0$ is defined at the end of the proof of Theorem \ref{dDNLSrThm}, by using the asymptotics \eqref{lambdaAsymp2}, while the exponent $\lambda_1$ is defined at the beginning of the proof of Theorem \ref{dDNLSrThm}.

We observe that the exponent $\alpha_0$ is defined in the statement of Proposition \ref{ApprPropdNLS}, see \eqref{alpha0Bound}, and is essential in order to ensure that the approximation estimates \eqref{apprDiscrCont3}, \eqref{LowModesApprNLS} and \eqref{HighModesApprNLS} hold true. The exponent $\alpha_1$ is defined due to Remark \ref{ApproxAssRem} and is needed in the proof of Theorem \ref{dDNLSrThm} in order to ensure that the remainders are smaller than the growth for the small-dispersion NLS stated in Theorem \ref{KuksinThm}. The precise upper bound on $\alpha_1$ is irrelevant, the important property is that $\alpha_1$ is bounded away from $\ell^{-1}$.

The threshold $\delta_0>0$ is defined at the end of the proof of Theorem \ref{dDNLSrThm}; the fact that the threshold is positive implies that in Proposition \ref{ApprPropdNLS} we prove that the normal modes of the lattice are close to the Fourier modes of the approximating PDE only if the associated frequency is not too high, see \eqref{LowModesApprNLS}. 

The time $T_m$ appearing in \eqref{Cascade} is related to the time for which the growth of Theorem \ref{KuksinThm} is achieved. Its estimate \eqref{TmCascade} is a consequance of \eqref{TEst2Kuk} (which in turn follows from \eqref{lambdaAsymp}) and Remark \ref{applKGrem}.

\end{remark}

\smallskip

\begin{proof}[Proof of Theorem \ref{dDNLSrThm}]

Let $m \geq m_0 > s+s_0+\frac{d}{2}$, where $s \gg 1$ and $s_0>0$ is the exponent for which Theorem \ref{gavthm} holds true; we will specify the value of $m_0$ later. Recalling the notations of Theorem \ref{KuksinThm}, we define $\nu(m) \coloneqq  2\ell + 1/m$, and we denote $\lambda_1 \coloneqq  \min( \lambda^\star(m) , \lambda^\star(m+1) )$. \\

We consider the lattice \eqref{Ham2KGs}, and we consider an initial datum of the form 
\begin{align*}
\mathcal{E}_{\kappa_0}(0) = C_0 \; h^{2\alpha} , &\; \; \mathcal{E}_{\kappa}(0) = 0 , \; \; \forall \kappa \neq \kappa_0 ,
\end{align*}
with $\alpha_0 < \alpha < \alpha_1 $ (recall that $\alpha_0$ is defined in \eqref{alpha0Bound}) and with corresponding Fourier coefficients $(\hat{Q}_{0,k},\hat{P}_{0,k})_{k \in \mathbb{Z}^d_N}$ given by \eqref{fourierQ}. 

Moreover, we define an interpolating function for the initial datum $(Q_0,P_0)$ by 
\begin{align*}
\mathsfit{Q}_0(y) &= \frac{1}{ (2N+1)^{d/2} } \sum_{K: h |K| \leq 1 } \hat{Q}_{0,K} \, e^{\ii\pi  \, h \, K \cdot y },
\end{align*}
and similarly for $y \mapsto \mathsfit{P}_0(y)$. \\

By applying \eqref{NLSr} and Lemma \ref{EnSpecNLSLemma}, the initial datum $(\hat{Q}_{0,k},\hat{P}_{0,k})_{k \in \mathbb{Z}^d_N}$ of the lattice \eqref{Ham2KGs} corresponds to $(\hat{q}_{0,J},\hat{p}_{0,J})_{J \in \mathbb{Z}^d} \in  \cW^{s+s_0} \cap (\ell^2_{m+1} \times \ell^2_{m+1})$, which, by using the change of variables \eqref{xieta}, takes the form of an initial datum in the Fourier space $(\widetilde{\zeta}_{a,0,L})_{L \in \mathbb{Z}^d} \in \cW^{s+s_0} \cap (\ell^2_{m+1} \times \ell^2_{m+1})$.

\begin{remark} \label{NormalizedRem}

Due to the specific form of the initial datum and to the change of variables \eqref{xieta}, we have that
\begin{align*}
\sum_{L \in \mathbb{Z}^d} \widetilde{\xi}_{a,0,L} \,\widetilde{\eta}_{a,0,L} &= 1 ,
\end{align*}
namely, the solution is normalized, in the language of \cite{kuksin1997oscillations}.
\end{remark}

By Remark \ref{KukRem}, Theorem \ref{KuksinThm} and Proposition \ref{KuksinProp} we have that for all positive $\lambda<\lambda_1$, if the dispersion is sufficiently small, namely
\begin{align} 
h &\leq h_{0,m} \; \coloneqq  \; \left[ K^{-2/\lambda} \; \left( \sum_{L \in \mathbb{Z}^d} |L|^{2m} \widetilde{\xi}_{a,0,L} \,\widetilde{\eta}_{a,0,L}  \right)^{-(1-2\ell\lambda)/\lambda} \right]^{ \frac{1}{4(1-\ell\alpha)}  }  , \label{mSmallDisp}
\end{align}
there exists a (large) $T_m^{\ast} >0$ satisfying
\begin{align*}
T_m^{\ast} &\leq C_1 \; K^{-\nu(m)} \;  m  \; h^{-2\ell\alpha - 2\nu(m) \, \lambda (1-\ell\alpha)} \; \left( \sum_{L \in \mathbb{Z}^d} |L|^{2m} \widetilde{\xi}_{a,0,L} \,\widetilde{\eta}_{a,0,L}  \right)^{- \ell (1-\nu(m) \, \lambda)} , 
\end{align*}
such that we have the growth
\begin{align} \label{mGrowth} 
\sum_{L \in \mathbb{Z}^d} |L|^{2m} \, \widetilde{\xi}_{a,L}(T_m^{\ast}) \,\widetilde{\eta}_{a,L}(T_m^{\ast}) &= c \;  h^{- 2(1-\ell\alpha) \frac{2\lambda}{1-2\ell\lambda}  }  ,  
\end{align}
where $c=c(\ell,\lambda)\coloneqq K^{-\frac{2}{1-2\ell\lambda}}>0$. Moreover, we can assume that $T_m^{\ast}$ is the \emph{minimal} time at which the growth \eqref{mGrowth} is attained. Notice that \eqref{mGrowth} implies
\begin{align*}  
\sum_{L \in \mathbb{Z}^d} |L|^{2(m+1)} \, \widetilde{\xi}_{a,L}(T_m^{\ast}) \,\widetilde{\eta}_{a,L}(T_m^{\ast}) &\geq c \;  h^{- 2(1-\ell\alpha) \frac{2\lambda}{1-2\ell\lambda}  }  . 
\end{align*}

Similarly, for all positive $\lambda<\lambda_1$, if the smallness condition
\begin{align} \label{msigmaSmallDisp}
h &\leq h_{0,m+1} \; \coloneqq  \; \left[ K^{-\frac{2}{\lambda}} \;  \left( \sum_{L \in \mathbb{Z}^d} |L|^{2(m+1)} \widetilde{\xi}_{a,0,L} \,\widetilde{\eta}_{a,0,L}  \right)^{-\frac{1-2\ell\lambda}{\lambda}} \right]^{ \frac{1}{4(1-\ell\alpha)}  }  
\end{align}
holds true, then there exists a positive (large) $T_{m+1}^{\ast} \leq T_m^{\ast}$ satisfying
\begin{align*}
T_{m+1}^{\ast} &\leq C_1 \; K^{-\nu(m+1)} \;  (m+1)  \; h^{-2\ell \alpha - 2\nu(m+1) \, \lambda (1-\ell\alpha)} \times \\
&\qquad\qquad \times \left( \sum_{L \in \mathbb{Z}^d} |L|^{2(m+1)} \widetilde{\xi}_{a,0,L} \,\widetilde{\eta}_{a,0,L}  \right)^{- \ell (1-\nu(m+1) \, \lambda)} , 
\end{align*}
such that we have the growth
\begin{align} \label{msigmaGrowth2}
\sum_{L \in \mathbb{Z}^d} |L|^{2(m+1)} \, \widetilde{\xi}_{a,L}(T_{m+1}^{\ast}) \,\widetilde{\eta}_{a,L}(T_{m+1}^{\ast}) &= c \;  h^{- 2(1-\ell\alpha) \frac{2\lambda}{1-2\ell\lambda} }. 
\end{align}

We now denote
\begin{align*}
h_0^{\ast} &\coloneqq  h_0^{\ast}(m) \, = \, \min( h_{0,m} , h_{0,m+1} ) ,
\end{align*}
where $h_{0,m}$ and $h_{0,m+1}$ are defined in \eqref{mSmallDisp} and \eqref{msigmaSmallDisp} , respectively.

Hence, for all positive $h < h_0^{\ast}$ and $\lambda < \lambda_1$, there exists $T_m \in [ T_{m+1}^{\ast}, T_{m}^{\ast}]$ such that the following property holds true: there exists $\vartheta \coloneqq  \vartheta(m+1,T_m) \geq 1$ such that
\begin{align}
\sum_{L \in \mathbb{Z}^d} |L|^{2m} \, \widetilde{\xi}_{a,L}(T_m) \,\widetilde{\eta}_{a,L}(T_m) &\geq \frac{3}{4} \, c \;  h^{- 2(1-\ell\alpha) \frac{2\lambda}{1-2\ell\lambda}  }  ,   \label{mGrowth2} \\ 
\sum_{L \in \mathbb{Z}^d} |L|^{2(m+1)} \, \widetilde{\xi}_{a,L}(T_m) \,\widetilde{\eta}_{a,L}(T_m) &= \vartheta \, c \;  h^{- 2(1-\ell\alpha) \frac{2\lambda}{1-2\ell\lambda}  }  . \label{msigmaGrowth3} 
\end{align}

We now exploit that for all $M>0$ and for all $\sigma>0$
\begin{align} \label{Estmmsigma}
\| \mathrm{id} - \Pi_M \|_{\ell^2_{m+\sigma} \times \ell^2_{m+\sigma} \to \ell^2_m \times \ell^2_m} &\leq (M+1)^{-\sigma} 
\end{align}
(see Lemma 3.4 in \cite{gallone2021metastability}). \\

If we fix $\delta \in ( 0 , 1) $ and we apply \eqref{Estmmsigma} to \eqref{msigmaGrowth3} with $M=h^{-(1-\delta)} -1$ and $\sigma=1$, we obtain 
\begin{align} \label{mHighEst}
\sum_{L \in \mathbb{Z}^d: |L| > h^{-(1-\delta)} } |L|^{2m} \, \widetilde{\xi}_{a,L}(T_m) \,\widetilde{\eta}_{a,L}(T_m) &\leq h^{ 2 (1-\delta) } \, \vartheta \, c \;  h^{- 2(1-\ell\alpha) \frac{2\lambda}{1-2\ell\lambda}  }  ,  
\end{align}
so that by \eqref{mGrowth2} and by \eqref{mHighEst} we get
\begin{align} \label{mLowEst1}
\sum_{L \in \mathbb{Z}^d: |L| \leq h^{-(1-\delta)} } |L|^{2m} \, \widetilde{\xi}_{a,L}(T_m) \,\widetilde{\eta}_{a,L}(T_m) &\geq \left( \frac{3}{4} - h^{ 2 (1-\delta) } \, \vartheta \right)  \,c \;  h^{- 2(1-\ell\alpha) \frac{2\lambda}{1-2\ell\lambda}  } .
\end{align}
Therefore, there exists $h_0^{\star}\coloneqq h_0^{\star}(m) < h_0^{\ast}$ such that for all positive $h < h_0^{\star}$ satisfying
\begin{align*}
h^{ 2(1-\delta) } &\leq \frac{1}{4 \, \vartheta}
\end{align*}
we have
\begin{align} \label{mLowBound}
\frac{3}{4} - h^{ 2 (1-\delta) } \, \vartheta &\geq \frac{1}{2} .
\end{align}

Therefore, by \eqref{mLowEst1} and \eqref{mLowBound} we can deduce that there exists $h_0^{\star}\coloneqq h_0^{\star}(m,\delta) \leq h_0^{\ast}$ such that for all positive $h < h_0^{\star}$ and $\lambda < \lambda_1$, there exists $T_m \in [ T_{m+1}^{\ast}, T_{m}^{\ast}]$ such that the following property holds true,
\begin{align} \label{mLowEst2}
\sum_{L \in \mathbb{Z}^d: |L| \leq h^{-(1-\delta)} } |L|^{2m} \, \widetilde{\xi}_{a,L}(T_m) \,\widetilde{\eta}_{a,L}(T_m) &\geq \frac{1}{2}  \,c \;  h^{- 2(1-\ell\alpha) \frac{2\lambda}{1-2\ell\lambda}  } .
\end{align}

Then by applying Proposition \ref{ApprPropdNLS}, we have that if $\lambda \in (0, \lambda_1 )$, then there exists a positive $h_0$ such that, if $h < h_0$, we have that there exists $C>0$ 
\begin{align*} 
\sup_j |Q_j(t) - \mathsfit{Q}_a(t,j)| + |P_j(t) - \mathsfit{P}_a(t,j)| &\leq C h^{\alpha+c_1/16}, \; \; \forall  \, |t| \leq T_m , 
\end{align*}
where $(\mathsfit{Q}_a,\mathsfit{P}_a)$ are given by considering $(\widetilde{\xi}_{a,0},\widetilde{\eta}_{a,0})$ and applying \eqref{xieta}, \eqref{qp} and \eqref{NLSr}; moreover, by \eqref{LowModesApprNLS} we have that
\begin{align*} 
\left| \frac{ \cE_\kappa{(t)} }{ h^{2\alpha} } - \frac{\widetilde{\xi}_{a,L}(t) \widetilde{\eta}_{a,L}(t)}{2} \right| &\leq C \, h^{ \delta/8 } \, , \; \; \forall |t| \leq T_m , 
\end{align*}
for $\kappa = h \, L$ and $|L| \leq h^{-(1-\delta)}$ (recall that  $\kappa \in h \,  \mathbb{Z}^d_N$).

By combining \eqref{mLowEst2} and \eqref{LowModesApprNLS} we have that
\begin{align*} 
\frac{c}{4} \;  h^{- 2(1-\ell\alpha) \frac{2\lambda}{1-2\ell\lambda}  } &\leq  \sum_{L \in \mathbb{Z}^d: |L| \leq h^{-(1-\delta)} } |L|^{2m} \, \frac{ \widetilde{\xi}_{a,L}(T_m) \,\widetilde{\eta}_{a,L}(T_m) }{2} \nonumber \\
&= \sum_{L \in \mathbb{Z}^d: |L| \leq h^{-(1-\delta)} } |L|^{2m} \, \left[  \frac{ \cE_{h  L}(T_m) }{ h^{2\alpha} } + \frac{ \widetilde{\xi}_{a,L}(T_m) \,\widetilde{\eta}_{a,L}(T_m) }{2} - \frac{ \cE_{h  L}(T_m) }{ h^{2\alpha} } \right] ,
\end{align*}
hence
\begin{align}
\frac{c}{4} \;  h^{- 2(1-\ell\alpha) \frac{2\lambda}{1-2\ell\lambda}  } &\leq \sum_{L \in \mathbb{Z}^d: \kappa =h \, L, |L| \leq h^{-(1-\delta)} } |L|^{2m} \, \frac{ \cE_\kappa(T_m) }{ h^{2\alpha} }  \label{SumCascade1} \\
&\;\;\;\; + \sum_{L \in \mathbb{Z}^d: \kappa =h \, L, |L| \leq h^{-(1-\delta)} } |L|^{2m} \,\left[ \frac{ \widetilde{\xi}_{a,L}(T_m) \,\widetilde{\eta}_{a,L}(T_m) }{2} - \frac{ \cE_\kappa(T_m) }{ h^{2\alpha} } \right] . \label{SumCascade2} 
\end{align}
Now, by \eqref{apprDiscrCont3} and by \eqref{LowModesApprNLS} we have that \eqref{SumCascade2} is bounded by $\mathcal{O}( h^{\delta/8-(1-\delta)(2m+d)} )$. Hence \eqref{SumCascade2} is of lower order compared to the left-hand side term of \eqref{SumCascade1} if
\begin{align*}
\left( \frac{1}{8}+2m+d \right)\delta -(2m+d) &> -2(1-\ell\alpha_1) \frac{2\lambda}{1-2\ell\lambda} ; 
\end{align*}
\begin{align*}
\delta &> \delta_{0}=\delta_{0}(m,d,\ell,\alpha_1,\lambda) \; \coloneqq  \; 1-\frac{1/8}{2m+d+1/8} \, \left[ 1+ (1-\ell\alpha_1) \frac{2\lambda}{1-2\ell\lambda} \right] .
\end{align*}
Therefore, by taking $\delta > \delta_0$ and by recalling the asymptotics \eqref{lambdaAsymp} (where we choose  $B_0=B_0(d,\ell,1)$) and \eqref{lambdaAsymp2} we can deduce that for all $\mathtt{c} \in (0, 1/(2B_0) )$ there exists a sufficiently large $m_0=m_0(d,\ell,\mathtt{c},s,s_0)$ such that for all $m \geq m_0$
\begin{align*}
\left| m^{-1} \, \frac{\lambda_1}{1-2\ell\lambda_1} - \frac{1}{B_0}  \right| &< \mathtt{c} .
\end{align*}
Hence there exists $\lambda_0=\lambda_0(m,d,\ell,\alpha_1)>0$ such that for any $\lambda \in (\lambda_0, \lambda_1)$, for any $\delta \in (\delta_0,1)$ and for any $\alpha \in (\alpha_0,\alpha_1)$ there exists a sufficiently small $h_0=h_0(m,d,\delta,\ell,\alpha,\lambda)>0$ such that for all positive $h < h_0$, if the initial datum is of the form \eqref{InDatumHyp}, then 
\begin{align} \label{CascadeTmp}
\sum_{L \in \mathbb{Z}^d: \kappa =h \, L, |L| \leq h^{-(1-\delta)} } |L|^{2m} \, \frac{ \cE_\kappa(T_m) }{ h^{2\alpha} } &\geq \frac{c}{8} \;  h^{- 2(1-\ell\alpha) \frac{2\lambda}{1-2\ell\lambda}  } .
\end{align}

Finally, from \eqref{CascadeTmp} we can deduce \eqref{Cascade} by observing that
\begin{align*}
\sum_{L \in \mathbb{Z}^d: \kappa=h \, L, |L| \leq h^{-(1-\delta)} } |L|^{2m} \, \cE_\kappa(T_m) &= \sum_{L \in \mathbb{Z}^d: \kappa=h \, L, |L| \leq h^{-(1-\delta)} } \frac{ |\kappa|^{2m} }{ h^{2m} } \, \cE_\kappa(T_m) \\
& =  \sum_{\kappa \in h \, \mathbb{Z}^d_N : |\kappa|  \leq h^{\delta} } \cE_\kappa(T_m) \, \frac{ |\kappa|^{2m} }{ h^{2m} } ,
\end{align*}
where the last equality holds true because $|L| \leq h^{-(1-\delta)}$ .

\end{proof}

\begin{remark} \label{rem:deltaRem}

In order to obtain \eqref{Cascade} it is essential to assume both that $ \frac{1}{\ell} -  \alpha_1$ is bounded away from zero and that $\lambda > \lambda_0$, otherwise we would obtain that $1-\delta_0 = \mathcal{O}(1/m)$ as $m \to \infty$.
\end{remark}

\section*{Acknowledgements} 

The author would like to thank Matteo Gallone, Patrick G\'erard, Tiziano Penati and Antonio Ponno for useful comments and suggestions. The author would like to thank everyone who gave suggestions which led to an improvement of the paper. This project has received funding from the European Union's Horizon 2020 research and innovation programme under the Marie Sk\l{}odowska-Curie grant agreement No 101034255.

\begin{figure}[!h]
\includegraphics[width = .1\textwidth]{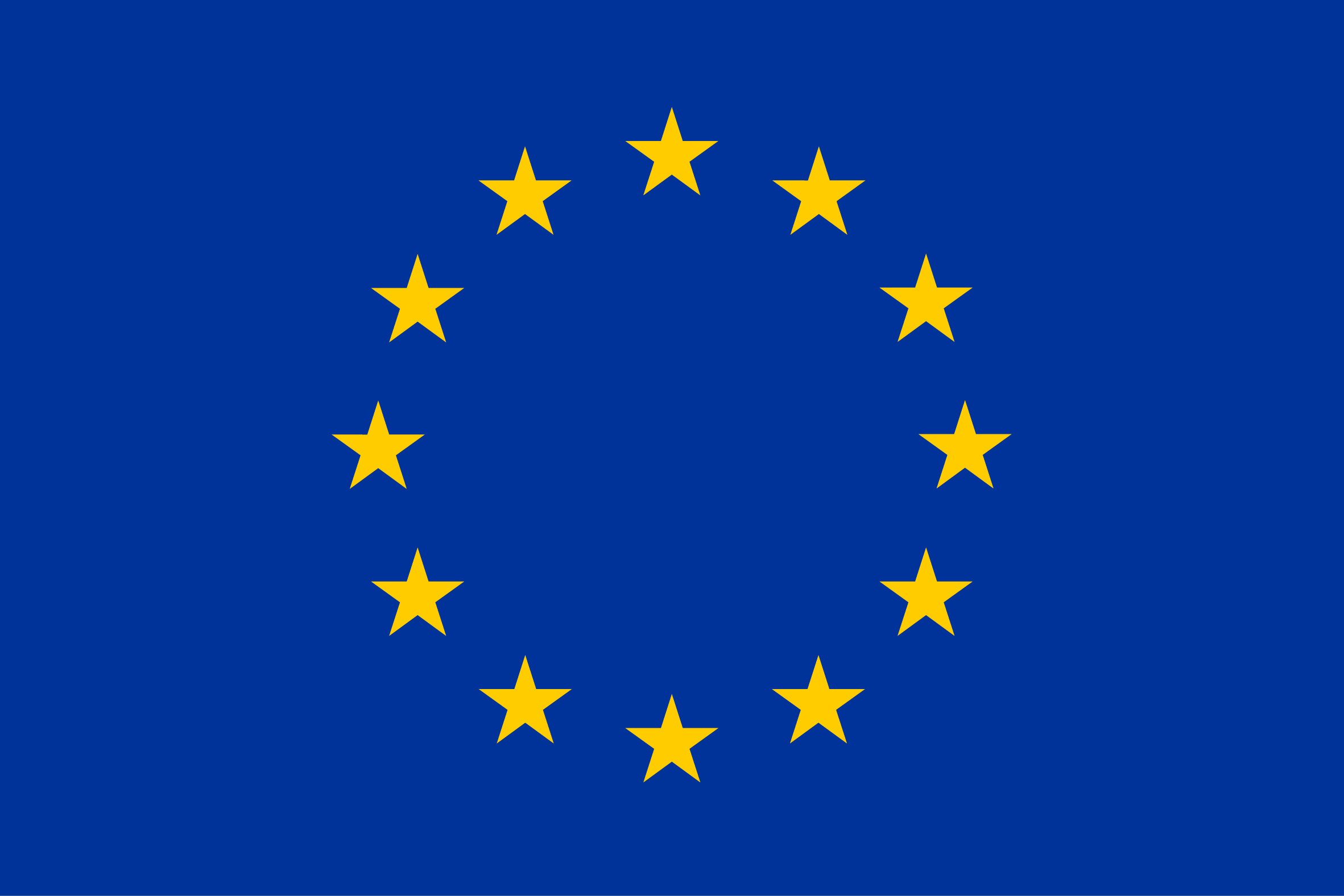}
\end{figure}

\begin{appendix}

\section{Proof of Lemma \ref{NFest}} \label{BNFest}

This appendix is devoted to the proof of the Lemma \ref{NFest}, which is a key step to normalize the system \eqref{truncsys}. Its proof is an adaptation of the proof of Theorem 4.4 in \cite{bambusi1999nekhoroshev} together with some results from Sec. III.3 and Sec.IV.1 of \cite{faou2012geometric}; it is based on the method of Lie transform, briefly recalled in the following. Throughout this Section, we consider  $s \geq s_1$ to be a fixed quantity.

Given $k_0 \geq 2$ and an auxiliary function $\chi \in \mathcal{P}_{k_0}$ analytic on $\cW^{s}$, 
we consider the auxiliary differential equation
\begin{align} \label{auxDE}
\dot \zeta &= X_\chi(\zeta)
\end{align}
and denote by $\Phi^t_\chi$ its flow at time $t$. 

\begin{definition}
The map $\Phi_\chi  \coloneqq  \Phi^1_\chi$ is called the \emph{Lie transform} 
generated by $\chi$.
\end{definition}

Given $k_1 \geq 2$ and $\Ham G \in \mathcal{P}_{k_1}$ analytic on $\cW^{s}$, 
let us consider the differential equation
\begin{align} \label{orDE}
\dot \zeta &= X_{\Ham G}(\zeta),
\end{align}
where by $X_{\Ham G}$ we denote the vector field of $\Ham G$. Now define
\begin{align*}
\Phi_\chi^\ast \Ham G(\tilde\zeta) & \coloneqq  \Ham G \circ \Phi_\chi(\tilde\zeta).
\end{align*}
By exploiting the fact that $\Phi_\chi$ is a canonical transformation, we have that in the new variable $\tilde\zeta$ defined by $\zeta=\Phi_\chi(\tilde\zeta)$ equation \eqref{orDE} is equivalent to
\begin{align} \label{pullbDE}
\dot{\tilde\zeta} &= X_{ \Phi_\chi^\ast \Ham G }(\tilde\zeta).
\end{align}

Using the relation
\begin{align}\label{CompositionRelation}
\frac{\di}{\di t} \Phi_\chi^\ast \Ham G &= \Phi_\chi^\ast \{\chi, \Ham G\}, 
\end{align}
and the Poisson bracket formalism $\{\Ham G_1,\Ham G_2\}(\zeta) \coloneqq  \di \Ham G_1(\zeta) [ X_{\Ham G_2}(\zeta)]$ we formally get
\begin{equation} \label{lieseries}
\Phi^\ast_\chi \Ham G = \sum_{l=0}^\infty \Ham G_l, \; \; \Ham G_0  \coloneqq  \Ham G, \; \; \Ham G_{l}  \coloneqq  \frac{1}{l} \{\chi,\Ham G_{l-1}\}, \; \; l \geq 1.
\end{equation}

In order to estimate the vector field of the terms appearing in \eqref{lieseries}, we exploit the following results (see also Proposition III.6 in \cite{faou2012geometric}). We recall that we defined the norm $\| \cdot \|$ in \eqref{defSupNorm}.

\begin{lemma}\label{lem:liebrest}
Let $k_0,k_1 \geq 2$, $R>0$, and assume that $\chi \in \mathcal{P}_{k_0}$, $\Ham G \in \mathcal{P}_{k_1}$ are analytic on $\cB_{s}(R)$, as well as their vector fields. Then we have that 
$\{\chi,\Ham G\} \in \mathcal{P}_{k_0+k_1-2}$ is analytic on $\cB_{s}(R)$, and
\begin{align} 
\| \{\chi,\Ham G\} \| &\leq 2 k_0 k_1 \|\chi\| \, \| \Ham G \| , \label{liebrest} 
\end{align}
\begin{align}
&\sup_{\zeta \in \cB_{s}(R)} \|X_{ \{\chi,\Ham G\} }(\zeta)\|_{\cW^{s}} \nonumber \\
&\leq 4 \, (k_0+k_1-2)(k_0+k_1-3)^s \, k_0 k_1 \; \|\chi\| \; \|\Ham G\| \; R \, \sup_{\zeta \in \cB_{s}(R)} \, \max_{n=1,\ldots,k_0+k_1-4} \|\zeta\|_{\cW^s}^n . \label{liebrest2} 
\end{align}
\end{lemma}

\begin{proof}
First we prove the thesis assuming that $\chi$ and $\Ham G$ are homogeneous polynomials of degrees $k_0$ and $k_1$ respectively, with coefficients $a_{\bfj}$, $\bfj \in \cI_{k_0}$, and $b_{\bfh}$, $\bfh \in \cI_{k_1}$, respectively. Then it follows from definition that $\{\chi,\Ham G\}$ is a monomial of degree $k_0+k_1-2$ satisfying the zero momentum condition. Moreover, write
\begin{align*}
\{\chi , \Ham G\}(\zeta) &= \sum_{\bfj \in \cI_{k_0+k_1-2}} c_{\bfj} \zeta_{\, \bfj} ,
\end{align*}
where $c_{\bfj}$ is a sum of coefficients $a_{\bfr} b_{\bfs}$, for which there exists $a \in \cN$ and $\sigma \in \{\pm 1\}$ such that
\begin{align*}
(a,\sigma) \subset \bfr \in \cI_{k_0}, &\; \; (a,-\sigma) \subset \bfs \in \cI_{k_1},
\end{align*}
and such that if for example $(a,\sigma)= \mathtt{r}_1$ and $(a,-\sigma) = \mathtt{s}_1$, we necessarily have
\begin{align*}
( \mathtt{r}_2,\ldots, \mathtt{r}_{k_0}, \mathtt{s}_2,\ldots, \mathtt{s}_{k_1}) &= \bfj .
\end{align*}

If $\chi = \sum_{i=2}^{k_0} \chi_i$ and $\Ham G = \sum_{j=2}^{k_0} \Ham G_j$, where $\chi_i$ and $\Ham G_j$ are homogeneous polynomials of degree $i$ and $j$ respectively, then we have
\begin{align*}
\{ \chi , \Ham G \} &= \sum_{n=2}^{k_0+k_1-2} \sum_{i+j-2=n} \{ \chi_i, \Ham G_j \} ,
\end{align*}
where all the polynomials $\{ \chi_i, \Ham G_j \}$ are homogeneous polynomials of degree $i+j-2$; hence we can conclude that $\{ \chi, \Ham G \} \in \cP_{k_0+k_1-2}$.

By applying \eqref{EstFG4} and \eqref{EstFG2} we can deduce the thesis.

\end{proof}

\begin{lemma}\label{lem:Liesrest}
Let $k_0,k_1 \geq 2$, $R>0$, and assume that $\chi \in \mathcal{P}_{k_0}$ and $\Ham G \in \mathcal{P}_{k_1}$ are analytic on $\cB_{s}(R)$, as well as their vector fields. 
Let $l \geq 1$, and consider $\Ham G_l$ as defined in \eqref{lieseries}, then  $\Ham G_l \in \mathcal{P}_{l(k_0-2)+k_1}$ is analytic on $\cB_{s}(R)$ as well as its vector field, and
\begin{align} 
\| \Ham G_l \| &\leq (2e \, k_0 \, \|\chi\|)^l \, k_1 \, \|\Ham G\| , \label{lieserest}
\end{align}
\begin{align}
\sup_{\zeta \in \cB_{s}(R)} \|X_{ \Ham G_\ell }(\zeta)\|_{\cW^{s}} & \leq 2 \, ( l(k_0-2)+k_1)^l \; ( l(k_0-2)+k_1-1)^{s} \; (2e \, k_0 \, \|\chi\|)^l \, k_1 \, \|\Ham G\| \; R \, \times \nonumber \\
&\qquad \qquad \times \sup_{\zeta \in \cB_{s}(R)} \max_{n=1,\ldots,l(k_0-2)+k_1-2} \|\zeta\|_{\cW^s}^n . \label{lieserest2} 
\end{align}
\end{lemma}

\begin{proof}
Fix $l \geq 1$. We look for a sequence $C^{(l)}_m$ such that
\begin{align*}
\|\Ham G_m\| &\leq C^{(l)}_m, \; \; \forall m \leq l.
\end{align*}
Lemma \ref{lem:liebrest} ensures that the following sequence satisfies this property.
\begin{align*}
C^{(l)}_0  \coloneqq  k_1 \, \|\Ham G\|, \; \; C^{(l)}_m &= \frac{2 l}{m} C^{(l)}_{m-1} \, k_0 \, \|\chi\| .
\end{align*}
One has
\begin{align*}
C^{(l)}_l &= \frac{1}{l!} \left( 2l \, k_0 \, \| \chi \| \right)^l  \; k_1 \, \|\Ham G\| ,
\end{align*}
and by using the inequality $l^l < l! \; e^l$ one obtains the estimate \eqref{lieserest}.
\end{proof}

\begin{lemma} \label{homeqlemma}
  Assume that $\Ham{G} \in \mathcal{P}_{k_1}$ is analytic on $\cB_{s}(R)$ as well as its vector field, 
and that $\Ham{h}_0$ satisfies (PER). 
Then there exists $\chi \in \mathcal{P}_{k_1}$ analytic on $\cB_{s}(R)$ and $\Ham{Z} \in \mathcal{P}_{k_1}$ analytic 
on $\cB_{s}(R)$ with $\Ham Z$ in normal form, namely $\{\Ham{h}_0,\Ham Z\}=0$, such that
\begin{align} \label{homeq}
\{ \chi,\Ham{h}_0 \} \; + \; \Ham G \; &= \; \Ham Z.
\end{align}
Such $\Ham Z$ and $\chi$ are given explicitly by
\begin{align}\label{eq:ExplicitZ}
\Ham Z(\zeta)\; = \; \frac{1}{T} \int_0^T \Ham G(\Phi^t_{\Ham{h}_0}(\zeta)) \, \di t \, ,
\end{align}
\begin{align}\label{eq:ExplicitChi}
\chi(\zeta) &= 
\frac{1}{T} \int_0^T t \, \left[ \Ham Z(\Phi^t_{\Ham{h}_0}(\zeta))-\Ham G(\Phi^t_{\Ham{h}_0}(\zeta)) \right] \, \di t \, .
\end{align}
Furthermore, we have that the vector fields of $\chi$ and $\Ham Z$ are analytic on $\cB_{s}(R)$, and satisfy 
\begin{align} 
\|\Ham Z\| \, \leq \, \|\Ham G\| , &\; \; \sup_{\zeta \in \cB_{s}(R)} \|X_{\Ham Z}(\zeta)\|_{\cW^{s}} \leq  \sup_{\zeta \in \cB_{s}(R)} \|X_{\Ham G}(\zeta)\|_{\cW^{s}}, \nonumber \\
\|\chi\| \, \leq \, 2T \, \|\Ham G\| , &\; \; \sup_{\zeta \in \cB_{s}(R)} \|X_\chi(\zeta)\|_{\cW^{s}} \leq  2T  \sup_{\zeta \in \cB_{s}(R)} \|X_{\Ham G}(\zeta)\|_{\cW^{s}} \label{vfhomeq}.
\end{align}
\end{lemma}

\begin{proof}
We check directly that the solution of \eqref{homeq} is \eqref{eq:ExplicitChi}. Indeed,
\begin{align*}
\{ \chi,\Ham{h}_0 \}(\zeta) &= \frac{\di}{\di s} \Big|_{s=0} \chi(\Phi^s_{\Ham{h}_0}(\zeta)) \\
&= \frac{1}{T} \int_0^{T} t \frac{\di}{\di s} \Big|_{s=0} \left[ \Ham Z(\Phi^{t+s}_{\Ham{h}_0}(\zeta))-\Ham G(\Phi^{t+s}_{\Ham{h}_0}(\zeta)) \right] \di t \\
&= \frac{1}{T} \int_0^{T} t \frac{\di}{\di t} \left[ \Ham Z(\Phi^{t}_{\Ham{h}_0}(\zeta))-\Ham G(\Phi^{t}_{\Ham{h}_0}(\zeta)) \right] \di t \\
&= \frac{1}{T} \left[ t \Ham Z(\Phi^{t}_{\Ham{h}_0}(\zeta))- t \Ham G(\Phi^{t}_{\Ham{h}_0}(\zeta)) \right]_{t=0}^{T} - \frac{1}{T} \int_0^{T} \left[ \Ham Z(\Phi^{t}_{\Ham{h}_0}(\zeta))-\Ham G(\Phi^{t}_{\Ham{h}_0}(\zeta)) \right] \di t \\
&= \Ham Z(\zeta)-\Ham G(\zeta).
\end{align*}
In the last step we used the explicit expression of $\Ham Z$ provided in \eqref{eq:ExplicitZ}. Finally, the first estimate in \eqref{vfhomeq} follows from the explicit expression of $\Ham Z$ in \eqref{eq:ExplicitZ} while for the second estimate we write explicitly the vector field $X_\chi$: 
\begin{align*}
X_\chi(\zeta) &= 
\frac{1}{T} \int_0^T t \, D\Phi^{-t}_{\Ham{h}_0}(\Phi_{\Ham{h}_0}^t(\zeta)) \circ X_{\Ham Z-\Ham G}(\Phi^t_{\Ham{h}_0}(\zeta)) \, \di t \, .
\end{align*}
Hypothesis (PER) guarantees that $\Phi_{\Ham{h}_0}^t$ as well as its derivatives and the inverses are uniformly bounded as operators from $\cW^{s}$ into itself. Moreover, for any $t \in \mathbb{R}$, the map $\zeta \mapsto \Phi_{\Ham{h}_0}^t(\zeta)$ is a diffeomorphism of $\cB_{s}(R)$ into itself. {Using the fact that $\Phi_{\Ham h_0}^t$ is an isometry, we have}
\[
	\begin{split}
		\sup_{\zeta \in \cB_{s}(R)} \Vert X_\chi(\zeta) \Vert_{\cW^{s}} &\leq T  \sup_{\zeta \in \cB_{s}(R)} \left(\Vert X_{\Ham Z}(\zeta) \Vert_{\cW^{s}} + \Vert X_{\Ham G}(\zeta) \Vert_{\cW^{s}} \right) \leq 2T \sup_{\zeta \in \cB_{s}(R)} \Vert X_{\Ham G}(\zeta) \Vert_{\cW^{s}}
	\end{split}
\]
where in the last step we used the first inequality in \eqref{vfhomeq}. 
\end{proof}

Now we introduce the following definition regarding Lie derivatives of functions defined on $\cW^s$ (see Definition IV.1-Proposition IV.3 in \cite{faou2012geometric}).

\begin{definition} \label{def:LieDer}
	
Let $\cU_{s} \subset \cW^{s}$ open, and let $\Ham H \in C^1(\cU_{s},\C)$ be a Hamiltonian. Now let $\Ham{g} \in C^\infty(\cW_s,\mathbb{C})$, we define the Lie derivative $\cL_{\Ham{H}}[\Ham{g}]$ in the following way,
\begin{align*}
\cL_{\Ham{H}}[\Ham{g}] & \coloneqq  \{ \Ham{H} , \Ham{g} \} .
\end{align*}
Similarly, let $\Ham{Y} \in C^\infty(\cW^s,\cW^s)$ with $Y=(Y_j)_{j \in \cZ}$, then we set
\begin{align*}
(\cL_{\Ham{H}}[\Ham{Y}])_j & \coloneqq  \{ \Ham{H} , \Ham{Y}_j \} , \; \; j \in \cZ.
\end{align*}
\end{definition}

\begin{lemma} \label{lem:LieDer}
	
Let $s_1 \geq s \geq 0$, and $k \geq 3$. Assume that $\Ham{P} \in \mathcal{P}_{k}$, $\Ham{Y} \in C^\infty(\cW^{s_1},\cW^s)$. Then $\cL_{\Ham{P}}[\Ham{Y}] \in C^\infty(\cW^{s_1},\cW^s)$ .
\end{lemma}
\begin{proof}
	
Observe that $\Ham{Y} = (Y_j)_{j \in \cZ} \in C^1(\cW^{s_1},\cW^s)$ and that
\begin{align*}
\cL_{\Ham{P}}[\Ham{Y}](\zeta) &= \nabla_{\zeta}\Ham{Y}(\zeta) \cdot X_{\Ham{P}}(\zeta),
\end{align*}
hence the fact that $X_{\Ham{P}}(\zeta) \in \cW^{s_1}$ for $\zeta \in \cW^{s_1}$ implies that $\cL_{\Ham{P}}[\Ham{Y}](\zeta) \in \cW^s$. The thesis follows by a similar argument on the higher-order derivatives of $\cL_{\Ham{P}}[\Ham{Y}]$.

\end{proof}

Let $\Phi^t_{\chi}$ be the flow associated to the equation
\begin{align*}
\zeta_t &= X_{\chi}(\zeta), \; \; \text{in} \; \; \cW^s,	
\end{align*}
for $|t| \leq 1$. Now if $\zeta(t) \in \cW^s$ for all $s \geq 0$, we can consider the Taylor expansion
around $t=0$.

\begin{lemma} \label{lem:TaylorExp}

Let $R>0$ and $k_0 \geq 3$. Let $\chi \in \cP_{k_0}$, and let $\zeta(t)=\Phi^t_{\chi}$ be a solution of the Hamiltonian vector field $X_{\chi}$ in $\cW^s$. Assume that $\zeta \in \cB_{s}(R)$, then there exists $T_0 \in  (0,1)$ such that there exists $C>0$ such that for all $t \in [0,T_0]$ we have
\begin{align*}
\left\| \Phi^t_{\Ham{P}}(\zeta) - \zeta - t \, \cL_{\Ham{P}}[\Ham{Id}](\zeta) \right\|_{\cW^s} &\leq C t^2.
\end{align*}
	
\end{lemma}
\begin{proof}
	
By Lemma \ref{lem:LieDer} we have that $\cL_{\Ham{P}}[\Ham{Id}] \in C^\infty(\cW^s,\cW^s)$, and the result is obtained as in the finite dimensional case.	
\end{proof}

We just mention that Lemma \ref{lem:TaylorExp} can be generalized to the Taylor expansion up to any finite order (see Proposition IV.3 in \cite{faou2012geometric}).

\begin{lemma}\label{lem:VFCanTr}

Let $s \geq s_1$, and assume that $\chi \in \mathcal{P}_{k_0}$, $\Ham F \in \mathcal{P}_{k_1}$  are analytic on $\cB_{s}(R)$. 
Then for $|t| \leq 1$
\begin{align}
\sup_{\cB_{s}(R)} \|X_{ (\Phi^t_\chi)^\ast \Ham{F} - \Ham{F} }(\zeta)\|_{\cW^{s}} &= \sup_{\cB_{s}(R)} \|X_{ \Ham{F} \circ \Phi^t_\chi - \Ham{F} }(\zeta)\|_{\cW^{s}} \nonumber \\
&\stackrel{\eqref{liebrest}}{\leq} 4 (k_0 + k_1 -2)(k_0 + k_1 -3)^s \, \|\chi\| \, \|\Ham{F}\| \, R \, \sup_{\zeta \in \cB_s(R)} \max_{n=1,\ldots,k_0+k_1-4} \|\zeta\|_{\cW^s}^n .  \label{vfest}
\end{align}

\end{lemma}
\begin{proof}
	
The thesis follows from the fact that
\begin{align*}
\Ham{F} \circ \Phi^t_\chi &= F + t \int_0^1 \cL_{\chi}[\Ham{F}] \circ \Phi_\chi^{t \tau} \, \mathrm{d}\tau,	
\end{align*}
and from \eqref{liebrest}-\eqref{liebrest2}, from the definition \ref{def:LieDer} of Lie derivative and from Lemma \ref{lem:TaylorExp}.
\end{proof}

\begin{lemma}
Let $k_1 \geq 3$, and assume that $\Ham{G} \in \mathcal{P}_{k_1}$ and its vector fields are analytic on $\cB_{s}(R)$, and that $\Ham{h}_0$ satisfies (PER). 
Let $\chi \in \mathcal{P}_{k_1}$ and its vector field be analytic on $\cB_{s}(R)$, and assume that $\chi$ solves \eqref{homeq}. For any $\ell\geq 1$ denote by $\Ham{h}_{0,\ell}$ the functions defined recursively as in \eqref{lieseries} from $\Ham{h}_0$.
Then $\Ham{h}_{0,\ell} \in \mathcal{P}_{(\ell+1)k_1-2\ell}$ and its vector field are analytic on $\cB_{s}(R)$, and
\begin{align}
 \sup_{\zeta \in \cB_{s}(R)} \|X_{ \Ham h_{0,\ell} }(\zeta)\|_{\cW^{s}} &\leq 
2 k_1 (k_1-1)^s \|\Ham{G}\| \,
\left[ 8 (k_1-1) (2k_1-3)^s \,  \|\chi\| \right]^\ell \, R^{\ell+1} \, \times \nonumber \\
&\qquad\qquad \times  \sup_{\zeta \in \cB_s(R)} \max_{n=1,\ldots,(\ell+1)(k_1-1)-2\ell-1} \|\zeta\|_{\cW^s}^{n} . \label{lieseriesh0}
\end{align}
\end{lemma}

\begin{proof}
By using \eqref{homeq} one gets that $\Ham h_{0,1} =\Ham Z- \Ham G$ is analytic on $\cB_{s}(R)$. Then by exploiting \eqref{vfest} and \eqref{EstFG2} one gets the result.
\end{proof}

\begin{lemma} \label{Poissonlemma}
Assume that $\Ham{G} \in \mathcal{P}_{k_1}$ and its vector field are analytic on $\cB_{s}(R)$, and that $\Ham{h}_0$ satisfies PER. 
Let $\chi \in \mathcal{P}_{k_1}$ be the solution of \eqref{homeq}, denote by $\Phi^t_\chi$ 
the flow of the Hamiltonian vector field associated to $\chi$ 
and by $\Phi_\chi$ the corresponding time-one map. Moreover, denote by 
\begin{align*}
\cF(\zeta) & \coloneqq  \Ham{h}_0(\Phi_\chi(\zeta)) - \Ham{h}_0(\zeta) - \{\chi,\Ham{h}_0\}(\zeta).
\end{align*}
Then we have that $\cF$ and its vector field are analytic on $\cB_{s}(R)$, and 
\begin{align} 
\sup_{\zeta \in \cB_{s}(R)} \|X_\cF(\zeta)\|_{\cW^{s}} &\leq \left[ \left( 2k_1(k_1-1)^s \|\chi\| \, R \, \sup_{\zeta \in \cB_{s}(R)}  \max_{n=1,\ldots,k_1-2} \|\zeta\|_{\cW^s}^n \right) + 6 \right] \times \nonumber \\
&\; \; \; \; \times \,  4 (k_1 -1)(2 k_1 -3)^s \, \|\chi\| \, \|\Ham{G}\| \, R \, \sup_{\zeta \in \cB_s(R)} \max_{n=1,\ldots,2k_1-4} \|\zeta\|_{\cW^s}^n   . \label{vfPois}
\end{align}
\end{lemma}

\begin{proof}
Since
\begin{align*}
\Ham{h}_0(\Phi_\chi(\zeta)) - \Ham{h}_0(\zeta) &= \int_0^1 \; \{\chi,\Ham{h}_0 \} \circ \Phi^{t}_{\chi}(\zeta) \; \di t \stackrel{\eqref{homeq}}{=} \int_0^1 \;  \Ham Z( \Phi^{t}_{\chi}(\zeta) ) - \Ham G( \Phi^{t}_{\chi}(\zeta) ) \; \di t,
\end{align*}
if we define $\Ham F(\zeta) \coloneqq \Ham Z(\zeta)-\Ham G(\zeta)$, we get
\begin{align*}
\cF(\zeta) &= \int_0^1 \left(\Ham F( \Phi^t_\chi(\zeta) ) -\Ham F(\zeta)\right) \, \di t.
\end{align*}

Now, we have
\begin{align*}
\sup_{\zeta \in \cB_{s}(R)} \Vert X_{\cF}(\zeta) \Vert_{\cW^{s}} &= \sup_{\zeta \in \cB_{s}(R)} \left\Vert \Omega_1^{-1} \di \left[ \; \int_0^1 \;\left( \Ham F(\Phi_\chi^t(\zeta))- \Ham F(\zeta)\right) \, \di t \; \right] \right\Vert_{\cW^{s}} \\
&\leq  \sup_{\zeta \in \cB_{s}(R)} \left\Vert \int_0^1 \; ( \di \Phi_{\chi}^{-t} ( \Phi_\chi^t(\zeta)) - \mathrm{id}) \Omega_1^{-1} \di \Ham F (\Phi^t_\chi) \; \di t \right\Vert_{\cW^{s}} \\
&\;\;\;\; + \sup_{\zeta \in \cB_{s}(R)} \left\Vert \int_0^1 \;  \left(X_{\Ham F}(\Phi_\chi^t(\zeta))-X_{\Ham F}(\zeta)\right) \; \di t \right\Vert_{\cW^{s}} 
\end{align*}
and by dominated convergence we can bound the last line by 
\begin{align*}
&\sup_{\zeta \in \cB_{s}(R)} \sup_{t \in[0,1]} \Vert \di \Phi_{\chi}^{-t}(\Phi_\chi^t(\zeta))- \mathrm{id} \Vert_{ \cW^{s} \to \cW^{s} } \sup_{\zeta \in \cB_{s}(R)} \Vert X_{\Ham F}(\Phi^t_\chi(\zeta))\Vert_{\cW^{s}}\\
& + \sup_{\zeta \in \cB_{s}(R)} \sup_{t \in [0,1]} \Vert  X_{\Ham F}(\Phi_\chi^t(\zeta))-X_{\Ham F}(\zeta) \Vert_{\cW^{s}} \\
&\leq \sup_{t \in [0,1]} \sup_{\zeta \in \cB_{s}(R)} \Vert \di \Phi_{\chi}^{-t}(\Phi_\chi^t(\zeta))- \mathrm{id} \Vert_{ \cW^{s} \to \cW^{s} } \sup_{\zeta \in\cB_{s}(R)} \Vert X_{\Ham F}(\Phi^t_\chi(\zeta))\Vert_{\cW^{s}}\\
& + \sup_{t \in [0,1]} \sup_{\zeta \in \cB_{s}(R)}  \Vert  X_{\Ham F}(\Phi_\chi^t(\zeta))-X_{\Ham F}(\zeta) \Vert_{\cW^{s}} ,
\end{align*}
where 
\begin{align*}
& \sup_{t \in [0,1]} \sup_{\zeta \in \cB_{s}(R)} \Vert \di \Phi_{\chi}^{-t}(\Phi_\chi^t(\zeta))- \mathrm{id} \Vert_{ \cW^{s} \to \cW^{s} } \sup_{\zeta \in\cB_{s}(R)} \Vert X_{\Ham F}(\Phi^t_\chi(\zeta))\Vert_{\cW^{s}} \\
&\leq \sup_{t \in [0,1]} \left[ \sup_{\zeta \in \cB_{s}(R)} \Vert \di \Phi_{\chi}^{-t}(\Phi_\chi^t(\zeta)) \Vert_{ \cW^{s} \to \cW^{s} }  + 1 \right] \, \sup_{\zeta \in\cB_{s}(R)} \Vert X_{\Ham F}(\Phi^t_\chi(\zeta))\Vert_{\cW^{s}}  \\
&\stackrel{  \eqref{EstFG2}  }{\leq}  \left[ \left( 2k_1(k_1-1)^s \|\chi\| \, \sup_{\zeta \in \cB_{s}(R)} \|\zeta\|_{\cW^s} \, \max_{n=1,\ldots,k_1-2} \|\zeta\|_{\cW^s}^n \right) + 1 \right] \, \times \\
&\qquad\qquad \times \sup_{t \in [0,1]}  \sup_{\zeta \in\cB_{s}(R)} \Vert X_{\Ham F}(\Phi^t_\chi(\zeta))\Vert_{\cW^{s}}  ,
\end{align*}
and where we can exploit Lemma \ref{lem:VFCanTr} so that
\begin{align*}
& \sup_{t \in [0,1]} \sup_{\zeta \in \cB_{s}(R)}  \Vert  X_{\Ham F}(\Phi_\chi^t(\zeta))-X_{\Ham F}(\zeta) \Vert_{\cW^{s}} \\
& \leq 4 (2k_1 -2)(2 k_1 -3)^s \, \|\chi\| \, \|\Ham{F}\| \, R \, \sup_{\zeta \in \cB_s(R)} \max_{n=1,\ldots,2k_1-4} \|\zeta\|_{\cW^s}^n .
\end{align*}

\end{proof}

\begin{lemma} \label{itlemma}
Let $s \geq s_1 \gg 1$, $R>0$, $-\frac{1}{2} < \nu \leq 0$, $k_1 \geq 3$, and consider the Hamiltonian 
\begin{align} \label{H0}
\Ham H^{(0)}(\zeta) &= \Ham{h}_0(\zeta) + \delta \, \Ham F^{(1,0)}(\zeta) + \delta^{1+\nu} \Ham F^{(2,0)}(\zeta)  + \delta^{ 2+\nu } \Ham R^{(0)}(\zeta).
\end{align}
Assume that $\Ham{h}_0$ satisfies (PER) and (INV), that $\Ham F^{(1,0)} \in \cP_{2}$, $\Ham F^{(2,0)} \in \cP_{k_1}$, and that either $\Ham R^{(0)} =0$ or $\Ham R^{(0)} \in \cP_{k_1}$. Moreover, assume that $\Ham F^{(1,0)}$, $\Ham F^{(2,0)}$, $\Ham R^{(0)}$ and their vector fields are analytic on $\cB_s(R)$, with
\begin{align*}
\| \Ham F^{(i,0)} \| \leq F_i, \; \; i=1,2, &\; \; \| \Ham R^{(0)} \| \leq \rho_0 .
\end{align*}

Then there exists $\delta_0=\delta_0(T,F_1,F_2,\rho_0,R,k_1,\nu)>0$ such that if $\delta \leq \delta_0$, then there exists a canonical transformation $\cT^{(0)}_\delta$ analytic on $\cB_{s}(R)$ such that 
\begin{align}
&\sup_{\zeta \in \cB_{s}(R)} \|\cT^{(0)}_\delta(\zeta)-\zeta\|_{\cW^{s}} \nonumber \\
&\leq 4T \, k_1 (k_1-1)^s \, \left( \delta^{1+\nu}  F_2 \, (1+\delta F_1) + \delta F_1 \right) \, R \;  \sup_{\zeta \in \cB_{s}(R)} \max_{n=1,\ldots,k_1-2} \|\zeta\|_{\cW^s}^n ,  \label{CTm}
\end{align}
and $\Ham H^{(1)}  \coloneqq  \Ham H^{(0)} \circ \cT^{(0)}_{\delta}$ has the form 
\begin{align} \label{H1}
\Ham H^{(1)}(\zeta) &= \Ham{h}_0(\zeta) + \delta \Ham Z^{(1,1)}(\zeta) + \delta^{1+\nu} \Ham Z^{(2,1)}(\zeta) + \delta^{ 2(1+\nu) } \Ham R^{(1)}(\zeta) ,
\end{align}
where
\begin{align}
\sup_{\zeta \in \cB_{s}(R)} \|X_{\Ham Z^{(i,1)}}(\zeta)\|_{\cW^{s}} &\leq 2k_1 (k_1-1)^s \; F_j \; R \; \sup_{\zeta \in \cB_{s}(R)} \max_{n=1,\ldots,k_1-2} \|\zeta\|_{\cW^s}^n , \label{Zj1Est} 
\end{align}
for all $\zeta \in \cW^s$, $i=1,2$, and where
\begin{align}
&\sup_{\zeta \in \cB_{s}(R)} \Vert X_{ \Ham R^{(1)}  } (\zeta) \Vert_{\cW^s} \nonumber \\
&\leq 2^4 \, T( k_1 -1)(2 k_1 -3)^s \, R \, \sup_{\zeta \in \cB_s(R)} \max_{n=1,\ldots,2(k_1-2)} \|\zeta\|_{\cW^s}^n \, \times \nonumber \\
&\;\; \times  \bigg\{ 4 \, F_1 \, \left[ \left( \left( 2T k_1(k_1-1)^s F_1 \, R \, \sup_{\zeta \in \cB_{s}(R)}  \max_{n=1,\ldots,k_1-2} \|\zeta\|_{\cW^s}^n \right) + 4 \right) F_1 + F_2  \right] \nonumber \\
&\;\;\;\;  + F_2^2 \left[ \left( 4T k_1(k_1-1)^s F_2 \, R \, \sup_{\zeta \in \cB_{s}(R)}  \max_{n=1,\ldots,k_1-2} \|\zeta\|_{\cW^s}^n \right) + 6 \right] + 2 F_2^2  \bigg\} \nonumber \\
&\;\;\;\; + 2k_1(k_1-1)^s \, ( \rho_0+ 4T \, F_1 F_2 )  \, R \, \sup_{\zeta \in \cB_{s}(R)} \max_{n=1,\ldots,k_1-2} \|\zeta\|_{\cW^s}^n \label{R1Estnu0}
\end{align}
for $\nu=0$, and
\begin{align}
&\sup_{\zeta \in \cB_{s}(R)} \Vert X_{ \Ham R^{(1)}  } (\zeta) \Vert_{\cW^s} \nonumber \\
&\leq \left[ \left( 4T k_1(k_1-1)^s F_2 \, R \, \sup_{\zeta \in \cB_{s}(R)}  \max_{n=1,\ldots,k_1-2} \|\zeta\|_{\cW^s}^n \right) + 8 \right] \times \nonumber \\
&\; \; \; \; \; \; \; \; \times \,  2^4 (k_1 -1)(2 k_1 -3)^s \; T F_2^2 \,  R \, \sup_{\zeta \in \cB_s(R)} \max_{n=1,\ldots,2(k_1-2)} \|\zeta\|_{\cW^s}^n   \label{R1Estnu}
\end{align}
for $-\frac{1}{2}  < \nu < 0$. 
\end{lemma}

\begin{proof}
Here we prove the case $\Ham R^{(0)} \in \cP_{k_1}$, the case $\Ham R^{(0)} =0$ being simpler. The key point of the proof is to look for $\cT^{(0)}_\delta$ as the time-one map of the Hamiltonian vector field of an analytic function. Hence, let $\chi_0 \in \cP_{2}$, and consider the differential equation
\begin{align} \label{chi0}
\dot\zeta &= X_{\delta \chi_0}(\zeta).
\end{align}
By standard theory we have that, if $\|X_{\delta \chi_0}\|_{\cB_{s}(R)}$ is small enough and $\zeta_0 \in \cB_{s}(R)$, then the solution of \eqref{chi0} exists for $|t| \leq 1$. 

Therefore we can define $\cT^t_{0,\delta}:\cB_{s}(R) \to \cB_{s}(R)$, and in particular the corresponding time-one map $\cT^1_{0,\delta}$, which is an analytic canonical transformation, $\delta$-close to the identity. We have 

\begin{align}
& (\cT^1_{0,\delta})^\ast \; (\Ham{h}_0 +  \delta \Ham F^{(1,0)} + \delta^{1+\nu} \Ham F^{(2,0)} + \delta^{ 2+\nu } \Ham R^{(0)}) \nonumber \\
&= \Ham{h}_0 + \delta \left[ \Ham F^{(1,0)} + \{ \chi_0,\Ham{h}_0 \} \right]  + \delta^{1+\nu} \Ham F^{(2,0)} \nonumber \\
&\; \; \; \;  + \left( \Ham{h}_0 \circ \cT^1_{0,\delta} - \Ham{h}_0  - \delta \{ \chi_0,\Ham{h}_0 \} \right) + \delta \left( \Ham F^{(1,0)} \circ \cT^1_{0,\delta} - \Ham F^{(1,0)} \right) \label{nonnorm11} \\
&\; \; \; + \delta^{1+\nu} \left( \Ham F^{(2,0)} \circ \cT^1_{0,\delta} - \Ham F^{(2,0)}  \right) + \delta^{ 2+\nu }  \Ham R^{(0)} + \delta^{ 2+\nu }  \left( \Ham R^{(0)} \circ \cT^1_{0,\delta} - \Ham R^{(0)} \right) \label{nonnorm12} \\
&=: \Ham{h}_0 + \delta \left[ \Ham F^{(1,0)} + \{ \chi_0,\Ham{h}_0 \} \right]  + \delta^{1+\nu} \Ham F^{(2,0)} + \delta^{2+\nu} \Ham R_{1,0}. \nonumber
\end{align}

We can observe that the three terms in square brackets in the second line can be normalized through the choice of a suitable $\chi_0$, while the terms in \eqref{nonnorm11}-\eqref{nonnorm12} contain all terms of order at least $\mathcal{O}(\delta^{2+\nu})$. This lead us to solve the homological equation
\begin{align*}
\{ \chi_0,\Ham{h}_0 \} + \Ham F^{(1,0)} \; &= \; \Ham Z^{(1,1)},
\end{align*}
with $\Ham Z^{(1,1)}$ in normal form.  Lemma \ref{homeqlemma} ensures the existence of $\chi_0$ and $\Ham Z^{(1,1)}$; using the analyticity of the flow $\Phi_{\Ham{h}_0}^t$ ensured by (PER) and \eqref{vfhomeq}, we have that the vector fields associated to $\Ham Z^{(1,1)}$ and $\chi_0$ satisfy

\begin{align} 
\sup_{\zeta \in\cB_{s}(R)} \Vert X_{\chi_0} (\zeta)\Vert_{\cW^{s}} &\leq 2T  \sup_{\zeta \in \cB_{s}(R)} \Vert X_{\Ham F^{(1,0)}}\Vert_{\cW^{s}} \leq 4T \, k_1 (k_1-1)^s \; F_1 \, R , \label{EstimateVectorField11} \\
\sup_{\zeta \in \cB_{s}(R)} \| X_{\Ham Z^{(1,1)}} \|_{\cW^{s}} &\leq   2k_1 (k_1-1)^s \; F_1 \, R  . \label{EstimateVectorField12}
\end{align}

We also estimate the vector field of the remainder $\Ham R_{1,0}$: due to \eqref{vfPois}, \eqref{nonnorm12} and \eqref{liebrest} we have	
\begin{align*}
& \sup_{\zeta \in \cB_{s}(R)} \Vert X_{  \Ham{h}_0 \circ \cT^1_{0,\delta} - \Ham{h}_0  - \delta \{ \chi_0,\Ham{h}_0 \} } (\zeta) \Vert_{\cW^s} \nonumber \\
&\stackrel{\eqref{vfPois}}{\leq} \left[ \left( 4T k_1(k_1-1)^s F_1 \, R \, \sup_{\zeta \in \cB_{s}(R)}  \max_{n=1,\ldots,k_1-2} \|\zeta\|_{\cW^s}^n \right) + 6 \right] \times \nonumber \\
&\; \; \; \; \times \,  4 (k_1 -1)(2 k_1 -3)^s \; 2T \delta^2 F_1^2 \, R , 
\end{align*}

\begin{align*}
\sup_{\zeta \in \cB_{s}(R)} \Vert X_{ \Ham F^{(1,0)} \circ \cT_{0,\delta}^1 - \Ham F^{(1,0)} } (\zeta) \Vert_{\cW^s} &\stackrel{\eqref{liebrest}}{\leq} 8 ( k_1 -1)(2 k_1 -3)^s \, 2T \delta F_1^2  \, R ,
\end{align*}
\begin{align*}
&\sup_{\zeta \in \cB_{s}(R)} \Vert X_{ \Ham F^{(2,0)} \circ \cT_{0,\delta}^1 - \Ham F^{(2,0)} } (\zeta) \Vert_{\cW^s} \\
&\stackrel{\eqref{liebrest}}{\leq} 8 ( k_1 -1)(2 k_1 -3)^s \, 2T \delta F_1 \, F_2  \, R \, \sup_{\zeta \in \cB_s(R)} \max_{n=1,\ldots,2(k_1-2)} \|\zeta\|_{\cW^s}^n , \\
&\sup_{\zeta \in \cB_{s}(R)} \Vert X_{ \Ham R^{(0)} \circ \cT_{0,\delta}^1 - \Ham R^{(0)} } (\zeta) \Vert_{\cW^s} \\
&\stackrel{\eqref{liebrest}}{\leq} 8 ( k_1 -1)(2 k_1 -3)^s \, 2T \delta F_1 \, \rho_0  \, R \, \sup_{\zeta \in \cB_s(R)} \max_{n=1,\ldots,2(k_1-2)} \|\zeta\|_{\cW^s}^n ,
\end{align*}
so that we can deduce from \eqref{nonnorm11}-\eqref{nonnorm12} that
\begin{align} 
& \sup_{\zeta \in \cB_{s}(R)} \Vert X_{ \Ham R_{1,0}  } (\zeta) \Vert_{\cW^s} \nonumber \\
&\leq 16T ( k_1 -1)(2 k_1 -3)^s \, F_1 \,  R \, \sup_{\zeta \in \cB_s(R)} \max_{n=1,\ldots,2(k_1-2)} \|\zeta\|_{\cW^s}^n  \times \nonumber \\
&\;\;\;\; \times \left[ \delta^{-\nu} \; \left( \left( 2T k_1(k_1-1)^s F_1 \, R \, \sup_{\zeta \in \cB_{s}(R)}  \max_{n=1,\ldots,k_1-2} \|\zeta\|_{\cW^s}^n \right) + 4 \right) F_1 + F_2 + \delta \; \rho_0 \right] \nonumber \\
&\;\;\;\; + 2k_1(k_1-1)^s \; \rho_0  \, R \, \sup_{\zeta \in \cB_{s}(R)} \max_{n=1,\ldots,k_1-2} \|\zeta\|_{\cW^s}^n . \label{R10Est}
\end{align}

Now let $\chi_1 \in \cP_{k_1}$, and consider the differential equation
\begin{align} \label{chi1}
\dot\zeta &= X_{\delta^{1+\nu} \chi_1}(\zeta).
\end{align}
Again by standard theory we have that, if $\|X_{\delta^{1+\nu} \chi_1}\|_{\cB_{s}(R)}$ is small enough and $\zeta_0 \in \cB_{s}(R)$, then the solution of \eqref{chi1} exists for $|t| \leq 1$. 

Therefore we can define $\cT^t_{1,\delta}:\cB_{s}(R) \to \cB_{s}(R)$, and in particular the corresponding time-one map $\cT^1_{1,\delta}$, which is an analytic canonical transformation, $\delta^{1+\nu}$-close to the identity. We have 
\begin{align}
& (\cT^1_{1,\delta})^\ast \; (\Ham{h}_0 +  \delta \Ham Z^{(1,1)} + \delta^{1+\nu} \Ham F^{(2,0)} + \delta^{ 2+ \nu } \Ham R_{1,0}) \nonumber \\
&= \Ham{h}_0 + \delta \Ham Z^{(1,1)} + \delta^{1+\nu} \left[ \Ham F^{(2,0)} + \{ \chi_1,\Ham{h}_0 \} \right]   \nonumber \\
&\; \; \; \; + \delta^{2+\nu}  \Ham R_{1,0} + \left( \Ham{h}_0 \circ \cT^1_{1,\delta} - \Ham{h}_0  - \delta^{1+\nu} \{ \chi_1,\Ham{h}_0 \} \right) + \delta \left( \Ham Z^{(1,1)} \circ \cT^1_{1,\delta} - \Ham Z^{(1,1)} \right) \label{nonnorm21} \\
&\; \; \; + \delta^{1+\nu} \left( \Ham F^{(2,0)} \circ \cT^1_{1,\delta} - \Ham F^{(2,0)} \right) + \delta^{2+\nu} \left( \Ham R_{1,0} \circ \cT^1_{1,\delta} - \Ham R_{1,0} \right) \label{nonnorm22} \\
&=: \Ham{h}_0 + \delta \Ham Z^{(1,1)} + \delta^{1+\nu} \left[ \Ham F^{(2,0)} + \{ \chi_1,\Ham{h}_0 \} \right] + \delta^{2(1+\nu)} \Ham R^{(1)}. \nonumber
\end{align}

Lemma \ref{homeqlemma} ensures the existence of $\chi_1$ and $\Ham Z^{(2,1)}$; using the analyticity of the flow $\Phi_{\Ham{h}_0}^t$ ensured by (PER) and \eqref{vfhomeq}, we have that the vector fields associated to $\Ham Z^{(2,1)}$ and $\chi_1$ satisfy
\begin{align} 
& \sup_{\zeta \in\cB_{s}(R)} \Vert X_{\chi_1} (\zeta)\Vert_{\cW^{s}} \nonumber \\
&\leq 2T  \sup_{\zeta \in \cB_{s}(R_m)} \Vert X_{\Ham F^{(2,0)}}\Vert_{\cW^{s}} \leq 4T \, k_1 (k_1-1)^s \; F_2 \, R \; \sup_{\zeta \in \cB_{s}(R)} \max_{n=1,\ldots,k_1-2} \|\zeta\|_{\cW^s}^n , \label{EstimateVectorField21} \\
& \sup_{\zeta \in \cB_{s}(R)} \| X_{\Ham Z^{(2,1)}} \|_{\cW^{s}} \leq   2 \, k_1 (k_1-1)^s \; F_2 \, R \; \sup_{\zeta \in \cB_{s}(R)} \max_{n=1,\ldots,k_1-2} \|\zeta\|_{\cW^s}^n . \label{EstimateVectorField22}
\end{align}

Defining now $\cT^{(0)}_\delta(\zeta)  \coloneqq  \Phi^1_{\delta^{ 1+\nu } \chi_1} \circ \Phi^1_{\delta \, \chi_0}(\zeta)$ we can apply Lemma \ref{lem:VFCanTr}, \eqref{EstimateVectorField11} and \eqref{EstimateVectorField21}, so that we obtain
\begin{align*}
 \sup_{\zeta \in \cB_{s}(R)} \Vert \cT^{(0)}_\delta(\zeta)-\zeta \Vert_{\cW^{s}} &\leq \sup_{\zeta \in \cB_{s}(R)}  \Vert \Phi^1_{\delta^{1+\nu} \chi_1} \left( \Phi^1_{\delta \, \chi_0}(\zeta) \right) - \Phi^1_{\delta \, \chi_0}(\zeta) \Vert_{\cW^{s}} +  \Vert \Phi^1_{\delta \, \chi_0}(\zeta) - \zeta \Vert_{\cW^{s}} \\
&\leq \sup_{\zeta \in \cB_{s}(R)} \Vert X_{\delta^{1+\nu} \chi_1} \circ \Phi^1_{\delta \, \chi_0} \Vert_{\cW^{s}} + \sup_{\zeta \in \cB_{s}(R)} \Vert X_{\delta \chi_0}  \Vert_{\cW^{s}} \\
&\leq 4T \, k_1 (k_1-1)^s \, \left( \delta^{1+\nu}  F_2 \, (1+\delta F_1) + \delta F_1 \right) \, R \; \sup_{\zeta \in \cB_{s}(R)} \max_{n=1,\ldots,k_1-2} \|\zeta\|_{\cW^s}^n ,
\end{align*}
which leads to \eqref{CTm}.

Now we estimate the vector fields of the terms in \eqref{nonnorm21}-\eqref{nonnorm22}: we have that
\begin{align*}
&\sup_{\zeta \in \cB_{s}(R)} \Vert X_{  \Ham{h}_0 \circ \cT^1_{1,\delta} - \Ham{h}_0  - \delta^{1+\nu} \{ \chi_1,\Ham{h}_0 \} } (\zeta) \Vert_{\cW^s} \\
&\stackrel{\eqref{vfPois}}{\leq} \left[ \left( 4T k_1(k_1-1)^s F_2 \, R \, \sup_{\zeta \in \cB_{s}(R)}  \max_{n=1,\ldots,k_1-2} \|\zeta\|_{\cW^s}^n \right) + 6 \right] \times \nonumber \\
&\; \; \; \; \times \,  4 (k_1 -1)(2 k_1 -3)^s \; 2T \delta^{2(1+\nu)} F_2^2 \, R \, \sup_{\zeta \in \cB_s(R)} \max_{n=1,\ldots,2k_1-4} \|\zeta\|_{\cW^s}^n   ,
\end{align*}

\begin{align*}
& \sup_{\zeta \in \cB_{s}(R)} \Vert X_{ \Ham Z^{(1,1)} \circ \cT_{1,\delta}^1 - \Ham Z^{(1,1)} } (\zeta) \Vert_{\cW^s} \\ &\stackrel{\eqref{EstimateVectorField12}}{\leq} 2k_1 (k_1-1)^s \; 2T \delta^{1+\nu} F_1 F_2 \, R \; \sup_{\zeta \in \cB_{s}(R)} \max_{n=1,\ldots,k_1-2} \|\zeta\|_{\cW^s}^n ,
\end{align*}
\begin{align*}
& \sup_{\zeta \in \cB_{s}(R)} \Vert X_{ \Ham F^{(2,0)} \circ \cT_{1,\delta}^1 - \Ham F^{(2,0)} } (\zeta) \Vert_{\cW^s} \\ &\stackrel{\eqref{liebrest}}{\leq} 8 ( k_1 -1)(2 k_1 -3)^s \, 2T \delta^{1+\nu} F_2^2  \, R \, \sup_{\zeta \in \cB_s(R)} \max_{n=1,\ldots,2(k_1-2)} \|\zeta\|_{\cW^s}^n ,
\end{align*}

\begin{align*}
& \sup_{\zeta \in \cB_{s}(R)} \Vert X_{ \Ham R_{1,0} \circ \cT_{1,\delta}^1 - \Ham R_{1,0} } (\zeta) \Vert_{\cW^s} \\
&\leq \delta^{1+\nu} \, \bigg\{ 32T^2 ( k_1 -1)(2 k_1 -3)^s \, F_1 \,  R \, \sup_{\zeta \in \cB_s(R)} \max_{n=1,\ldots,2(k_1-2)} \|\zeta\|_{\cW^s}^n  \times \nonumber \\
&\;\;\;\; \times \left[ \delta^{-\nu} \; \left( \left( 2T k_1(k_1-1)^s F_1 \, R \, \sup_{\zeta \in \cB_{s}(R)}  \max_{n=1,\ldots,k_1-2} \|\zeta\|_{\cW^s}^n \right) + 4 \right) F_1 +  F_2 + \delta \; \rho_0 \right] \nonumber \\
&\;\;\;\; + 4T k_1(k_1-1)^s \; \rho_0  \, R \,\sup_{\zeta \in \cB_{s}(R)} \max_{n=1,\ldots,k-2} \|\zeta\|_{\cW^s}^n \bigg\} ,
\end{align*}
from which we deduce that
\begin{align}
& \sup_{\zeta \in \cB_{s}(R)} \Vert X_{ \Ham R^{(1)}  } (\zeta) \Vert_{\cW^s} \nonumber \\
&\leq 16T ( k_1 -1)(2 k_1 -3)^s \, F_1 \,  R \, \left( \sup_{\zeta \in \cB_s(R)} \max_{n=1,\ldots,2(k_1-2)} \|\zeta\|_{\cW^s}^n \right) \, \delta^{-\nu} \, \times \nonumber \\
&\;\;\;\; \; \; \; \; \times  \left[ \delta^{-\nu} \; \left( \left( 2T k_1(k_1-1)^s F_1 \, R \, \sup_{\zeta \in \cB_{s}(R)}  \max_{n=1,\ldots,k_1-2} \|\zeta\|_{\cW^s}^n \right) + 4 \right) F_1 + F_2 + \delta \; \rho_0 \right] \nonumber \\
&\;\;\;\; + 2k_1(k_1-1)^s \; \delta^{-\nu}  \; \rho_0  \, R \, \sup_{\zeta \in \cB_{s}(R)} \max_{n=1,\ldots,k-2} \|\zeta\|_{\cW^s}^n \nonumber \\
&\; \; \; \; +\left[ \left( 4T k_1(k_1-1)^s F_2 \, R \, \sup_{\zeta \in \cB_{s}(R)}  \max_{n=1,\ldots,k_1-2} \|\zeta\|_{\cW^s}^n \right) + 6 \right] \times \nonumber \\
&\; \; \; \; \; \; \; \; \times \,  4 (k_1 -1)(2 k_1 -3)^s \; 2T F_2^2 \,  R \, \sup_{\zeta \in \cB_s(R)} \max_{n=1,\ldots,2(k_1-2)} \|\zeta\|_{\cW^s}^n   \nonumber \\
&\; \; \; \; + 2k_1 (k_1-1)^s \; 2T \delta^{-\nu} F_1 F_2 \, R \; \sup_{\zeta \in \cB_{s}(R)} \max_{n=1,\ldots,k_1-2} \|\zeta\|_{\cW^s}^n \nonumber \\
&\; \; \; \; + 8 ( k_1 -1)(2 k_1 -3)^s \, 2T \, F_2^2  \,  R \, \sup_{\zeta \in \cB_s(R)} \max_{n=1,\ldots,2(k_1-2)} \|\zeta\|_{\cW^s}^n \nonumber \\
&\; \; \; \; +32T^2 ( k_1 -1)(2 k_1 -3)^s \, F_1 \,  R \, \sup_{\zeta \in \cB_s(R)} \max_{n=1,\ldots,2(k_1-2)} \|\zeta\|_{\cW^s}^n  \times \nonumber \\
&\;\;\;\; \; \; \; \; \times  \delta \, \left[ \delta^{-\nu} \; \left( \left( 2T k_1(k_1-1)^s F_1 \, R \, \sup_{\zeta \in \cB_{s}(R)}  \max_{n=1,\ldots,k_1-2} \|\zeta\|_{\cW^s}^n \right) + 4 \right) F_1 +  F_2 + \delta \; \rho_0 \right] \nonumber \\
&\;\;\;\; + 4T \, \delta \, k_1(k_1-1)^s \;  \rho_0  \, R \, \sup_{\zeta \in \cB_{s}(R)} \max_{n=1,\ldots,k_1-2} \|\zeta\|_{\cW^s}^n . \label{R1Esttmp}
\end{align}
In particular, if we evaluate \eqref{R1Esttmp} for $\nu=0$ we have
\begin{align*}
& \sup_{\zeta \in \cB_{s}(R)} \Vert X_{ \Ham R^{(1)}  } (\zeta) \Vert_{\cW^s} \nonumber \\
&\leq 8T( k_1 -1)(2 k_1 -3)^s \, R \, \sup_{\zeta \in \cB_s(R)} \max_{n=1,\ldots,2(k_1-2)} \|\zeta\|_{\cW^s}^n \, \times \nonumber \\
&\;\; \times  \bigg\{2 \, F_1 \, \left[ \left( \left( 2T k_1(k_1-1)^s F_1 \, R \, \sup_{\zeta \in \cB_{s}(R)}  \max_{n=1,\ldots,k_1-2} \|\zeta\|_{\cW^s}^n \right) + 4 \right) F_1 + F_2 + \delta \; \rho_0 \right] \nonumber \\
&\;\;\;\;  + F_2^2 \left[ \left( 4T k_1(k_1-1)^s F_2 \, R \, \sup_{\zeta \in \cB_{s}(R)}  \max_{n=1,\ldots,k_1-2} \|\zeta\|_{\cW^s}^n \right) + 6 \right] \nonumber \\
&\;\;\;\; + 2 F_2^2 +4T \, F_1 \, \delta \, \left[ \left( \left( 2T k_1(k_1-1)^s F_1 \, R \, \sup_{\zeta \in \cB_{s}(R)}  \max_{n=1,\ldots,k_1-2} \|\zeta\|_{\cW^s}^n \right) + 4 \right) F_1 +  F_2 + \delta \; \rho_0 \right] \bigg\} \nonumber \\
&\;\;\;\; + 2k_1(k_1-1)^s \, ( \rho_0+2T \, ( \delta \, \rho_0+F_1 F_2) )  \, R \, \sup_{\zeta \in \cB_{s}(R)} \max_{n=1,\ldots,k_1-2} \|\zeta\|_{\cW^s}^n ,
\end{align*}
so that there exists $\delta_0=\delta_0(T,F_1,F_2,\rho_0,R,k_1)>0$ such that if $\delta \leq \delta_0$
\begin{align*}
\sup_{\zeta \in \cB_{s}(R)} \Vert X_{ \Ham R^{(1)}  } (\zeta) \Vert_{\cW^s} &\leq 2^4 \, T( k_1 -1)(2 k_1 -3)^s \, R \, \sup_{\zeta \in \cB_s(R)} \max_{n=1,\ldots,2(k_1-2)} \|\zeta\|_{\cW^s}^n \, \times \nonumber \\
&\;\; \times  \bigg\{ 4 \, F_1 \, \left[ \left( \left( 2T k_1(k_1-1)^s F_1 \, R \, \sup_{\zeta \in \cB_{s}(R)}  \max_{n=1,\ldots,k_1-2} \|\zeta\|_{\cW^s}^n \right) + 4 \right) F_1 + F_2  \right] \nonumber \\
&\;\;\;\;  + F_2^2 \left[ \left( 4T k_1(k_1-1)^s F_2 \, R \, \sup_{\zeta \in \cB_{s}(R)}  \max_{n=1,\ldots,k_1-2} \|\zeta\|_{\cW^s}^n \right) + 6 \right] + 2 F_2^2  \bigg\} \nonumber \\
&\;\;\;\; + 2k_1(k_1-1)^s \, ( \rho_0+ 4T \, F_1 F_2 )  \, R \, \sup_{\zeta \in \cB_{s}(R)} \max_{n=1,\ldots,k_1-2} \|\zeta\|_{\cW^s}^n . \nonumber 
\end{align*}
Similarly, if $\nu \in \left( - \frac{1}{2},0 \right)$, then there exists $\delta_0=\delta_0(T,F_1,F_2,\rho_0,R,k_1,\nu)>0$ such that if $\delta \leq \delta_0$
\begin{align*}
&\sup_{\zeta \in \cB_{s}(R)} \Vert X_{ \Ham R^{(1)}  } (\zeta) \Vert_{\cW^s} \nonumber \\
&\leq \left[ \left( 4T k_1(k_1-1)^s F_2 \, R \, \sup_{\zeta \in \cB_{s}(R)}  \max_{n=1,\ldots,k_1-2} \|\zeta\|_{\cW^s}^n \right) + 8 \right] \times \nonumber \\
&\; \; \; \; \; \; \; \; \times \,  2^4 (k_1 -1)(2 k_1 -3)^s \; T F_2^2 \,  R \, \sup_{\zeta \in \cB_s(R)} \max_{n=1,\ldots,2(k_1-2)} \|\zeta\|_{\cW^s}^n   .  \nonumber 
\end{align*}

\end{proof}

\begin{proof}[Proof of Lemma \ref{NFest}]
The Hamiltonian \eqref{truncsys} satisfies the assumptions of Lemma \ref{itlemma}, $\Ham F_{j1} \circ \Pi_M$ in place of $\Ham F_{(i,0)}$, ($i=1,2$), $F_i = \| \Ham F_{i1} \| \, M^2$ ($i=1,2$), $\Ham R^{(0)}=0$, $\rho_0=0$ and $k_1=k$.

So by applying Lemma \ref{itlemma}, we obtain that for $\delta \leq \delta_0$ (where $\delta_0$ is the quantity appearing in the statement of Lemma \ref{itlemma}) there exists an analytic canonical transformation $\cT^{(0)}_{\delta,M}: \cB_{s}(R) \to \cB_{s}(R)$ such that 
\begin{align}
& \Ham H_{1,M} \circ \cT^{(0)}_{\delta,M} = \Ham{h}_0 + \delta \Ham Z^{(1)}_M + \delta^{1+\nu} \Ham Z^{(2)}_M + \delta^{2(1+\nu)} \cR^{(1)}_M, \label{step1} \\
& \Ham Z^{(i)}_M  \coloneqq  \la \Ham F_{i1} \circ \Pi_M \ra, \; \; i=1,2, \; \; \delta^{2(1+\nu)} \cR^{(1)}_M \nonumber = \delta^{2(1+\nu)} \Ham R^{(1)},
\end{align}
where $\Ham R^{(1)}$ is the remainder appearing in the statement of Lemma \ref{itlemma}. 

\end{proof}

\section{Proof of Proposition \ref{NLSphiProp}  } \label{ApprEstSec11}

In order to prove Proposition \ref{NLSphiProp} we first discuss the specific energies associated to the high modes, and then the ones associated to the low modes. In the rest of the section $K = (K_1,\ldots,K_d) \in \mathbb{Z}^d_N$.

First we remark that for $\kappa = (h K_1,\ldots,h K_d)$ we have
\begin{align}
\left| \omega^2_K \right| &\stackrel{\eqref{FreqNormModeKG}}{=}  1+ 4 \sum_{j=1}^d \sin^2\left(\frac{h K_j \pi}{2}\right)  \;  \leq \; \pi^2 (1 + h^2 |K|^2 ) . \label{EstFreqNLSr}
\end{align}
moreover, for $K \neq 0$ and for all $m \geq s > 0$
\begin{align}
\frac{|\hat{p}_K|^2+\pi^2(1+h^2 |K|^2)|\hat{q}_K|^2}{2} &\leq \pi^2 \, |K|^{-2s} \; \frac{|\hat{p}_K|^2+(1+h^2 |K|^2 )|\hat{q}_K|^2}{2} |K|^{2s} \nonumber \\
&\leq \pi^2 \, |K|^{-2s} \; \left(1+h^2 |K|^2  \right) \; \|(\xi,\eta)\|_{\cW^s}^2 , \label{NormModeNLSEst1} 
\end{align}
where in the last inequality we used the fact that $\ell^1_s \subset \ell^2_s$ (see Proposition III.2 in \cite{faou2012geometric}). \\

In order to estimate $\cE_\kappa$ for large values of $K_j$ ($j=1,\ldots,d$), it is convenient to divide the frequency-space in different regions and bound the terms supported in each region separately. Indeed, in the multi-dimensional case the introduction of different regions in the frequency space will help us estimating most of the terms in an efficient way. Let us define
\begin{align} 
\mathscr{L} &=  \mathscr{L}_{h,\delta,\alpha,s}  \coloneqq \left\{ L=(L_1,\ldots,L_d) \in \mathbb{Z}^d : h L_j \in 2\mathbb{Z} \; \forall j=1,\ldots,d, |K_1|+\cdots+|K_d| > M  \right\} , \label{eq:scrL} \\
M  &= M(h,\delta) = h^{ - (1-\delta) }  , \; \; \delta>0, \; s>d/2 . \nonumber 
\end{align}

Hence we obtain that for $\kappa = (h K_1,\ldots,h K_d)$  and $|K_1|+\cdots+|K_d| > M$
\begin{align}
\frac{ \cE_\kappa }{h^{2\alpha}} &= \frac{1}{2} \sum_{ \substack{ L \in \mathscr{L} \\ |K+L|=0 }  } \left( |\hat{p}_{0}|^2 + \omega_{K}^2 \left| \hat{q}_{0} \right|^2 \right)  + \frac{1}{2} \sum_{  \substack{ L \in \mathscr{L} \\ |K+L| \neq 0 } } \left( |\hat{p}_{K+L}|^2 + \omega_K^2 \left| \hat{q}_{K+L} \right|^2 \right)  \nonumber \\
&= \frac{1}{2} \sum_{  \substack{ L \in \mathscr{L} \\ |K+L| \neq 0 } } \left( |\hat{p}_{K+L}|^2 + \omega_K^2 \left| \hat{q}_{K+L} \right|^2 \right)  , \nonumber 
\end{align}
where
\begin{align}
& \sum_{  \substack{ L \in \mathscr{L} \\ |K+L| \neq 0 } } \left( |\hat{p}_{K+L}|^2 + \omega_K^2 \left| \hat{q}_{K+L} \right|^2 \right)  \nonumber \\
&\stackrel{ \eqref{NormModeNLSEst1} }{\leq} \pi^2 \; \|(\xi,\eta)\|_{\cW^{s}}^2 \; 2 \sum_{\substack{ L \in \mathscr{L} \\ |K+L| \neq 0 } } |K+L|^{-2s}  \; \left[ 1 + h^2 \; \sum_{j=1}^d (K_j+L_j)^2  \right] \nonumber \\
&\leq 2\pi^2 \; \|(\xi,\eta)\|_{\cW^{s}}^2 \bigg[  \sum_{ \substack{ L \in \mathscr{L} \\ |K+L| \neq 0 } } |K+L|^{-2s} + h^2 \; \sum_{ \substack{ L \in \mathscr{L} \\ |K+L| \neq 0 } } |K+L|^{-2s} \; \sum_{j=1}^d  (K_j+L_j)^2  \bigg] . \label{HighFreqTerm}
\end{align}

Now, 
\begin{align}
&\sum_{ \substack{ L \in \mathscr{L} \\ |K+L| \neq 0  } } |K+L|^{-2s} \nonumber \\
&\leq |K|^{-2s} \nonumber \\
&\;\;\;\; + \sum_{ \substack{  L=(L_1,\ldots,L_d) \in \mathscr{L} \\ |K+L| \neq 0 \\ L_1=0, (L_2,\ldots,L_d) \neq (0,\ldots,0) } } |K+L|^{-2s} +\cdots+ \sum_{ \substack{  L=(L_1,\ldots,L_d) \in \mathscr{L} \\ |K+L| \neq 0 \\ (L_1,\ldots,L_{d-1}) \neq (0,\ldots,0), L_d=0 } } |K+L|^{-2s} \nonumber \\
&\;\;\;\; + \sum_{ \substack{  L=(L_1,\ldots,L_d) \in \mathscr{L} \\ |K+L| \neq 0 \\ L_1,\ldots,L_d \neq 0 } } |K+L|^{-2s} . \label{HighFreqTerm1}
\end{align}

We now estimate the last sum in \eqref{HighFreqTerm1}; we point out that for $L_1,\ldots,L_d \neq 0$ we have $|L| \geq 2d \; h^{-1}$, hence $d \; |K| \leq |L|$.

Therefore, for any $s>\frac{d}{2}$ and for $\kappa = (h K_1,\ldots,h K_d)$ and $|K_1|+\cdots+|K_d| > M$
\begin{align}
\sum_{ \substack{  L=(L_1,\ldots,L_d) \in \mathscr{L} \\ |K+L| \neq 0 \\ L_1,\ldots,L_d \neq 0  } } |K+L|^{-2s} &\leq \sum_{ \substack{  L=(L_1,\ldots,L_d) \in \mathscr{L} \\ |K+L| \neq 0 \\ L_1,\ldots,L_d \neq 0 } }  | \, |K|-|L| \, |^{-2s} \nonumber \\
&\leq \sum_{ \substack{  L=(L_1,\ldots,L_d) \in \mathscr{L} \\ |K+L| \neq 0 \\ L_1,\ldots,L_d \neq 0  } } |K|^{-2s} \left(  |L|/|K| - 1 \right)^{-2s} , \nonumber 
\end{align}
where for $d=2,3$ the above sum can be bounded by
\begin{align}
 |K|^{-2s} \, 2\pi \, \pi^{d-2} \, \int_{d|K|}^{+\infty} R^{d-1} (R/|K|-1)^{-2s} \di R &= 2\pi^{d-1} \; C(s,d) \; |K|^{-(2s-d)} , \label{Est1HighFreqTerm1}
\end{align}
where
\begin{align*}
C(s,2) = (4s-3) \, \frac{ \Gamma(2s-2) }{ \Gamma(2s) }, &\; \; C(s,3) = \frac{ (6s-11) 3s+13 }{4^{s-1}} \, \frac{ \Gamma(2s-3) }{ \Gamma(2s) } , \\
\end{align*}
($\Gamma$ is Euler's gamma function), while for $d=1$ the above sum can be bounded by
\begin{align}
&\sum_{ \substack{  L=(L_1,\ldots,L_d) \in \mathscr{L} \\ |K+L| \neq 0 \\ L_1,\ldots,L_d \neq 0  } } |K+L|^{-2s} \nonumber \\
&\leq  \sum_{ \substack{  L=(L_1,\ldots,L_d) \in \mathscr{L} \\ |K+L| \neq 0 \\ L_1,\ldots,L_d \neq 0  \\ |K| \leq |L| \leq 2|K| } } |K+L|^{-2s}  + 2 |K|^{-2s} \,  \int_{2|K|}^{+\infty}  (R/|K|-1)^{-2s} \di R  , \label{Est1HighFreqTerm11}
\end{align}
where the above integral is equal to $\frac{1}{2s-1} |K|$, and where the sum in the right-hand side of \eqref{Est1HighFreqTerm11} is bounded by
\begin{align*}
|K|^{-2s} \sum_{ \substack{  L=(L_1,\ldots,L_d) \in \mathscr{L} \\ |K+L| \neq 0 \\ L_1,\ldots,L_d \neq 0  \\ |K| \leq |L| \leq 2|K| } } |1+L/K|^{-2s} &\leq C_s \; |K|^{-2s+1} ,
\end{align*}
for some $C_s>0$ (as $|K+L| \neq 0$ implies $|K+L| \geq 1$). 

Next we estimate the second sum in \eqref{HighFreqTerm1}; of course the sum vanishes for $d=1$, while for $d=2,3$
\begin{align}
\sum_{ \substack{  L=(L_1,\ldots,L_d) \in \mathscr{L} \\ |K+L| \neq 0 \\ L_1=0, (L_2,\ldots,L_d) \neq (0,\ldots,0) } } |K+L|^{-2s} &\leq |K|^{-2s} \sum_{ \substack{  L=(0,L_2,\ldots,L_d) \in \mathscr{L} \\ |K+L| \neq 0 \\ (L_2,\ldots,L_d) \neq (0,\ldots,0) }  } \left| \frac{ 2|L| }{ h } - 1 \right|^{-2s}  , \label{Est2HighFreqTerm1}
\end{align}
which is polynomially small with respect to $h$, and similarly for the other sums in the second line of \eqref{HighFreqTerm1}. 

Now we want to estimate the term
\begin{equation*}
\sum_{ \substack{ L=(L_1,\ldots,L_d) \in \mathscr{L} \\ |K+L| \neq 0 } } |K+L|^{-2s} \sum_{j+1}^d (K_j+L_j)^2 
\end{equation*}
in \eqref{HighFreqTerm}. We have for all $j \in \{1,\ldots,d\}$ 
\begin{align*}
& \sum_{ \substack{ L=(L_1,\ldots,L_d) \in \mathscr{L} \\ |K+L| \neq 0 } } |K+L|^{-2s} \, (K_j+L_j)^2 \nonumber \\
&\leq 2 \sum_{ \substack{ L=(L_1,\ldots,L_d) \in \mathscr{L} \\ |K+L| \neq 0 } } |K+L|^{-2s} \, (K_j^2+L_j^2) \nonumber \\
&\leq K_j^2 \; \sum_{ \substack{ L=(L_1,\ldots,L_d) \in \mathscr{L} \\ |K+L| \neq 0 } } |K+L|^{-2s} \, + \sum_{ \substack{ L=(L_1,\ldots,L_d) \in \mathscr{L} \\ |K+L| \neq 0 } } |K+L|^{-2s} \, L_j^2 ,
\end{align*}
where the first sum can be bounded as before, while the second sum can be decomposed in the following way,
\begin{align*}
&  \sum_{ \substack{ L=(L_1,\ldots,L_d) \in \mathscr{L} \\ |K+L| \neq 0 \\ L_1=0,(L_2,\ldots,L_d) \neq (0,\ldots,0) } } |K+L|^{-2s} \, L_j^2 + \cdots + \sum_{ \substack{ L=(L_1,\ldots,L_d) \in \mathscr{L} \\ |K+L| \neq 0 \\ L_d=0,(L_1,\ldots,L_{d-1}) \neq (0,\ldots,0) } } |K+L|^{-2s} \, L_j^2 \nonumber \\
& \; \; \; \; + \sum_{ \substack{ L=(L_1,\ldots,L_d) \in \mathscr{L} \\ |K+L| \neq 0 \\ L_1,\ldots,L_d \neq 0 } } |K+L|^{-2s} \, L_j^2 .
\end{align*}

Then, for any $s>\frac{d}{2}+1$ and for $\kappa = (h K_1,\ldots, h K_d)$ and $|K_1|+\ldots+|K_d| > M$ we have
\begin{align}
\sum_{ \substack{ L=(L_1,\ldots,L_d) \in \mathscr{L} \\ |K+L| \neq 0 \\ L_1, \ldots,L_d \neq 0 } }  |K+L|^{-2s} L_j^2 &\leq \sum_{ \substack{ L=(L_1,\ldots,L_d) \in \mathscr{L} \\ |K+L| \neq 0 \\ L_1, \ldots,L_d \neq 0 } }  |K+L|^{-2s} |L|^2 \nonumber 
\end{align}
which for $d=2,3$ can be bounded by
\begin{align}
 |K|^{-2s} \sum_{ \substack{ L=(L_1,\ldots,L_d) \in \mathscr{L} \\ |K+L| \neq 0 \\ L_1, \ldots,L_d \neq 0 } } \left| \frac{|L|}{|K|} - 1  \right|^{-2s} \, |L|^2 &\leq 2\pi^{d-1} \; |K|^{-2s} \; \int_{d |K|}^{+\infty} R^{d+1} (R/|K| -1)^{-2s}  \; \mathrm{d}R \nonumber \\
&= 2\pi^{d-1} \; C_1(s,d) \; |K|^{-(2s-d-2)} , \label{HighFreqTerm2}
\end{align}

\begin{align*}
C_1(s,2) &= 2( 4s(2s-5)+9 ) (4s-5) \frac{ \Gamma(2s-4) }{ \Gamma(2s) } , \\
C_1(s,3) &= (-1)^{2s+1} \left[ B_3(5,1-2s) + 24 \frac{ \Gamma(2s-5) }{ \Gamma(2s) } \right] , \; \; B_z(a,b)  \coloneqq \int_0^z t^{a-1} (1-t)^{b-1} \mathrm{d}t, 
\end{align*}
($(z,a,b) \mapsto B_z(a,b)$ is the generalized Euler Beta function); for $d=1$ the above sum can be bounded by
\begin{align*}
\sum_{ \substack{ L=(L_1,\ldots,L_d) \in \mathscr{L} \\ |K+L| \neq 0 \\ L_1, \ldots,L_d \neq 0 } } L_j^2 \, |K+L|^{-2s} + 2|K|^{-2s} \; \int_{2|K|}^{+\infty} R^2 \left( \frac{R}{|K|} -1 \right)^{-2s} \; \mathrm{d}R &\leq C_s \; |K|^{-(2s-3)} .
\end{align*}
Similarly, we now estimate the terms
\begin{equation*}
\sum_{ \substack{ L=(L_1,\ldots,L_d) \in \mathscr{L} \\ |K+L| \neq 0 \\ L_1=0,(L_2,\ldots,L_d) \neq (0,\ldots,0) } } |K+L|^{-2s} \, L_j^2 + \cdots + \sum_{ \substack{ L=(L_1,\ldots,L_d) \in \mathscr{L} \\ |K+L| \neq 0 \\ L_d=0,(L_1,\ldots,L_{d-1}) \neq (0,\ldots,0) } } |K+L|^{-2s} \, L_j^2 ; 
\end{equation*}
of course the sums vanish for $d=1$, while for $d=2,3$ and for any $j \in \{1,\ldots,d\}$
\begin{align}
\sum_{ \substack{  L=(L_1,\ldots,L_d) \in \mathscr{L} \\ |K+L| \neq 0 \\ L_1=0, (L_2,\ldots,L_d) \neq (0,\ldots,0) } } |K+L|^{-2s} \; L_j^2 &\leq |K|^{-2s} \sum_{ \substack{  L=(0,L_2,\ldots,L_d) \in \mathscr{L} \\ |K+L| \neq 0 \\ (L_2,\ldots,L_d) \neq (0,\ldots,0) }  } \left| \frac{ 2|L| }{ h } - 1 \right|^{-2s}  L_j^2 , \label{Est2HighFreqTerm2}
\end{align}
which is polynomially small with respect to $h$. 

On the other hand, for $\kappa = (h K_1,\ldots,h K_d)$ and $|K_1|+\cdots+|K_d| \leq M$ and for any $s > \frac{d}{2} +1$ we have
\begin{align}
&\left| \frac{\cE_\kappa}{h^{2\alpha}} - \frac{\xi_K \; \eta_K}{2} \right| \nonumber \\
&\leq \left| \omega_K^2-1 \right| \; |\hat{q}_K|^2 + \frac{1}{2} \sum_{ \substack{  L=(L_1,\ldots,L_d) \in \mathbb{Z}^d\setminus \{0\} \\ h L_1,\ldots, h L_d \in 2 \mathbb{Z}   } } \left( |\hat{p}_{K+L}|^2 + \omega_K^2 |\hat{q}_{K+L}|^2  \right) \nonumber \\
&\leq h^2 \; \pi^2 |K|^2\; |\hat{q}_K|^2 + \frac{1}{2} \sum_{ \substack{  L=(L_1,\ldots,L_d) \in \mathbb{Z}^d\setminus \{0\} \\ h L_1,\ldots, h L_d \in 2 \mathbb{Z}   } } \left( |\hat{p}_{K+L}|^2 +  |\hat{q}_{K+L}|^2 + h^2 \; \pi^2 \; |K+L|^2  \; |\hat{q}_{K+L}|^2 \right) \nonumber \\
&\leq h^2 \; \pi^2 |K|^2\; |\hat{q}_K|^2 + \| (\xi,\eta) \|_{\cW^s}^2 \sum_{ \substack{  L=(L_1,\ldots,L_d) \in \mathbb{Z}^d\setminus \{0\} \\ h L_1,\ldots, h L_d \in 2 \mathbb{Z}   } } |K+L|^{-2s} \left( 1 + h^2 \; \pi^2 \; |K+L|^2   \right) \nonumber \\
&\leq \pi^2 \; h^{2 \delta } \; \| (\xi,\eta) \|_{\cW^0}^2 + \| (\xi,\eta) \|_{\cW^s}^2 \sum_{ \substack{  L=(L_1,\ldots,L_d) \in \mathbb{Z}^d\setminus \{0\} \\ h L_1,\ldots, h L_d \in 2 \mathbb{Z}   } } |K+L|^{-2s} \left( 1 + h^2 \; \pi^2 \; |K+L|^2   \right) , \label{EstLowFreqTerm1} 
\end{align}
while for any $s > \frac{d}{2}+1$ we estimate the sum in  \eqref{EstLowFreqTerm1} in the following way,
\begin{align}
& \sum_{ \substack{  L=(L_1,\ldots,L_d) \in \mathbb{Z}^d\setminus \{0\} \\ h L_1,\ldots, h L_d \in 2 \mathbb{Z}   } } |K+L|^{-2s} \left( 1 + h^2 \; \pi^2 \; |K+L|^2   \right) \nonumber \\
&\leq \sum_{ \substack{  L=(L_1,\ldots,L_d) \in \mathbb{Z}^d\setminus \{0\} \\ h L_1,\ldots, h L_d \in 2 \mathbb{Z}   } } |K+L|^{-2s} \left( 1 + h^2 \; 2\pi^2 \; \sum_{j=1}^d (K_j^2 + L_j^2)   \right)  \nonumber \\
&\leq |K|^{-2s} \; 2\pi^{d-1} \; I(h,|K|,s,d) , \nonumber \\
& I(h,|K|,s,d)  \coloneqq \int_{2/h}^{+\infty} \left( \frac{R}{|K|} -1 \right)^{-2s} \left(1+2\pi^2 \; h^2 (|K|^2+R^2) \right) \; R^{d-1} \; \mathrm{d}R , \label{EstLowFreqTerm2} 
\end{align}
where
\begin{align*}
I(h,|K|,s,1) &= \left| 1- \frac{2}{|K|h} \right|^{-2s} \; \left( \frac{2}{h} - |K| \right) \; \frac{1}{(s-1)(2s-1)(2s-3)} \times \nonumber \\
&\; \times \bigg[ 3 - 5 s + 2 \left( s^2 + \pi^2 \left( 4+4s(2s-3)+2|K|h - 4 s h |K| + (4+s(2s-5)) |K|^2 h^2  \right) \right)  \bigg] ,
\end{align*}

\begin{align*}
 I(h,|K|,s,2) &= \left| 1- \frac{2}{|K|h} \right|^{-2s} \; \left( \frac{2}{h} - |K| \right) \; \frac{1}{2(s-1)(s-2)(2s-1)(2s-3) h} \times \nonumber \\
&\; \; \times \bigg[  (s-2)(2s-3) (4s-2-|K| h)  \nonumber \\
&\;\;\;\; + 2 \pi^2 \bigg(  8(s-1)(2s-1)(2s-3) - 12(s-1)(2s-1) \; |K|h \nonumber \\
&\; \; \; \;  + 2(2s-1)(9+s(2s-7)) \; |K|^2 h^2 - (9+ s(2s-7)) \; |K|^3 h^3  \bigg) \bigg] , \nonumber \\
I(h,|K|,s,3) &= \left| 1- \frac{2}{|K|h} \right|^{-2s} \;  \frac{1}{(s-1)(s-2)(2s-1)(2s-3)(2s-5) \; h^3} \times \nonumber \\
&\; \; \; \; \times \bigg[ 2\pi^2 \bigg(  32(s-1)(s-2)(2s-1)(2s-3) -16 s(s-1)(2s-1)(2s-3) \; |K|h \nonumber \\
&\; \; \; \; \; \; +8 (s-1)(2s-1) (10+s(2s-5)) \; |K|^2h^2 - 4 s(2s-1)(16+s(2s-9)) \; |K|^3h^3 \nonumber \\
&\; \; \; \; \; \; + 4s (16 +s(2s-9)) \; |K|^4 h^4 +(-16+s(2s-9)) \; |K|^5 h^5  \bigg) \nonumber \\
&\; \; \; \; \; \; + (s-2)(2s-5)(2-|K|h) \left( 4+8s^2 + (2+|K|h) \; |K|h - 4s (3+|K|h)  \right) \bigg] ,
\end{align*}

and since 
\begin{align*}
I(h,|K|,s,d) &\leq C_2(s) \; h^{2s-d} \; |K|^{2s} , \; \; d=1,2,3 
\end{align*}
we obtain
\begin{align}
 \sum_{ \substack{  L=(L_1,\ldots,L_d) \in \mathbb{Z}^d\setminus \{0\} \\ h L_1,\ldots, h L_d \in 2 \mathbb{Z}   } } |K+L|^{-2s} \left( 1 + h^2 \; \pi^2 \; |K+L|^2   \right) &\leq 2\pi^{d-1} \; C_2(s) \; h^{2s-d}. \label{EstLowFreqTerm21} 
\end{align}

Finally, we observe that
\begin{align*}
\min \left( 2\delta, 2s-d \right) &= 2\delta , \; \; \forall s > \frac{d}{2} +\delta ,
\end{align*}

so that we can deduce the thesis with the constants specified in \eqref{constAppr}.

\end{appendix}

\bibliography{P_KG_2021}
\bibliographystyle{alpha}

\end{document}